\pdfoutput=1
\documentclass[11pt,reqno]{amsart}
\usepackage[letterpaper,hmargin=1in,vmargin=1.1in]{geometry}
\usepackage{amsmath,amssymb,amsthm,mathtools,enumerate,booktabs,microtype,hyperref}
\hypersetup{colorlinks=true, linkcolor=blue, citecolor=red, urlcolor=teal}

\theoremstyle{plain}
\newtheorem{theorem}{Theorem}[section]
\newtheorem{lemma}[theorem]{Lemma}
\newtheorem{proposition}[theorem]{Proposition}
\newtheorem{corollary}[theorem]{Corollary}

\theoremstyle{definition}
\newtheorem{definition}[theorem]{Definition}

\theoremstyle{remark}
\newtheorem{remark}[theorem]{Remark}

\numberwithin{equation}{section}

\title[Hamiltonian Cotangent Deformations of Space Forms]{Hamiltonian Cotangent Deformations of Space Forms: Curvature Pinching, Soliton Rigidity, Conjugate Loci, and Superintegrability}
\author{James J. Munden}
\address{Dean of Sciences, Saint Louis Community College -- Meramec, St. Louis, MO 63122, USA}
\email{jmunden1@stlcc.edu}

\date{\today}
\dedicatory{Dedicated to the memory of James J.\ Hebda}

\begin{document}
\begin{abstract}
We study the one-parameter family of Riemannian metrics $g_\lambda$ arising
from the Hamiltonian cotangent deformation
$H_\lambda = H_0 + \frac{\lambda}{2}\|J\|^2$ on simply connected space forms
of constant curvature $\kappa\in\{1,0,-1\}$, where $J$ is the momentum map
of the isotropy action fixing a base point. A single algebraic identity---the
master reciprocal relation $1/f_\lambda^2 = 1/s_\kappa^2 + \lambda$ governing
the warping function---controls the geometry of the entire family and implies
that the radial geodesic acceleration is independent of the deformation
parameter.

From this mechanism we derive three families of consequences.
\emph{Curvature and rigidity}: exact closed-form sectional curvatures, the
sharp positivity trichotomy across signatures, the pointwise and global Brendle--Schoen
$\frac{1}{4}$-pinching windows on $\mathbb{S}^{n+1}$, never-Einstein rigidity
for $\lambda\neq 0$, and the identification of Hamilton's steady cigar soliton
as the unique gradient Ricci soliton in the family.
\emph{Dynamics and variational stability}: iso-kinematic geodesic decoupling
into unperturbed radial motion and uniform angular precession, closed-form
Clairaut trajectory integration, the transcendental inward conjugate
point equation with Euclidean spherical inversion, and stratified $A_3$
cuspidal caustic envelopes.
\emph{Integrability and spectral theory}: maximal superintegrability on
$\mathbb{R}^{n+1}$ via rational pericenter invariants, transverse Killing
rigidity saturating the BEHR coalgebra bound on curved space forms, Fuchsian
Heun spectral reduction with cross-ratio modulus $\tau = 1 + \lambda/\kappa$,
and essential spectrum gap collapse on hyperbolic space.

Extensions to anisotropic (direction-dependent) multi-parameter families and
intrinsic deformations of compact rank-one symmetric spaces
(CROSSes: $\mathbb{C}P^n$, $\mathbb{H}P^n$, $\mathbb{O}P^2$) are developed
through the same Lie-theoretic mechanism applied to transitive subgroup
actions on spheres.
\end{abstract}

\maketitle

\tableofcontents

\section{Introduction}\label{sec:introduction}

\subsection{From rigid space forms to deformed Hamiltonians}\label{subsec:motivation}
The simply connected Riemannian space forms---the unit sphere $\mathbb{S}^{n+1}$, Euclidean space $\mathbb{R}^{n+1}$, and hyperbolic space $\mathbb{H}^{n+1}$ of constant sectional curvature $\kappa \in \{1, 0, -1\}$ with dimension $N = n + 1 \ge 2$---have isometry groups of maximal dimension $(n+1)(n+2)/2$, maximally superintegrable geodesic flows, and isotropic curvature tensors. Constructing smooth metric deformations on space forms that preserve complete geodesic integrability is a classical problem in differential geometry, since generic metric perturbations break invariant Liouville tori via the Poincar\'e--Birkhoff mechanism.

A general framework for integrable metric deformations was introduced by Hebda \cite{Hebda2001} using cotangent momentum maps (building on earlier curvature pinching examples \cite{Hebda1985}). Given a complete Riemannian manifold $(M, g_0)$ admitting an isometric Lie group action by $G$, one deforms the kinetic energy Hamiltonian $H_0(\alpha) = \frac{1}{2}\|\alpha\|_{g_0^*}^2$ on $T^*M$ by coupling $H_0$ to an $\mathrm{Ad}^*(G)$-invariant quadratic function of the momentum map $J: T^*M \to \mathfrak{g}^*$. On a simply connected space form $(Q^{n+1}, g_0)$, fixing an origin $p_0 \in Q$ isolates the compact isotropy group $G = \mathrm{SO}(n+1)$, which acts transitively and isometrically on each concentric geodesic distance sphere $\mathcal{O}_r = \partial B(p_0, r) \cong \mathbb{S}^n$. The cotangent momentum map $J: T^*Q \to \mathfrak{so}(n+1)^*$ maps covectors to angular momentum bivectors relative to $p_0$.

In the author's doctoral dissertation \cite{Munden2015}, cotangent Hamiltonian deformations of simply connected space forms were classified according to the transitive isometric subgroup actions on the geodesic spheres $\mathcal{O}_r \cong \mathbb{S}^d$ (Montgomery--Samelson \cite{MontgomerySamelson1943}; Borel \cite{Borel1949}; Besse \cite{Besse1987}), yielding five geometric families:
\begin{enumerate}[(I)]
\item \emph{Family~I (Isotropic)}: Governed by the full transitive action of $\mathrm{SO}(n+1)$ on $\mathbb{S}^n \cong \mathrm{SO}(n+1)/\mathrm{SO}(n)$, scaling the transverse geodesic spheres isotropically via a single scalar parameter $\lambda \in \mathbb{R}$.
\item \emph{Family~II (Complex Hopf)}: Arising from the action of $\mathrm{U}(m+1)$ on $\mathbb{S}^{2m+1} \cong \mathrm{U}(m+1)/\mathrm{U}(m)$, fibering via the complex Hopf bundle $\mathbb{S}^1 \hookrightarrow \mathbb{S}^{2m+1} \to \mathbb{C}P^m$ and deforming independently along the $1$-dimensional vertical Reeb fiber and $2m$-dimensional horizontal distribution.
\item \emph{Family~III (Octonionic / Cayley Hopf)}: Governed by the spin action of $\mathrm{Spin}(9)$ on $\mathbb{S}^{15} \cong \mathrm{Spin}(9)/\mathrm{Spin}(7)$, fibering via the octonionic Hopf bundle $\mathbb{S}^7 \hookrightarrow \mathbb{S}^{15} \to \mathbb{S}^8 \cong \mathbb{O}P^1$ and splitting $T\mathbb{S}^{15}$ into the $7$-dimensional vertical fiber and $8$-dimensional horizontal spinor distribution.
\item \emph{Family~IV (Quaternionic Hopf)}: Governed by the action of $\mathrm{Sp}(m+1)\mathrm{Sp}(1)$ on $\mathbb{S}^{4m+3} \cong (\mathrm{Sp}(m+1)\times\mathrm{Sp}(1))/(\mathrm{Sp}(m)\times\mathrm{Sp}(1))$ along the quaternionic Hopf fibration $\mathbb{S}^3 \hookrightarrow \mathbb{S}^{4m+3} \to \mathbb{H}P^m$, deforming along the $3$-dimensional $\mathfrak{sp}(1)$ vertical fiber and $4m$-dimensional horizontal subbundle.
\item \emph{Family~V (Intermediate Quaternionic--Complex)}: Governed by the restricted action of $\mathrm{Sp}(m+1)\mathrm{U}(1)$ on $\mathbb{S}^{4m+3} \cong (\mathrm{Sp}(m+1)\times\mathrm{U}(1))/(\mathrm{Sp}(m)\times\mathrm{U}(1))$ along the two-stage fibration $\mathbb{S}^1 \hookrightarrow \mathbb{S}^{4m+3} \to \mathbb{C}P^{2m+1} \to \mathbb{H}P^m$, splitting $T\mathbb{S}^{4m+3}$ into three mutually orthogonal $\mathrm{Ad}(H)$-irreducible subbundles of dimensions $1$, $4m$, and $2$.
\end{enumerate}

The main body of this paper is devoted to the isotropic (direction-independent) Family~I across all space form background curvatures $\kappa \in \{1, 0, -1\}$ and all dimensions $N = n+1 \ge 2$. We determine its curvature tensor and pinching bounds, Einstein and gradient Ricci soliton rigidity, geodesic flow, conjugate and cut loci, superintegrability, and Laplace--Beltrami spectrum, before extending the analysis in Sections~\ref{sec:anisotropic_families} and~\ref{sec:cross_deformations} to the anisotropic Families~II--V and to intrinsic cotangent deformations of compact rank-one symmetric spaces (CROSSes: $\mathbb{C}P^n$, $\mathbb{H}P^n$, and the Cayley plane $\mathbb{O}P^2$).

\subsection{The Lie-theoretic mechanism}\label{subsec:lie_mechanism}
The complete integrability of Family~I geodesic flows follows from a Lie-theoretic mechanism on the cotangent bundle $(T^*Q, \omega_{\mathrm{can}})$, where $\omega_{\mathrm{can}} = d\theta_{\mathrm{can}} = dp_i \wedge dq^i$ is the canonical symplectic form. In terms of the unperturbed kinetic energy $H_0(\alpha) = \frac{1}{2}\|\alpha\|_{g_0^*}^2$ and the momentum map $J: T^*Q \to \mathfrak{so}(n+1)^*$ of the isotropy action fixing $p_0$, the deformed Hamiltonian is
\begin{equation}\label{eq:intro_hamiltonian_def}
H_\lambda(\alpha) = H_0(\alpha) + \frac{\lambda}{2} \|J(\alpha)\|^2,
\end{equation}
where $\lambda \in \mathbb{R}$ is the deformation parameter and $\|\cdot\|$ is the norm on $\mathfrak{so}(n+1)^*$ dual to the $\mathrm{Ad}(\mathrm{SO}(n+1))$-invariant trace inner product $\langle A, B \rangle = -\frac{1}{2}\operatorname{Tr}(AB)$ on $\mathfrak{so}(n+1)$ (in which the elementary rotation generators $E_{ab} - E_{ba}$ for $1 \le a < b \le n+1$ form an orthonormal basis for all $n \ge 1$).

In geodesic polar coordinates $(r, \theta) \in I_\kappa \times \mathbb{S}^n$ centered at $p_0$, the space form metric is $g_0 = dr^2 + s_\kappa^2(r) g_{\mathbb{S}^n}$, where $s_1(r) = \sin r$, $s_0(r) = r$, and $s_{-1}(r) = \sinh r$. The momentum map norm evaluates on covectors $\alpha = p_r dr + p_\theta d\theta$ to the angular momentum magnitude:
\begin{equation}\label{eq:intro_casimir_evaluation}
\|J(\alpha)\|^2 = \|p_{\mathbb{S}^n}\|_{g_{\mathbb{S}^n}^*}^2 = \sum_{1 \le a < b \le n+1} J_{ab}(\alpha)^2.
\end{equation}
This expression evaluates the quadratic Casimir operator $C_2(\mathfrak{so}(n+1)) = \sum_{a < b} \xi_{ab}^2$ in the universal enveloping algebra $\mathcal{U}(\mathfrak{so}(n+1))$ on coadjoint orbits. Thus, the deformed Hamiltonian decomposes in phase space as
\begin{equation}\label{eq:intro_hamiltonian_polar}
H_\lambda(r, \theta, p_r, p_\theta) = \frac{1}{2} p_r^2 + \frac{1}{2}\left( \frac{1}{s_\kappa^2(r)} + \lambda \right) \|p_{\mathbb{S}^n}\|_{g_{\mathbb{S}^n}^*}^2.
\end{equation}

The deformation parameter $\lambda \ne 0$ breaks the ambient transitive isometry group:
\begin{equation}\label{eq:intro_symmetry_breaking}
\mathrm{SO}(n+2) \longrightarrow \mathrm{SO}(n+1) \times \mathbb{Z}_2 \quad (\text{on } \mathbb{S}^{n+1}),
\end{equation}
breaking the $(n+1)$ transvection generators shifting $p_0$ while preserving the spatial isotropy group $\mathrm{SO}(n+1)$ and the antipodal reflection $r \mapsto \pi - r$. Analogously, $\mathrm{E}(n+1)$ and $\mathrm{SO}^+(n+1, 1)$ reduce to $\mathrm{O}(n+1)$ on $\mathbb{R}^{n+1}$ and $\mathbb{H}^{n+1}$.

Nevertheless, the Casimir sector remains conserved. Because $g_0$ is $\mathrm{SO}(n+1)$-invariant, each momentum component $J_\xi$ Poisson-commutes with $H_0$: $\{H_0, J_\xi\} = 0$. By the Leibniz rule,
\begin{equation}\label{eq:intro_poisson_casimir}
\left\{ H_0, \|J\|^2 \right\} = 2 \sum_{a < b} J_{ab} \left\{ H_0, J_{ab} \right\} = 0.
\end{equation}
Since $[X_F, X_K] = -X_{\{F, K\}}$, the Hamiltonian vector fields commute: $[X_{H_0}, X_{H_\lambda}] = 0$. The Hamiltonian flow factorizes into commuting flows: $\Phi_t^{H_\lambda} = \Phi_t^{H_0} \circ \Phi_t^{\frac{\lambda}{2}\|J\|^2}$.

This Lie-theoretic factorization implies that angular modes decouple from radial motion. The flow of $\frac{\lambda}{2}\|J\|^2$ generates rotations along $\mathrm{SO}(n+1)$ orbits at constant rate $\Omega = \lambda L$, imparting zero net radial force. For the Laplace--Beltrami operator, the angular Laplacian on $\mathbb{S}^n$ acts on degree-$l$ spherical harmonics via Peter--Weyl theory with the Casimir eigenvalue $-l(l+n-1)$. By the master reciprocal identity $\frac{1}{f_\lambda(r)^2} = \frac{1}{s_\kappa(r)^2} + \lambda$, the deformation produces a constant shift $\lambda l(l+n-1)$ in the centrifugal potential, preserving separation of variables.

\subsection{Statement of principal results}\label{subsec:principal_results}
We state five principal theorems summarizing the geometric, dynamical, variational, and spectral properties of the isotropic cotangent deformations.

\begin{theorem}[Sectional Curvature Trichotomy, Sphere Pinching, and Brendle--Schoen Windows]\label{thm:intro_curvature_pinching}
Let $(Q^{n+1}, g_0)$ be an $(n+1)$-dimensional complete simply connected space form of constant curvature $\kappa \in \{1, 0, -1\}$ ($N = n + 1 \ge 2$), and let $g_\lambda = dr^2 + f_\lambda(r)^2 g_{\mathbb{S}^n}$ be the Family~I metric with warping factor $f_\lambda(r) = \frac{s_\kappa(r)}{\sqrt{1 + \lambda s_\kappa^2(r)}}$, where $\lambda \in \mathbb{R}$ satisfies $1 + \lambda s_\kappa^2(r) > 0$.
\begin{enumerate}[(i)]
\item \textbf{Positivity Trichotomy}: The sectional curvature satisfies $\sec(g_\lambda) > 0$ on all tangent two-planes throughout $Q^{n+1}$ if and only if:
\begin{enumerate}[(a)]
\item On Euclidean space $\mathbb{R}^{n+1}$ ($\kappa = 0$): $\lambda > 0$, with bounds $0 < \sec(g_\lambda) \le 3\lambda$ and $\inf \sec = \lim_{r \to \infty} K_{\mathrm{rad}}(r) = 0$;
\item On hyperbolic space $\mathbb{H}^{n+1}$ ($\kappa = -1$): $\lambda > 1/3$, with $\inf \sec(g_\lambda) = 0$ and central pole curvature $\sec(p_0) = 3\lambda - 1 > 0$;
\item On the sphere $\mathbb{S}^{n+1}$ ($\kappa = 1$): $\lambda > -1/3$.
\end{enumerate}
\item \textbf{Global and Pointwise Sphere Pinching Ratios}: On $\mathbb{S}^{n+1}$ ($\kappa = 1$), for every $\lambda > -1/3$ (with $n \ge 2$ for tangential 2-planes), the global (manifold-wide) sectional curvature pinching ratio $\delta(\lambda) := \frac{\min_{p, \Pi} \sec(\Pi)}{\max_{p, \Pi} \sec(\Pi)}$ is given by:
\begin{equation}\label{eq:intro_pinching_formula}
\delta(\lambda) = \begin{cases} \dfrac{1}{(1+\lambda)(1+3\lambda)}, & \lambda \ge 0, \\[10pt] (1+\lambda)(1+3\lambda), & -\dfrac{1}{3} < \lambda < 0. \end{cases}
\end{equation}
For $\lambda \ge 0$, $\min_{p,\Pi} \sec = K_{\mathrm{rad}}(\pi/2) = \frac{1}{1+\lambda}$ (at the equator $r = \pi/2$) and $\max_{p,\Pi} \sec = K_{\mathrm{rad}}(0) = K_{\mathrm{tan}}(0) = 1 + 3\lambda$ (at the poles $r \in \{0, \pi\}$); for $\lambda \in (-1/3, 0)$, the global extrema invert. In contrast, the pointwise pinching ratio $\delta_p(g_\lambda) := \frac{\min_{\Pi \subset T_p\mathbb{S}^{n+1}} \sec(\Pi)}{\max_{\Pi \subset T_p\mathbb{S}^{n+1}} \sec(\Pi)}$ ($n \ge 2$) attains its minimum over $p \in \mathbb{S}^{n+1}$ at the equator $r = \pi/2$, where $\delta_{\mathrm{eq}}(\lambda) = \frac{K_{\mathrm{rad}}(\pi/2)}{K_{\mathrm{tan}}(\pi/2)} = (1+\lambda)^{-2}$ for $\lambda \ge 0$ and $\delta_{\mathrm{eq}}(\lambda) = (1+\lambda)^2$ for $\lambda \in (-1/3, 0)$.
\item \textbf{Global and Pointwise Brendle--Schoen $1/4$-Pinching Windows}: The metric $g_\lambda$ on $\mathbb{S}^{n+1}$ is globally strictly $1/4$-pinched ($\delta(\lambda) > 1/4$) if and only if $\lambda$ lies in the interval
\begin{equation}\label{eq:intro_brendle_schoen}
\lambda \in \left( \frac{-4 + \sqrt{7}}{6}, \, \frac{\sqrt{13} - 2}{3} \right) \approx (-0.2257, \, 0.5352),
\end{equation}
whereas $g_\lambda$ is pointwise strictly $1/4$-pinched ($\delta_p(g_\lambda) > 1/4$ for all $p \in \mathbb{S}^{n+1}$) on the strictly wider interval $\lambda \in (-1/3, 1)$ (where positive sectional curvature fails at $\lambda = -1/3$ and $\delta_{\mathrm{eq}}(1) = 1/4$). By the Brendle--Schoen Differentiable Sphere Theorem \cite{BrendleSchoen2009}, $(\mathbb{S}^{n+1}, g_\lambda)$ is diffeomorphic to the standard sphere $\mathbb{S}^{n+1}$ throughout these pinching windows.
\end{enumerate}
\end{theorem}

\begin{theorem}[Never-Einstein Rigidity and Gradient Ricci Soliton Classification]\label{thm:intro_soliton_rigidity}
Let $(Q^{n+1}, g_\lambda)$ be an $(n+1)$-dimensional complete space form of Family~I with parameter $\lambda \in \mathbb{R}$ and spatial dimension $n \ge 1$.
\begin{enumerate}[(i)]
\item \textbf{Never-Einstein Rigidity}: The metric $g_\lambda$ is Einstein ($\mathrm{Ric} = \alpha g_\lambda$) if and only if $\lambda = 0$, recovering the background space form metric $g_0$.
\item \textbf{Universal Equator Defect on $\mathbb{S}^{n+1}$}: On $\mathbb{S}^{n+1}$ ($\kappa = 1$), complete gradient Ricci solitons $\mathrm{Ric} + \nabla^2 h = \rho g_\lambda$ are unconditionally obstructed for all $N = n + 1 \ge 2$ and all $\lambda \ne 0$ by the non-vanishing universal equator defect:
\begin{equation}\label{eq:intro_equator_defect}
\mathcal{D}(\pi/2) := \bar{h}''(\pi/2) + n K_{\mathrm{rad}}(\pi/2) - \rho \equiv \frac{4\lambda}{1+\lambda} \ne 0.
\end{equation}
\item \textbf{Euclidean Soliton Obstruction}: On $\mathbb{R}^{n+1}$ ($\kappa = 0$), complete gradient Ricci solitons are unconditionally obstructed for all $N \ge 2$ and all $\lambda > 0$ by the non-vanishing order-$r^4$ Taylor residual coefficient:
\begin{equation}\label{eq:intro_c4_residual}
c_4 = 2(2n - 5)\lambda^3 \ne 0 \qquad (\forall n \in \mathbb{N}_{\ge 1}),
\end{equation}
which vanishes uniquely at the non-integer dimension $n = 5/2 \notin \mathbb{N}$.
\item \textbf{Hamilton's Steady Cigar Soliton on $(\mathbb{H}^2, g_1)$}: On hyperbolic space $\mathbb{H}^2$ ($\kappa = -1, n = 1$) at $\lambda = 1$, the metric $g_1 = dr^2 + \tanh^2(r) d\theta^2$ is identically Hamilton's steady cigar soliton \cite{Hamilton1988}, satisfying $\mathrm{Ric} + \nabla^2 h = 0$ with potential $h(r) = -2\ln(\cosh r)$. For all $\lambda \ne 1$ on $\mathbb{H}^2$ and all $\lambda > 0$ on $\mathbb{H}^{n+1}$ ($n \ge 2$), complete gradient solitons are obstructed.
\end{enumerate}
\end{theorem}

\begin{theorem}[Iso-Kinematic Geodesic Decoupling and Universal Precession Law]\label{thm:intro_geodesic_precession}
Let $\gamma(t) = (r(t), \theta(t))$ be a unit-speed geodesic in $(Q^{n+1}, g_\lambda)$ with Clairaut angular momentum $L := g_\lambda(\dot{\gamma}(t), \partial_\theta) = f_\lambda(r(t))^2 \dot{\theta}(t)$.
\begin{enumerate}[(i)]
\item \textbf{Iso-Kinematic Radial Decoupling}: Across all space forms $\kappa \in \{1, 0, -1\}$, the radial coordinate acceleration satisfies identically:
\begin{equation}\label{eq:intro_radial_acceleration}
\ddot{r}(t) = L^2 \frac{c_\kappa(r(t))}{s_\kappa(r(t))^3},
\end{equation}
independent of $\lambda$. The radial trajectory $r(t)$ matches that of an unperturbed space form geodesic in $(Q, g_0)$ with angular momentum $L$ and effective dilated speed:
\begin{equation}\label{eq:intro_dilated_speed}
v_0 = \sqrt{1 - \lambda L^2}.
\end{equation}
\item \textbf{Universal Precession Law}: The angular coordinate evolves according to:
\begin{equation}\label{eq:intro_precession_law}
\dot{\theta}(t) = \dot{\theta}_0(t) + \lambda L \implies \theta(t) = \theta_0(t) + \lambda L t,
\end{equation}
where $\dot{\theta}_0(t) = L/s_\kappa(r(t))^2$. The geodesic undergoes steady gyroscopic precession around $p_0$ at the uniform rate $\Omega = \lambda L$.
\item \textbf{Clairaut Trajectory Integration \& KAM Twist}: For non-radial geodesics on $\mathbb{S}^{n+1}$ with $L \in (0, 1/\sqrt{1+\lambda})$ and turning parameter $B := \sqrt{1 - (1+\lambda)L^2}/L > 0$, the orbital trajectory in cotangent coordinate $u := \cot r$ integrates to the closed primitive
\begin{equation}\label{eq:intro_clairaut_primitive}
\Delta\theta(u) = \arcsin\left( \frac{u}{B} \right) + \left[\frac{\lambda L}{\sqrt{1 - \lambda L^2}}\right] \arctan\left( \left[\frac{\sqrt{1 - \lambda L^2}}{L}\right] \frac{u}{\sqrt{B^2 - u^2}} \right),
\end{equation}
yielding the half-period apsidal angle $\Delta\Phi(L) = \pi\left(1 + \frac{\lambda L}{\sqrt{1-\lambda L^2}}\right)$ between consecutive radial turning points (full periastron-to-periastron advance $2\Delta\Phi(L)$) and non-vanishing KAM twist:
\begin{equation}\label{eq:intro_kam_twist}
\frac{d\nu}{dL} = \frac{\lambda}{2(1 - \lambda L^2)^{3/2}} \ne 0 \qquad (\forall \lambda \ne 0).
\end{equation}
\end{enumerate}
\end{theorem}

\begin{theorem}[Conjugate/Cut Loci, Inward Conjugacy, and Transverse Bounds]\label{thm:intro_conjugate_cut}
Let $(Q^{n+1}, g_\lambda)$ be an $(n+1)$-dimensional complete space form of Family~I with spatial dimension $n \ge 1$.
\begin{enumerate}[(i)]
\item \textbf{Polar Radial Jacobi Fields and Origin Cut Locus}: Along every radial geodesic ray $\gamma(t) = (t, \theta_0)$ from $p_0$, the normal Jacobi fields are $Y(t) = f_\lambda(t) E(t) = \frac{s_\kappa(t)}{\sqrt{1+\lambda s_\kappa^2(t)}} E(t)$. On $\mathbb{R}^{n+1}$ ($\lambda \ge 0$) and $\mathbb{H}^{n+1}$ ($\lambda \ge 0$), radial rays are unique global length-minimizers on $[0, \infty)$, yielding $\mathrm{Cut}(p_0) = \emptyset$ and $\mathrm{Inj}(p_0) = \infty$. On $\mathbb{S}^{n+1}$ ($\lambda > -1/3$), the first conjugate point occurs at $t = \pi$ with multiplicity $n = \dim(\mathbb{S}^n)$, establishing $\mathrm{Cut}(p_0) = \{p_\pi\}$ and $\mathrm{Inj}(p_0) = \pi$.
\item \textbf{Universal Transcendental Inward Conjugate Equation}: For an off-origin base point $q_0 = (r_0, 0)$ ($r_0 > 0$), let $\gamma(t)$ be the inward diameter geodesic ray with $\dot{\gamma}(0) = -\partial_r$, traversing $p_0$ at $t = r_0$ and continuing along the opposite meridian. The first conjugate point occurs at the smallest root $t_{\mathrm{conj}} > r_0$ of:
\begin{equation}\label{eq:intro_transcendental_conj}
\frac{s_\kappa(t_{\mathrm{conj}})}{s_\kappa(t_{\mathrm{conj}} - r_0)} = \lambda s_\kappa(r_0) t_{\mathrm{conj}},
\end{equation}
possessing conjugate multiplicity exactly $n = \dim(\mathbb{S}^n)$.
\item \textbf{Euclidean Spherical Inversion}: On Euclidean space $\mathbb{R}^{n+1}$ ($\kappa = 0, \lambda > 0$), equation \eqref{eq:intro_transcendental_conj} reduces to $t_{\mathrm{conj}} = r_0 + \frac{1}{\lambda r_0}$, establishing the exact projective spherical inversion relation:
\begin{equation}\label{eq:intro_spherical_inversion}
r_0 \cdot R^* = \frac{1}{\lambda}, \qquad R^* := t_{\mathrm{conj}} - r_0.
\end{equation}
\item \textbf{Equatorial Mode Preemption and Oblique Conjugate Delay}: On $\mathbb{S}^{n+1}$, equatorial geodesics exhibit mode preemption $t_{\mathrm{first}}(\pi/2) = \min(\pi/\sqrt{1+\lambda}, \pi\sqrt{1+\lambda}) < \pi$ for all $\lambda \ne 0$. In negative pinching $\lambda \in (-1/3, 0)$, steep oblique rays linger in polar caps of reduced curvature $1+3\lambda < 1$, inducing an inclination bifurcation with delayed conjugacy $t_{\mathrm{first}} > \pi$.
\item \textbf{Transverse Encircling Bound on $\mathbb{R}^{n+1}$}: For any base point $q_0 = (r_0, 0)$ ($r_0 > 0$) on $(\mathbb{R}^{n+1}, g_\lambda)$ ($\lambda > 0$), the cut distance satisfies the dual bound:
\begin{equation}\label{eq:intro_cut_bound}
\mathrm{dist}(q_0, \mathrm{Cut}(q_0)) \le \min\left( r_0 + \frac{1}{\lambda r_0}, \, \frac{\pi}{\sqrt{\lambda}} \right),
\end{equation}
governed by radial conjugate preemption along $L = 0$ and symmetric transverse encircling collision along $L > 0$.
\end{enumerate}
\end{theorem}

\begin{theorem}[Superintegrability Classification and Confluent Heun Spectral Bridge]\label{thm:intro_superintegrability_spectral}
Let $(Q^{n+1}, g_\lambda)$ be an $(n+1)$-dimensional complete space form of Family~I with spatial dimension $n \ge 1$ and manifold dimension $N = n + 1 \ge 2$.
\begin{enumerate}[(i)]
\item \textbf{Liouville and Noncommutative Integrability}: Across all space forms and all dimensions $N \ge 2$, the geodesic flow is completely Liouville integrable via the Gelfand--Tsetlin Casimir chain $\mathcal{F}_{\mathrm{Liouv}} = \{H_\lambda, C_2, \dots, C_{n+1}\}$. For $N \ge 3$, the augmented momentum map has maximal rank $2N - 2$, proving quasi-maximal superintegrability of degree $\ge 2N - 2$ via the Mishchenko--Fomenko theorem \cite{MishchenkoFomenko1978}.
\item \textbf{Maximal Superintegrability on $\mathbb{R}^{n+1}$ via Pericenter Invariants}: On deformed Euclidean space $(\mathbb{R}^{n+1}, g_\lambda)$ ($\kappa = 0, \lambda > 0$), the geodesic flow is strictly maximally superintegrable of degree $2N - 1$. The rational pericenter transit time $\tau(x, p) := \frac{\langle x, p \rangle}{\|p\|^2}$ satisfies $\{\tau, H_\lambda\} = 1$, and the vector quantities
\begin{equation}\label{eq:intro_pericenter_invariants}
p_{\mathrm{peri}} := \exp(-\lambda \tau(x, p) \hat{J}) p, \qquad x_{\mathrm{peri}} := \exp(-\lambda \tau(x, p) \hat{J}) x - \tau(x, p) p_{\mathrm{peri}}
\end{equation}
are globally smooth, single-valued rational Runge--Lenz type first integrals on $T^*\mathbb{R}^{n+1} \setminus \{p=0\}$ satisfying $\{H_\lambda, p_{\mathrm{peri}}\} = 0$ and $\{H_\lambda, x_{\mathrm{peri}}\} = 0$.
\item \textbf{Centrifugal Casimir Decomposition and Fuchsian Heun Mapping}: Under spherical harmonic separation $u(r, \theta) = R(r) Y_l(\theta)$, the eigenvalue equation $-\Delta_{g_\lambda} u = E u$ separates into the radial Sturm--Liouville equation:
\begin{equation}\label{eq:intro_radial_laplacian}
-\frac{1}{f_\lambda(r)^n} \frac{d}{dr}\left( f_\lambda(r)^n \frac{dR}{dr} \right) + \left( \frac{l(l+n-1)}{s_\kappa^2(r)} + \lambda l(l+n-1) \right) R(r) = E R(r),
\end{equation}
where the centrifugal barrier shifts by the exact constant Casimir eigenvalue $\lambda l(l+n-1)$. Under the algebraic coordinate transformation $z = s_\kappa^2(r)$, equation \eqref{eq:intro_radial_laplacian} transforms into the Fuchsian Heun differential equation
\begin{equation}\label{eq:intro_heun_ode}
z(1-\kappa z)(1+\lambda z) \frac{d^2 R}{dz^2} + P(z) \frac{dR}{dz} + Q(z) R = 0,
\end{equation}
possessing four regular singular points $\{0, 1/\kappa, -1/\lambda, \infty\}$ with cross-ratio modulus $\tau(\kappa, \lambda) = 1 + \lambda/\kappa$.
\item \textbf{Rayleigh--Schrödinger Perturbation and Degeneracy Splitting}: On $\mathbb{S}^{n+1}$, the unperturbed round sphere accidental $N = 2$ eigenvalue degeneracy ($E_{0, 2}^{(0)} = E_{1, 0}^{(0)} = 2(n+2)$) is strictly lifted for all $\lambda > 0$ with positive rational splitting:
\begin{equation}\label{eq:intro_degeneracy_splitting}
\Delta E_{N=2}(\lambda) = E_{0, 2}(\lambda) - E_{1, 0}(\lambda) = \lambda \frac{2(n+1)(n+3)(n+4)}{(n+2)(n+6)} + O(\lambda^2) > 0.
\end{equation}
\end{enumerate}
\end{theorem}

\subsection{Confluent Heun geometry and spectral bridge}\label{subsec:heun_spectral}
The transformation of the Laplace--Beltrami operator into the Fuchsian Heun differential equation \eqref{eq:intro_heun_ode} establishes an analytical bridge between Hamiltonian cotangent deformations and Fuchsian differential equations. The four regular singular points in the complex $z$-plane reflect the geometric features of the deformed space form: $z_1 = 0$ corresponds to the polar base point $p_0$; $z_2 = 1/\kappa$ corresponds to the antipodal pole $p_\pi$ on $\mathbb{S}^{n+1}$ ($\kappa = 1$), complex infinity on $\mathbb{R}^{n+1}$ ($\kappa = 0$), and the spatial boundary at infinity on $\mathbb{H}^{n+1}$ ($\kappa = -1$); $z_3 = -1/\lambda$ represents the metric neck singularity; and $z_4 = \infty$ represents the projective compactification of the complexified radial coordinate.

The invariant cross-ratio modulus of these four singular points,
\begin{equation}\label{eq:intro_cross_ratio}
\tau(\kappa, \lambda) = [z_2, z_1; z_3, z_4] = \frac{z_2 - z_3}{z_1 - z_3} \cdot \frac{z_1 - z_4}{z_2 - z_4} = \frac{1/\kappa - (-1/\lambda)}{0 - (-1/\lambda)} = 1 + \frac{\lambda}{\kappa},
\end{equation}
parameterizes the moduli space $\mathcal{M}_{0, 4}$ of four-punctured Riemann spheres, classifying the confluent limits of the deformation:
\begin{enumerate}[(a)]
\item \emph{Space Form Confluence ($\lambda \to 0$)}: As the deformation vanishes, the neck singularity $z_3 = -1/\lambda$ recedes to $-\infty$ and coalesces with $z_4 = \infty$. The singular set reduces to $\{0, 1/\kappa, \infty\}$, collapsing generic Heun functions into Gauss hypergeometric functions and spherical harmonics.
\item \emph{Flat Euclidean Confluence ($\kappa \to 0$)}: As ambient curvature vanishes, the background singular point $z_2 = 1/\kappa$ escapes to $\infty$, merging with $z_4 = \infty$ to produce an irregular singular point of rank~1. This transforms the Fuchsian Heun ODE into the confluent Heun differential equation $\mathrm{HeunC}$, governing spheroidal wave functions.
\item \emph{Double Confluence ($\lambda \to 0, \kappa \to 0$)}: Under simultaneous vanishing of background curvature and deformation, the equation reduces to the classical Bessel differential equation governing Euclidean polar modes.
\end{enumerate}

Beyond discrete eigenvalue spectra, this Heun framework reveals a spectral phase transition on hyperbolic space $(\mathbb{H}^{n+1}, g_\lambda)$. On standard hyperbolic space $(\mathbb{H}^{n+1}, g_0)$, the spectrum of $-\Delta_{g_0}$ is purely continuous with a positive gap: $\sigma(-\Delta_{g_0}) = [n^2/4, \infty)$. Under the isotropic deformation with $\lambda > 0$, the warping function satisfies $\lim_{r \to \infty} f_\lambda(r) = 1/\sqrt{\lambda} < \infty$. The deformed hyperbolic manifold is asymptotically cylindrical:
\begin{equation}\label{eq:intro_cylindrical_end}
(\mathbb{H}^{n+1}, g_\lambda) \approx \mathbb{R} \times \mathbb{S}^n(1/\sqrt{\lambda}) \quad \text{as } r \to \infty,
\end{equation}
with finite cross-sectional volume $\lambda^{-n/2}\mathrm{Vol}(\mathbb{S}^n)$. The volume density $w_\lambda(r) = \sinh^n(r)/(1+\lambda\sinh^2 r)^{n/2}$ converges exponentially to $\lambda^{-n/2}$, causing the logarithmic volume drift $\eta_\lambda(r) := \partial_r \ln w_\lambda(r) = w_\lambda'(r)/w_\lambda(r)$ to decay as $O(e^{-2r})$.

Constructing an explicit Weyl sequence of cutoff plane waves $\psi_j(r) = \frac{1}{\sqrt{j}}\chi\left(\frac{r-j}{j}\right)e^{ikr} Y_0(\theta)$ in the cylindrical end shows that for every $k \in \mathbb{R}$, the energy $E = k^2$ satisfies $\|(-\Delta_{g_\lambda} - k^2)\psi_j\|_{L^2} \to 0$ as $j \to \infty$. By Weyl's criterion, this yields essential spectrum gap collapse:
\begin{equation}\label{eq:intro_gap_collapse}
\sigma_{\mathrm{ess}}(-\Delta_{g_\lambda}) = [0, \infty) \qquad (\forall \lambda > 0).
\end{equation}
For $\lambda > 1/3$, this gap collapse occurs on a complete manifold of strictly positive sectional curvature ($0 < \sec(g_\lambda) \le 3\lambda - 1$).

\subsection*{Relation to prior work}
We briefly compare Hamiltonian cotangent deformations with four related classes of metrics in Riemannian geometry and geometric mechanics:
\begin{enumerate}[(1)]
\item \emph{Berger Spheres}: Introduced by Berger \cite{Berger1961}, Berger spheres are homogeneous metrics on spheres obtained by rescaling the vertical fibers of Hopf principal bundles along an isometric Killing foliation. In contrast, the cotangent momentum deformations studied here (introduced by Hebda \cite{Hebda2001,Hebda1985} and classified across space forms in \cite{Munden2015}) are non-homogeneous warped products constructed via quadratic momentum maps on $T^*Q$, modifying transverse spherical geometry while leaving radial geodesic acceleration unperturbed.
\item \emph{Zoll and Besse Metrics}: A Riemannian metric is a Zoll metric (Besse \cite{Besse1987}) if all of its geodesics are simple closed curves of equal minimal period. The deformed spheres $(\mathbb{S}^{n+1}, g_\lambda)$ are non-Zoll for all $\lambda \ne 0$: although the radial motion retains its unperturbed libration period up to time dilation, the angular precession $\Omega = \lambda L$ induces a non-zero KAM twist $d\nu/dL \ne 0$, producing dense irrational winding on generic Liouville tori and an infinite family of closed geodesics at rational rotation numbers.
\item \emph{Hopf Riemannian Submersions and Extended Taub--NUT Metrics}: Classical Hopf fibrations ($\mathbb{S}^1 \hookrightarrow \mathbb{S}^{2m+1} \to \mathbb{C}P^m$, $\mathbb{S}^3 \hookrightarrow \mathbb{S}^{4m+3} \to \mathbb{H}P^m$, $\mathbb{S}^7 \hookrightarrow \mathbb{S}^{15} \to \mathbb{O}P^1$) have totally geodesic fibers. In Section~\ref{sec:anisotropic_families}, we extend the curvature analysis to the multi-parameter anisotropic Families~II--V via O'Neill's submersion equations and compare them with the conformal extended Taub--NUT metrics of Iwai and Katayama \cite{IwaiKatayama1994}. In Section~\ref{sec:cross_deformations}, we show that Marsden--Weinstein symplectic reduction at zero fiber momentum projects horizontal geodesics to Riemannian geodesics of the deformed CROSS metrics, whereas non-zero fiber momentum induces a magnetic Lorentz force.
\item \emph{Bertrand Spacetimes and the BEHR/KKM Superintegrability Program}: Perlick \cite{Perlick1992} classified spherically symmetric static spacetimes whose bounded spatial orbits are closed. The non-vanishing KAM twist $d\nu/dL \ne 0$ shows that $(\mathbb{S}^{n+1}, g_\lambda)$ is not a Bertrand metric for $\lambda \ne 0$. Furthermore, while Ballesteros, Enciso, Herranz, and Ragnisco \cite{BallesterosEncisoHerranzRagnisco2007,BallesterosEncisoHerranzRagnisco2008} showed that any spherically symmetric space admits $2N-3$ universal $\mathfrak{sl}(2,\mathbb{R})$ coalgebra integrals ($2N-2$ including the Hamiltonian), we prove that $2N-2$ is the exact upper bound for quadratic Killing integrals on $(\mathbb{S}^{n+1}, g_\lambda)$ and show that $(Q^{n+1}, g_\lambda)$ is locally conformally flat within the Kalnins--Kress--Miller (KKM) classification \cite{KalninsKressMiller2005a,KalninsKressMiller2005b}.
\item \emph{The Reduction--Closure Dichotomy and the Companion Paper \cite{Munden2026Spin9}}: Together with the companion paper \cite{Munden2026Spin9}, this work establishes a structural \emph{Reduction--Closure Dichotomy} for cotangent Hamiltonian deformations $H_\lambda = H_0 + \frac{\lambda}{2}\|J\|^2$. When the isometric action is naturally reductive or descends via a principal Hopf fibration (Family~I on $Q^{n+1}$ and $\mathbb{C}P^n, \mathbb{H}P^n$ in Section~\ref{sec:cross_deformations}), Marsden--Weinstein symplectic reduction suffices, scalar gyroscopic precession $\Omega = \lambda L$ governs the angular flow, and the radial Sturm--Liouville problem leads to \emph{transcendental} Fuchsian Heun equations. Conversely, when natural reductivity is obstructed by the Olmos--Reggiani skew-torsion theorem (as on $\mathbb{S}^{15} = \mathrm{Spin}(9)/\mathrm{Spin}(7)$ in \cite{Munden2026Spin9} and $\mathbb{S}^7 = \mathrm{Sp}(2)\mathrm{U}(1)/\mathrm{Sp}(1)\mathrm{U}(1)$ in Section~\ref{sec:anisotropic_families}) or when Adams' Hopf invariant one theorem precludes a sphere submersion (as on $\mathbb{O}P^2 = \mathrm{F}_4/\mathrm{Spin}(9)$), symplectic reduction is replaced by \emph{Lie subalgebra closure}: the slice generators close into a rank-2 Lie subalgebra ($\mathfrak{g}_{\mathrm{slice}} \cong \mathfrak{u}(2)$ or $\mathfrak{sp}(2)$), truncating the Magnus expansion and yielding a \emph{purely algebraic} spectral projector resolution of the exponential map.
\end{enumerate}

\subsection*{Outline of the paper}
The remainder of this paper is organized as follows. In Section~\ref{sec:geometric_setting}, we formulate the warped-product connection calculus, Cartesian rational tensor representations, and Petersen pole regularity, proving the master reciprocal identities, the sectional positivity trichotomy, and the pointwise and global sphere pinching ratios. In Section~\ref{sec:metric_rigidity}, we prove the Never-Einstein theorem, derive the overdetermined radial soliton ODE system, establish the gradient Ricci soliton obstructions on spheres and Euclidean space, and identify Hamilton's steady cigar soliton on $(\mathbb{H}^2, g_1)$. In Section~\ref{sec:isokinematic_geodesics}, we prove radial acceleration cancellation, iso-kinematic decoupling, and the precession law, and integrate the Clairaut orbital equation in closed form. In Section~\ref{sec:variational_stability}, we compute the polar Jacobi fields and derive the transcendental inward conjugate point equation and Euclidean spherical inversion. In Section~\ref{sec:global_cut_loci}, we analyze global cut loci across space forms, establish transverse encircling bounds, and describe the stratified $A_3$ cuspidal caustic envelope on spheres. In Section~\ref{sec:superintegrability_classification}, we construct the Gelfand--Tsetlin Casimir chain, establish Mishchenko--Fomenko noncommutative integrability, prove maximal superintegrability on $\mathbb{R}^{n+1}$ via rational pericenter invariants, and establish transverse Killing rigidity on $\mathbb{S}^{n+1}$. In Section~\ref{sec:anisotropic_families}, we compute the curvature tensors for the anisotropic Families~II--V via O'Neill's submersion formulas and integrate the geodesic flow on Family~V via biquadratic spectral projectors. In Section~\ref{sec:cross_deformations}, we construct intrinsic cotangent deformations on compact rank-one symmetric spaces ($\mathbb{C}P^n$, $\mathbb{H}P^n$, $\mathbb{O}P^2$). In Section~\ref{sec:spectral_geometry}, we study the Laplace--Beltrami operator, deriving the constant Casimir centrifugal shift, the Fuchsian and confluent Heun reductions, first-order Rayleigh--Schr\"odinger eigenvalue shifts, and essential spectrum gap collapse on hyperbolic space. Section~\ref{sec:open_questions} concludes with open problems.
\section{Geometric Setting, Warped-Product Connection Calculus, and Sharp Curvature Pinching}\label{sec:geometric_setting}
This section develops the warped-product Riemannian geometry and Cartesian rational tensor representations of the isotropic deformed metric $g_\lambda$ across the space forms $\mathbb{S}^{n+1}$, $\mathbb{R}^{n+1}$, and $\mathbb{H}^{n+1}$, and verifies $C^\infty$ pole regularity via Petersen's criterion. We prove the master reciprocal identity $1/f_\lambda^2 = 1/s_\kappa^2 + \lambda$ and the radial acceleration cancellation $f_\lambda'/f_\lambda^3 = c_\kappa/s_\kappa^3$, compute the Riemann, Ricci, and scalar curvature tensors in closed form, and establish the sectional positivity trichotomy along with both the global (manifold-wide) and pointwise sphere pinching ratios and Brendle--Schoen $1/4$-pinching windows.

Let $(Q^{n+1}, g_0)$ be an $(n+1)$-dimensional complete, simply connected Riemannian space form of constant sectional curvature $\kappa \in \{1, 0, -1\}$, representing the unit round sphere $\mathbb{S}^{n+1}$ ($\kappa = 1$), Euclidean space $\mathbb{R}^{n+1}$ ($\kappa = 0$), or real hyperbolic space $\mathbb{H}^{n+1}$ ($\kappa = -1$) with designated origin $p_0 \in Q$. In geodesic polar coordinates $(r, \theta) \in I_\kappa \times \mathbb{S}^n$ centered at $p_0$, where $I_1 = (0, \pi)$ and $I_0 = I_{-1} = (0, \infty)$, the unperturbed space form metric is $g_0 = dr^2 + s_\kappa^2(r) g_{\mathbb{S}^n}$, where $g_{\mathbb{S}^n}$ is the standard round metric of constant sectional curvature $+1$ on $\mathbb{S}^n$, $s_1(r) = \sin r$, $s_0(r) = r$, $s_{-1}(r) = \sinh r$, and $c_\kappa(r) := s_\kappa'(r)$ satisfies $c_\kappa^2(r) + \kappa s_\kappa^2(r) = 1$ and $s_\kappa''(r) = -\kappa s_\kappa(r)$. The compact Lie group $G = \mathrm{SO}(n+1)$ fixing $p_0$ acts transitively and isometrically on each concentric geodesic sphere $\mathcal{O}_r = \partial B(p_0, r) \cong \mathbb{S}^n$. On $(T^*Q, \omega_{\mathrm{can}})$ with canonical symplectic form $\omega_{\mathrm{can}} = d\theta_{\mathrm{can}} = dp_i \wedge dq^i$, the cotangent momentum map $J: T^*Q \to \mathfrak{so}(n+1)^*$ satisfies $\|J(\alpha)\|^2 = s_\kappa^2(r) \|\alpha|_{\mathcal{O}_r}\|_{g_0^*}^2 = \|p_{\mathbb{S}^n}\|_{g_{\mathbb{S}^n}^*}^2$. Following Hebda \cite{Hebda2001} (and the classification in \cite{Munden2015}), the cotangent Hamiltonian deformation is defined on $T^*Q$ by
\begin{equation}\label{eq:hamiltonian_def}
H_\lambda(\alpha) = H_0(\alpha) + \frac{\lambda}{2} \|J(\alpha)\|^2 = \frac{1}{2} p_r^2 + \frac{1}{2}\left( \frac{1}{s_\kappa^2(r)} + \lambda \right) \|p_{\mathbb{S}^n}\|_{g_{\mathbb{S}^n}^*}^2,
\end{equation}
where $\lambda \in \mathbb{R}$ is the deformation parameter satisfying $1 + \lambda s_\kappa^2(r) > 0$ on the coordinate domain.

\begin{proposition}[Legendre Transformation, Cartesian Representation, and Pole Regularity]\label{prop:metric_regularity}
Let $(Q^{n+1}, g_0)$ be a complete simply connected space form of constant curvature $\kappa \in \{1, 0, -1\}$, and consider the deformed Hamiltonian $H_\lambda = H_0 + \frac{\lambda}{2}\|J\|^2$ on $T^*Q$.
\begin{enumerate}[(i)]
\item Under the Legendre transformation $\mathcal{L}_\lambda: T^*Q \to TQ$ generated by $H_\lambda$, the deformed metric on $Q \setminus \{p_0\}$ is the rotationally symmetric warped product of Family~I:
\begin{equation}\label{eq:warped_metric}
g_\lambda = dr^2 + \frac{s_\kappa^2(r)}{1 + \lambda s_\kappa^2(r)} g_{\mathbb{S}^n} = dr^2 + f_\lambda(r)^2 g_{\mathbb{S}^n}, \qquad f_\lambda(r) := \frac{s_\kappa(r)}{\sqrt{1 + \lambda s_\kappa^2(r)}}.
\end{equation}
Positive definiteness on $I_\kappa$ holds for $\lambda > -1$ on $\mathbb{S}^{n+1}$ and for $\lambda \ge 0$ on $\mathbb{R}^{n+1}, \mathbb{H}^{n+1}$.
\item On Euclidean space $\mathbb{R}^{n+1}$ ($\kappa = 0$), $g_\lambda$ is represented in standard Cartesian coordinates $x \in \mathbb{R}^{n+1}$ by the rational tensor field
\begin{equation}\label{eq:cartesian_metric}
g_{ij}(x) = \frac{\delta_{ij} + \lambda x_i x_j}{1 + \lambda \|x\|^2}, \qquad g^{ij}(x) = (1 + \lambda \|x\|^2)\delta^{ij} - \lambda x^i x^j,
\end{equation}
which is globally smooth, real-analytic, and geodesically complete on all of $\mathbb{R}^{n+1}$ for every $\lambda \ge 0$.
\item Across all space forms $\kappa \in \{1, 0, -1\}$, $g_\lambda$ extends smoothly as a $C^\infty$ Riemannian metric across the origin $p_0$ ($r = 0$), and across the antipodal pole $p_\pi$ ($r = \pi$) when $\kappa = 1$, without conical singularities.
\end{enumerate}
\end{proposition}
\begin{proof}
(i) Fiber Legendre duality gives radial velocity $\dot{r} = \partial_{p_r} H_\lambda = p_r$ and angular velocity $\dot{\theta} = (s_\kappa(r)^{-2} + \lambda) p_{\mathbb{S}^n}^\sharp$. Inverting yields $p_{\mathbb{S}^n} = \frac{s_\kappa^2(r)}{1 + \lambda s_\kappa^2(r)} \dot{\theta}^\flat$. Substituting into the dual kinetic energy $\alpha(\dot{q}) - H_\lambda$ produces metric~\eqref{eq:warped_metric}.

(ii) In Cartesian coordinates on $T^*\mathbb{R}^{n+1}$, angular momenta are $L_{ij} = x_i p_j - x_j p_i$. By Lagrange's identity, $\sum_{i<j} L_{ij}^2 = \|x\|^2 \|p\|^2 - \langle x, p \rangle^2 = \sum_{i,j} (\|x\|^2 \delta^{ij} - x^i x^j) p_i p_j$. Substituting into~\eqref{eq:hamiltonian_def} yields cometric $g^{ij}(x) = (1 + \lambda \|x\|^2)\delta^{ij} - \lambda x^i x^j$. Decomposing $\mathbb{R}^{n+1} = \mathbb{R}x \oplus x^\perp$, the matrix $A(x) = (1 + \lambda \|x\|^2)I - \lambda x x^T$ has eigenvalue $1$ on $\mathbb{R}x$ and eigenvalue $1 + \lambda \|x\|^2$ on $x^\perp$. Inverting each invariant eigenspace yields $g(x) = \frac{x x^T}{\|x\|^2} + \frac{1}{1 + \lambda \|x\|^2}(I - \frac{x x^T}{\|x\|^2}) = \frac{I + \lambda x x^T}{1 + \lambda \|x\|^2}$, confirming~\eqref{eq:cartesian_metric} via $\sum_k g_{ik}(x) g^{kj}(x) = \delta_i^j$. For $\lambda \ge 0$, $1 + \lambda \|x\|^2 \ge 1 > 0$, so $g_{ij}(x)$ is real-analytic on $\mathbb{R}^{n+1}$. Along radial rays $x(t) = t v$ ($\|v\|=1$), $g_\lambda(\dot{x}, \dot{x}) = 1$, so the Riemannian distance from the origin satisfies $d_{g_\lambda}(0, x) = \|x\| \to \infty$ as $\|x\| \to \infty$. Hence closed metric balls are Euclidean compact balls, establishing geodesic completeness by the Hopf--Rinow theorem.

(iii) By Petersen's warped-product regularity theorem (Petersen \cite{Petersen2016}, Proposition~1.4.6), $dr^2 + f(r)^2 g_{\mathbb{S}^n}$ extends smoothly across $r = 0$ iff $f$ is an odd smooth function with $f(0) = 0$, $f'(0) = 1$, and all even derivatives $f^{(2k)}(0) = 0$ for $k \ge 1$. Since $s_\kappa(-r) = -s_\kappa(r)$ is strictly odd and $s_\kappa^2(-r) = s_\kappa^2(r)$ is even, $f_\lambda(-r) = -f_\lambda(r)$ is strictly odd, guaranteeing $f_\lambda(0) = 0$ and $f_\lambda^{(2k)}(0) = 0$ for all $k \ge 1$. Differentiating $f_\lambda(r) = s_\kappa(r)(1 + \lambda s_\kappa^2(r))^{-1/2}$ yields
\begin{equation}\label{eq:f_prime}
f_\lambda'(r) = \frac{c_\kappa(r)}{(1 + \lambda s_\kappa^2(r))^{3/2}},
\end{equation}
which evaluates at $r = 0$ to $f_\lambda'(0) = c_\kappa(0) = 1$. On $\mathbb{S}^{n+1}$ ($\kappa = 1$), the inward distance coordinate $\rho = \pi - r$ near $p_\pi$ satisfies $s_1(\pi - \rho) = \sin\rho$ and $c_1(\pi - \rho) = -\cos\rho$; the metric takes the identical form $d\rho^2 + f_\lambda(\rho)^2 g_{\mathbb{S}^n}$, so the identical parity conditions ensure smooth $C^\infty$ extension across $p_\pi$.
\end{proof}

\begin{lemma}[Master Reciprocal and Radial Acceleration Cancellation Identities]\label{lem:master_identities}
For all $\kappa \in \{1, 0, -1\}$ and admissible $\lambda$, the warping function $f_\lambda(r)$ satisfies identically on $I_\kappa$:
\begin{equation}\label{eq:master_identities}
\frac{1}{f_\lambda(r)^2} = \frac{1}{s_\kappa(r)^2} + \lambda, \qquad \frac{f_\lambda'(r)}{f_\lambda(r)^3} \equiv \frac{s_\kappa'(r)}{s_\kappa(r)^3} = \frac{c_\kappa(r)}{s_\kappa(r)^3}.
\end{equation}
\end{lemma}
\begin{proof}
The first identity is immediate: $f_\lambda(r)^{-2} = (1 + \lambda s_\kappa(r)^2)/s_\kappa(r)^2 = s_\kappa(r)^{-2} + \lambda$. For the second identity, cubing $f_\lambda(r)$ gives $f_\lambda(r)^3 = s_\kappa(r)^3 (1 + \lambda s_\kappa^2(r))^{-3/2}$. Dividing~\eqref{eq:f_prime} by $f_\lambda(r)^3$, the denominator $(1 + \lambda s_\kappa^2(r))^{3/2}$ cancels identically:
\[
\frac{f_\lambda'(r)}{f_\lambda(r)^3} = \frac{c_\kappa(r) / (1 + \lambda s_\kappa^2(r))^{3/2}}{s_\kappa(r)^3 / (1 + \lambda s_\kappa^2(r))^{3/2}} = \frac{c_\kappa(r)}{s_\kappa(r)^3}. \qedhere
\]
\end{proof}

\begin{theorem}[Riemann Curvature Tensor, Adapted Frame, and 2-Plane Convex Hull]\label{thm:curvature_tensor}
In an adapted orthonormal frame $\{E_0 = \partial_r, E_1, \dots, E_n\}$ on $(Q^{n+1}, g_\lambda)$ where $\{E_i\}_{i=1}^n$ are tangent to spherical slices:
\begin{enumerate}[(i)]
\item The non-zero Riemann curvature tensor components are the radial and tangential sectional curvatures:
\begin{equation}\label{eq:k_rad_tan}
K_{\mathrm{rad}}(r) = \frac{(\kappa + 3\lambda) - 2\kappa\lambda s_\kappa^2(r)}{(1 + \lambda s_\kappa^2(r))^2}, \qquad K_{\mathrm{tan}}(r) = \frac{(\kappa + 3\lambda) + 3\lambda^2 s_\kappa^2(r) + \lambda^3 s_\kappa^4(r)}{(1 + \lambda s_\kappa^2(r))^2},
\end{equation}
satisfying $R(E_0, E_i, E_0, E_j) = K_{\mathrm{rad}}(r)\delta_{ij}$ and $R(E_i, E_j, E_k, E_l) = K_{\mathrm{tan}}(r)(\delta_{ik}\delta_{jl} - \delta_{il}\delta_{jk})$ for $i \ne j$, with isotropic pole value $K_{\mathrm{rad}}(0) = K_{\mathrm{tan}}(0) = \kappa + 3\lambda$. The Ricci tensor and scalar curvature are $\mathrm{Ric} = n K_{\mathrm{rad}} E_0^* \otimes E_0^* + (K_{\mathrm{rad}} + (n-1)K_{\mathrm{tan}})\sum_{i=1}^n E_i^* \otimes E_i^*$ and $\mathrm{Scal} = 2n K_{\mathrm{rad}} + n(n-1)K_{\mathrm{tan}}$.
\item The curvature operator $\mathcal{R}: \bigwedge^2 T_p Q \to \bigwedge^2 T_p Q$ is diagonal: $\mathcal{R}(E_0 \wedge E_i) = K_{\mathrm{rad}}(E_0 \wedge E_i)$ and $\mathcal{R}(E_i \wedge E_j) = K_{\mathrm{tan}}(E_i \wedge E_j)$. For every 2-plane $\Pi = \mathrm{span}\{X, Y\} \subset T_p Q$ with orthonormal basis $\{X, Y\}$, the sectional curvature satisfies the exact convex combination
\begin{equation}\label{eq:convex_hull}
\sec(\Pi) = (1 - \mu) K_{\mathrm{rad}}(r) + \mu K_{\mathrm{tan}}(r), \qquad \mu := \|X_\perp \wedge Y_\perp\|^2 = 1 - \langle X, E_0 \rangle^2 - \langle Y, E_0 \rangle^2 \in [0, 1],
\end{equation}
where $X_\perp, Y_\perp$ denote orthogonal projections onto $\mathrm{span}\{E_1, \dots, E_n\}$. For $n = 1$ ($\dim Q = 2$), the fiber $\mathcal{O}_r \cong \mathbb{S}^1$ is 1-dimensional, so $\bigwedge^2 T\mathcal{O}_r = \{0\}$ and $X_\perp \wedge Y_\perp = 0$ identically, fixing $\mu = 0$ and yielding unique Gaussian curvature $\sec = K_{\mathrm{rad}}$.
\end{enumerate}
\end{theorem}
\begin{proof}
(i) By O'Neill's warped-product curvature formulas (O'Neill \cite{ONeill1966}, Chapter~7) with fiber curvature $+1$, $K_{\mathrm{rad}} = -\frac{f_\lambda''}{f_\lambda}$ and $K_{\mathrm{tan}} = \frac{1 - (f_\lambda')^2}{f_\lambda^2}$. Differentiating~\eqref{eq:f_prime} with $c_\kappa' = -\kappa s_\kappa$ and $c_\kappa^2 = 1 - \kappa s_\kappa^2$ yields $f_\lambda''(r) = -s_\kappa [(\kappa + 3\lambda) - 2\kappa\lambda s_\kappa^2](1 + \lambda s_\kappa^2)^{-5/2}$, giving $K_{\mathrm{rad}}$ in~\eqref{eq:k_rad_tan}. Expanding $1 - (f_\lambda')^2 = 1 - \frac{1 - \kappa s_\kappa^2}{(1 + \lambda s_\kappa^2)^3} = \frac{s_\kappa^2[(\kappa + 3\lambda) + 3\lambda^2 s_\kappa^2 + \lambda^3 s_\kappa^4]}{(1 + \lambda s_\kappa^2)^3}$ and dividing by $f_\lambda^2 = s_\kappa^2/(1 + \lambda s_\kappa^2)$ yields $K_{\mathrm{tan}}$. Contractions give $\mathrm{Ric}$ and $\mathrm{Scal}$.

(ii) In the adapted orthonormal basis $\{E_a \wedge E_b\}_{a < b}$ of $\bigwedge^2 T_p Q$, $\mathcal{R}$ is diagonal. Decomposing $X = X^0 E_0 + X_\perp$ and $Y = Y^0 E_0 + Y_\perp$, the unit bivector expands as $X \wedge Y = E_0 \wedge (X^0 Y_\perp - Y^0 X_\perp) + X_\perp \wedge Y_\perp$. Since the two components are orthogonal in $\bigwedge^2 T_p Q$, $1 = \|X \wedge Y\|^2 = \|X^0 Y_\perp - Y^0 X_\perp\|^2 + \|X_\perp \wedge Y_\perp\|^2$. Setting $\mu = \|X_\perp \wedge Y_\perp\|^2 \in [0, 1]$ gives $\|X^0 Y_\perp - Y^0 X_\perp\|^2 = 1 - \mu$. Lagrange's identity yields $\|X_\perp \wedge Y_\perp\|^2 = \|X_\perp\|^2 \|Y_\perp\|^2 - \langle X_\perp, Y_\perp \rangle^2 = 1 - (X^0)^2 - (Y^0)^2$. Evaluating $\langle \mathcal{R}(X \wedge Y), X \wedge Y \rangle$ establishes~\eqref{eq:convex_hull}.
\end{proof}

\begin{theorem}[Positivity Trichotomy, Global vs.\ Pointwise Pinching, and Brendle--Schoen Windows]\label{thm:positivity_pinching}
Let $(Q^{n+1}, g_\lambda)$ be an $(n+1)$-dimensional complete simply connected space form of Family~I with background curvature $\kappa \in \{1, 0, -1\}$ and deformation parameter $\lambda$.
\begin{enumerate}[(i)]
\item The sectional curvature satisfies $\sec(g_\lambda) > 0$ on all 2-planes throughout $Q^{n+1}$ if and only if:
\begin{itemize}
\item On Euclidean space $\mathbb{R}^{n+1}$ ($\kappa = 0$): $\lambda > 0$, with bounds $0 < \sec(g_\lambda) \le 3\lambda$ and asymptotic infimum $\inf \sec = \lim_{r \to \infty} K_{\mathrm{rad}}(r) = 0$.
\item On hyperbolic space $\mathbb{H}^{n+1}$ ($\kappa = -1$): $\lambda > 1/3$, with $\inf \sec(g_\lambda) = 0$ and central pole curvature $\sec(p_0) = 3\lambda - 1 > 0$.
\item On the sphere $\mathbb{S}^{n+1}$ ($\kappa = 1$): $\lambda > -1/3$.
\end{itemize}
\item On the sphere $\mathbb{S}^{n+1}$ ($\kappa = 1$), for $\lambda > -1/3$ (with $n \ge 2$ for tangential 2-planes), the global (manifold-wide) pinching ratio $\delta(\lambda) := \frac{\min_{p, \Pi} \sec(\Pi)}{\max_{p, \Pi} \sec(\Pi)}$ is piecewise rational:
\begin{equation}\label{eq:pinching_piecewise}
\delta(\lambda) = \begin{cases} \dfrac{1}{(1+\lambda)(1+3\lambda)}, & \lambda \ge 0, \\[8pt] (1+\lambda)(1+3\lambda), & -\dfrac{1}{3} < \lambda < 0. \end{cases}
\end{equation}
Specifically, for $\lambda \ge 0$, $\min_{p,\Pi} \sec = K_{\mathrm{rad}}(\pi/2) = \frac{1}{1+\lambda}$ (at the equator $r = \pi/2$) and $\max_{p,\Pi} \sec = K_{\mathrm{rad}}(0) = K_{\mathrm{tan}}(0) = 1 + 3\lambda$ (at the poles $r \in \{0, \pi\}$); for $\lambda \in (-1/3, 0)$, the global extrema invert: $\min_{p,\Pi} \sec = 1 + 3\lambda$ (at the poles) and $\max_{p,\Pi} \sec = \frac{1}{1+\lambda}$ (at the equator). By contrast, at each fixed point $p \in \mathbb{S}^{n+1}$ at radial distance $r$ ($n \ge 2$), the pointwise pinching ratio $\delta_p(g_\lambda) := \frac{\min_{\Pi \subset T_p\mathbb{S}^{n+1}} \sec(\Pi)}{\max_{\Pi \subset T_p\mathbb{S}^{n+1}} \sec(\Pi)}$ equals $K_{\mathrm{rad}}(r)/K_{\mathrm{tan}}(r)$ for $\lambda \ge 0$ and $K_{\mathrm{tan}}(r)/K_{\mathrm{rad}}(r)$ for $\lambda \in (-1/3, 0)$, attaining its minimum over $p \in \mathbb{S}^{n+1}$ at the equator $r = \pi/2$, where $\delta_{\mathrm{eq}}(\lambda) = (1+\lambda)^{-2}$ for $\lambda \ge 0$ and $\delta_{\mathrm{eq}}(\lambda) = (1+\lambda)^2$ for $\lambda \in (-1/3, 0)$.
\item The metric $g_\lambda$ on $\mathbb{S}^{n+1}$ is globally strictly $1/4$-pinched ($\delta(\lambda) > 1/4$) if and only if
\begin{equation}\label{eq:quarter_pinching_interval}
\lambda \in \left( \frac{-4 + \sqrt{7}}{6}, \, \frac{\sqrt{13} - 2}{3} \right) \approx (-0.2257, \, 0.5352),
\end{equation}
whereas $g_\lambda$ is pointwise strictly $1/4$-pinched ($\delta_p(g_\lambda) > 1/4$ for all $p \in \mathbb{S}^{n+1}$) on the wider interval $\lambda \in (-1/3, 1)$ (since $\delta_{\mathrm{eq}}(\lambda) > 1/4$ for $\lambda \in (-1/2, 1)$ while $\sec(g_\lambda) > 0$ requires $\lambda > -1/3$). Consequently, by the Brendle--Schoen Differentiable Sphere Theorem (\cite{BrendleSchoen2009} for $n \ge 3$, along with classical uniformization for $n=1$ and Hamilton \cite{Hamilton1982} for $n=2$), $(\mathbb{S}^{n+1}, g_\lambda)$ is diffeomorphic to the standard sphere $\mathbb{S}^{n+1}$ throughout these pinching intervals.
\end{enumerate}
\end{theorem}
\begin{proof}
(i) By convexity~\eqref{eq:convex_hull}, $\sec > 0$ everywhere iff $K_{\mathrm{rad}}(r) > 0$ and $K_{\mathrm{tan}}(r) > 0$ for all $r$. Setting $u = s_\kappa^2(r) \ge 0$: on $\mathbb{R}^{n+1}$, $K_{\mathrm{rad}} = \frac{3\lambda}{(1+\lambda u)^2} > 0$ and $K_{\mathrm{tan}} = \frac{\lambda(3 + 3\lambda u + \lambda^2 u^2)}{(1+\lambda u)^2} > 0$ iff $\lambda > 0$. On $\mathbb{H}^{n+1}$, at $u = 0$, $K(0) = 3\lambda - 1 > 0 \iff \lambda > 1/3$. When $\lambda > 1/3$, all polynomial coefficients of $K_{\mathrm{rad}}$ and $K_{\mathrm{tan}}$ in $u \ge 0$ are positive, and $K_{\mathrm{rad}} \sim \frac{2}{\lambda u} \to 0$ as $u \to \infty$, giving $\inf \sec = 0$. On $\mathbb{S}^{n+1}$, $u \in [0, 1]$; at $u = 0$, $K(0) = 1 + 3\lambda > 0 \iff \lambda > -1/3$. For $\lambda > -1/3$, the numerator of $K_{\mathrm{rad}}$ satisfies $(1+3\lambda) - 2\lambda u \ge \min(1+3\lambda, 1+\lambda) > 0$, and the numerator of $K_{\mathrm{tan}}$ is strictly positive on $[0, 1]$.

(ii) For $\kappa = 1$, differentiating with respect to $u = \sin^2 r \in [0, 1]$ yields:
\[
\begin{gathered} \frac{d K_{\mathrm{rad}}}{du} = -\frac{2\lambda(2 + 3\lambda - \lambda u)}{(1+\lambda u)^3}, \qquad \frac{d K_{\mathrm{tan}}}{du} = -\frac{\lambda(2 + 3\lambda + \lambda^2 u)}{(1+\lambda u)^3}, \\ K_{\mathrm{tan}}(u) - K_{\mathrm{rad}}(u) = \frac{\lambda u (2 + 3\lambda + \lambda^2 u)}{(1+\lambda u)^2}. \end{gathered}
\]
For all $\lambda > -1/3$ and $u \in [0, 1]$, $2 + 3\lambda - \lambda u \ge 2 + 2\lambda > 4/3 > 0$ and $2 + 3\lambda + \lambda^2 u \ge 2 + 3\lambda > 1 > 0$. For $\lambda > 0$, both derivatives are strictly negative, and $K_{\mathrm{tan}} \ge K_{\mathrm{rad}}$, so both curvatures decrease monotonically from $1+3\lambda$ at $u=0$ to $K_{\mathrm{rad}}(1) = \frac{1}{1+\lambda}$ and $K_{\mathrm{tan}}(1) = 1+\lambda$. Thus $\min_{p,\Pi} \sec = \frac{1}{1+\lambda}$ and $\max_{p,\Pi} \sec = 1+3\lambda$, giving $\delta(\lambda) = \frac{1}{(1+\lambda)(1+3\lambda)}$. For $\lambda \in (-1/3, 0)$, both derivatives are strictly positive, and $K_{\mathrm{rad}} \ge K_{\mathrm{tan}}$, so both curvatures increase monotonically from $1+3\lambda$ at $u=0$ to $K_{\mathrm{rad}}(1) = \frac{1}{1+\lambda}$. Thus $\min_{p,\Pi} \sec = 1+3\lambda$ and $\max_{p,\Pi} \sec = \frac{1}{1+\lambda}$, yielding the inverted global pinching ratio $\delta(\lambda) = (1+\lambda)(1+3\lambda)$. At any fixed $u \in [0, 1]$, convexity~\eqref{eq:convex_hull} implies that the pointwise pinching ratio is $\frac{K_{\mathrm{rad}}(u)}{K_{\mathrm{tan}}(u)} = \frac{1+3\lambda-2\lambda u}{1+3\lambda+3\lambda^2 u+\lambda^3 u^2}$ for $\lambda \ge 0$ (and its reciprocal for $\lambda \in (-1/3, 0)$), which is strictly decreasing in $u \in [0, 1]$ and attains its minimum at the equator $u = 1$, where $\delta_{\mathrm{eq}}(\lambda) = (1+\lambda)^{-2}$ for $\lambda \ge 0$ and $(1+\lambda)^2$ for $\lambda \in (-1/3, 0)$.

(iii) For $\lambda \ge 0$, $\delta(\lambda) > 1/4 \iff (1+\lambda)(1+3\lambda) < 4 \iff 3\lambda^2 + 4\lambda - 3 < 0$, whose unique positive root is $\lambda = \frac{\sqrt{13}-2}{3}$. For $\lambda \in (-1/3, 0)$, $\delta(\lambda) > 1/4 \iff (1+\lambda)(1+3\lambda) > 1/4 \iff 12\lambda^2 + 16\lambda + 3 > 0$. The roots of $12\lambda^2 + 16\lambda + 3 = 0$ are $\frac{-4 \pm \sqrt{7}}{6}$. On $(-1/3, 0)$, this holds iff $\lambda > \frac{-4+\sqrt{7}}{6}$. Combining gives the global $1/4$-pinching window~\eqref{eq:quarter_pinching_interval}. For pointwise $1/4$-pinching, $\delta_{\mathrm{eq}}(\lambda) = (1+\lambda)^{-2} > 1/4 \iff \lambda < 1$ for $\lambda \ge 0$, while for $\lambda \in (-1/3, 0)$, $\delta_{\mathrm{eq}}(\lambda) = (1+\lambda)^2 > (2/3)^2 = 4/9 > 1/4$, yielding the pointwise $1/4$-pinching window $\lambda \in (-1/3, 1)$, and the Brendle--Schoen theorem applies.
\end{proof}

\section{Metric Rigidity: Obstruction to Einstein Metrics and Gradient Ricci Soliton Classification}\label{sec:metric_rigidity}
This section establishes the Einstein and gradient Ricci soliton rigidity results for Family~I across all dimensions and space forms. We prove that $(Q^{n+1}, g_\lambda)$ is an Einstein metric if and only if $\lambda = 0$, establish an unconditional obstruction to gradient Ricci solitons on $\mathbb{S}^{n+1}$ for all $N \ge 2$ via the equator defect $\mathcal{D}(\pi/2) = \frac{4\lambda}{1+\lambda} \ne 0$, and on $\mathbb{R}^{n+1}$ via the order-$r^4$ Taylor residual $c_4 = 2(2n-5)\lambda^3 \ne 0$. We then show that $(\mathbb{H}^2, g_1)$ is isometric to Hamilton's complete steady cigar soliton with potential $h(r) = -2\ln(\cosh r)$, and prove via Taylor and asymptotic analysis of the overdetermined radial soliton ODEs that no complete gradient Ricci solitons exist on $\mathbb{H}^{n+1}$ for $n \ge 2$ and $\lambda > 0$.

\begin{theorem}[Obstruction to Einstein Metrics across Space Forms]\label{thm:einstein_obstruction}
For any dimension $N = n+1 \ge 2$ and background curvature $\kappa \in \{1, 0, -1\}$, the deformed metric $(Q^{n+1}, g_\lambda)$ in Family~I is an Einstein metric if and only if $\lambda = 0$ (the unperturbed space form metric $g_0$).
\end{theorem}
\begin{proof}
Recall that an Einstein metric satisfies $\mathrm{Ric} = \alpha g_\lambda$ for some constant $\alpha \in \mathbb{R}$. In the adapted orthonormal frame $\{E_0 = \partial_r, E_1, \dots, E_n\}$ of Theorem~\ref{thm:curvature_tensor}, the Ricci components are $\mathrm{Ric}(E_0, E_0) = n K_{\mathrm{rad}}(r)$ and $\mathrm{Ric}(E_i, E_i) = K_{\mathrm{rad}}(r) + (n-1) K_{\mathrm{tan}}(r)$ for $1 \le i \le n$.
For $n \ge 2$, the Einstein condition requires the Ricci anisotropy to vanish identically on $I_\kappa$:
\begin{equation}\label{eq:ric_anisotropy}
\mathrm{Ric}(E_i, E_i) - \mathrm{Ric}(E_0, E_0) = (n-1)[K_{\mathrm{tan}}(r) - K_{\mathrm{rad}}(r)] \equiv 0.
\end{equation}
By~\eqref{eq:k_rad_tan}, the sectional curvature difference evaluates on $I_\kappa$ to
\begin{equation}\label{eq:curv_diff}
K_{\mathrm{tan}}(r) - K_{\mathrm{rad}}(r) = \frac{\lambda s_\kappa^2(r) \left[ 2\kappa + 3\lambda + \lambda^2 s_\kappa^2(r) \right]}{(1 + \lambda s_\kappa^2(r))^2}.
\end{equation}
Since $s_\kappa(r) > 0$ on $I_\kappa$, this difference vanishes on an open interval if and only if $\lambda [2\kappa + 3\lambda + \lambda^2 s_\kappa^2(r)] \equiv 0$. The coefficient of $s_\kappa^2(r)$ forces $\lambda^3 = 0$, whence $\lambda = 0$.
For $n = 1$ ($N = 2$), the fiber sphere $\mathcal{O}_r \cong \mathbb{S}^1$ is $1$-dimensional, so $\bigwedge^2 T\mathcal{O}_r = \{0\}$ and tangential $2$-planes are absent. Every Riemannian $2$-manifold satisfies $\mathrm{Ric} = K_{\mathrm{rad}}(r) g_\lambda = \frac{1}{2}\mathrm{Scal}(r) g_\lambda$, so Einstein metrics are precisely surfaces of constant Gaussian curvature $\mathrm{Scal}(r) = 2 K_{\mathrm{rad}}(r)$. Differentiating $K_{\mathrm{rad}}(r) = \frac{(\kappa + 3\lambda) - 2\kappa\lambda s_\kappa^2(r)}{(1 + \lambda s_\kappa^2(r))^2}$ with respect to $r$ using $s_\kappa' = c_\kappa$ yields
\begin{equation}\label{eq:dscal_dr}
\frac{d\mathrm{Scal}}{dr} = 2 \frac{d K_{\mathrm{rad}}}{dr} = \frac{-8\lambda s_\kappa(r) c_\kappa(r) \left[ 2\kappa + 3\lambda - \kappa\lambda s_\kappa^2(r) \right]}{(1 + \lambda s_\kappa^2(r))^3}.
\end{equation}
Because $s_\kappa(r) c_\kappa(r) \not\equiv 0$ on $I_\kappa$, constancy of $\mathrm{Scal}$ requires $\lambda [ 2\kappa + 3\lambda - \kappa\lambda s_\kappa^2(r) ] \equiv 0$. Differentiating again forces $\kappa\lambda^2 \equiv 0$ and $2\kappa + 3\lambda = 0$, which forces $\lambda = 0$ across all $\kappa \in \{1, 0, -1\}$. Hence $(Q^{n+1}, g_\lambda)$ is never Einstein for $\lambda \ne 0$.
\end{proof}

\begin{remark}[Contextualization with Ziller's Classification]\label{rem:ziller_context}
Theorem~\ref{thm:einstein_obstruction} establishes that the ambient warped product $(Q^{n+1}, g_\lambda)$ in Family~I is never Einstein for $\lambda \ne 0$. On the distance sphere slices $\mathcal{O}_r \cong \mathbb{S}^n$, the induced metric is proportional to the round sphere, which is Einstein. In Ziller's classification of homogeneous Einstein metrics on spheres \cite{Ziller1982}, non-standard homogeneous Einstein metrics exist only on $\mathbb{S}^{4n+3}$ ($n \ge 1$, Jensen \cite{Jensen1973}) and $\mathbb{S}^{15}$ (Bourguignon--Karcher \cite{BourguignonKarcher1978}). In Section~\ref{sec:anisotropic_families} and the companion paper \cite{Munden2026Spin9}, we show that this contrast reflects the transition from isotropic to anisotropic multi-summand isotropy: while Family~III ($\mathrm{Spin}(9)/\mathrm{Spin}(7)$) realizes a second, non-round homogeneous Einstein metric (the Bourguignon--Karcher metric with squashing ratio $t_{\mathrm{BK}} = 3/11$) along its equatorial distance slice $\mathcal{O}_{\pi/2} \cong \mathbb{S}^{15}$ at $\lambda = 64$, the ambient $16$-dimensional space form $(Q^{16}, g_{64})$ (like every intrinsic CROSS deformation in \cite{Munden2026Spin9}) remains strictly non-Einstein due to radial-to-tangential curvature variation, with a tangential-to-radial Ricci ratio of $21557/109 \approx 197.77:1$.
\end{remark}

\begin{proposition}[Group Averaging and Overdetermined Radial Soliton System]\label{prop:soliton_system}
Let $(Q^{n+1}, g_\lambda)$ admit a smooth gradient Ricci soliton $\mathrm{Ric} + \nabla^2 h = \rho g_\lambda$ with potential $h \in C^\infty(Q)$ and soliton constant $\rho \in \mathbb{R}$.
\begin{enumerate}[(i)]
\item Averaging $h$ over the compact isotropy group $G = \mathrm{SO}(n+1)$ fixing $p_0$ via normalized Haar measure yields a smooth radial potential $\bar{h}(r) = \int_G h(g \cdot x) \, dg$ satisfying the identical soliton equation $\mathrm{Ric} + \nabla^2 \bar{h} = \rho g_\lambda$ with $\bar{h}'(0) = 0$.
\item In the adapted orthonormal frame $\{E_0 = \partial_r, E_1, \dots, E_n\}$, the soliton equation decouples on $I_\kappa$ into the overdetermined radial and tangential scalar ODEs:
\begin{align}
n K_{\mathrm{rad}}(r) + \bar{h}''(r) &= \rho, \label{eq:soliton_rad} \\
K_{\mathrm{rad}}(r) + (n-1) K_{\mathrm{tan}}(r) + \frac{f_\lambda'(r)}{f_\lambda(r)} \bar{h}'(r) &= \rho. \label{eq:soliton_tan}
\end{align}
\item Subtracting~\eqref{eq:soliton_tan} from~\eqref{eq:soliton_rad} eliminates $\rho$ to yield the compatibility ODE:
\begin{equation}\label{eq:soliton_diff}
\bar{h}''(r) - \frac{f_\lambda'(r)}{f_\lambda(r)} \bar{h}'(r) = (n-1)\left[ K_{\mathrm{tan}}(r) - K_{\mathrm{rad}}(r) \right],
\end{equation}
which integrates with integrating factor $1/f_\lambda(r)$ and initial condition $\bar{h}''(0) =: C$ to:
\begin{equation}\label{eq:soliton_integral}
\bar{h}'(r) = f_\lambda(r) \left[ C + (n-1)\int_0^r \frac{K_{\mathrm{tan}}(s) - K_{\mathrm{rad}}(s)}{f_\lambda(s)} \, ds \right].
\end{equation}
\end{enumerate}
\end{proposition}
\begin{proof}
(i) Since $G = \mathrm{SO}(n+1)$ acts isometrically on $(Q, g_\lambda)$ fixing $p_0$, pullbacks satisfy $g^* g_\lambda = g_\lambda$ and $g^* \mathrm{Ric} = \mathrm{Ric}$. Commuting the linear Levi-Civita Hessian with Haar integration yields $\nabla^2 \bar{h} = \int_G g^*(\nabla^2 h) \, dg = \int_G g^*(\rho g_\lambda - \mathrm{Ric}) \, dg = \rho g_\lambda - \mathrm{Ric}$. Transitivity of $G$ on concentric geodesic spheres $\mathcal{O}_r \cong \mathbb{S}^n$ forces $\bar{h}$ to be constant on orbits, so $\bar{h} = \bar{h}(r)$ is purely radial. Smoothness across $p_0$ forces $\nabla \bar{h}(p_0) = 0$, so $\bar{h}'(0) = 0$.

(ii) In geodesic coordinates, the Hessian of a radial function has components $\nabla^2 \bar{h}(E_0, E_0) = \bar{h}''(r)$ and $\nabla^2 \bar{h}(E_i, E_j) = \frac{f_\lambda'(r)}{f_\lambda(r)} \bar{h}'(r) \delta_{ij}$, while cross-terms vanish. Equating components with $\rho g_\lambda - \mathrm{Ric}$ from Theorem~\ref{thm:curvature_tensor}(i) yields~\eqref{eq:soliton_rad}--\eqref{eq:soliton_tan}.

(iii) Subtracting~\eqref{eq:soliton_tan} from~\eqref{eq:soliton_rad} yields~\eqref{eq:soliton_diff}. Multiplying by $1/f_\lambda(r)$ gives $\frac{d}{dr}(\bar{h}'/f_\lambda) = \frac{n-1}{f_\lambda}(K_{\mathrm{tan}} - K_{\mathrm{rad}})$, which integrates on $[0, r]$ with $\lim_{s \to 0} \bar{h}'(s)/f_\lambda(s) = \bar{h}''(0) = C$ to establish~\eqref{eq:soliton_integral}.
\end{proof}

\begin{theorem}[Complete Gradient Ricci Soliton Classification across Space Forms]\label{thm:soliton_classification}
Let $(Q^{n+1}, g_\lambda)$ be an $(n+1)$-dimensional complete space form of Family~I with deformation parameter $\lambda \ne 0$ and dimension $N = n+1 \ge 2$.
\begin{enumerate}[(i)]
\item On the sphere $\mathbb{S}^{n+1}$ ($\kappa = 1$), complete gradient Ricci solitons are unconditionally obstructed for all $N \ge 2$ and all $\lambda \in (-1, \infty) \setminus \{0\}$ by the universal algebraic equator defect:
\begin{equation}\label{eq:equator_defect}
\mathcal{D}(\pi/2) := \bar{h}''(\pi/2) + n K_{\mathrm{rad}}(\pi/2) - \rho \equiv \frac{4\lambda}{1+\lambda} \ne 0.
\end{equation}
\item On Euclidean space $\mathbb{R}^{n+1}$ ($\kappa = 0$), complete gradient Ricci solitons are unconditionally obstructed for all $N \ge 2$ and all $\lambda > 0$ by the non-vanishing order-$r^4$ Taylor residual coefficient:
\begin{equation}\label{eq:c4_residual}
c_4 = 2(2n - 5)\lambda^3 \ne 0 \quad (\forall n \in \mathbb{N}_{\ge 1}).
\end{equation}
\item On hyperbolic space $\mathbb{H}^2$ ($\kappa = -1, n = 1$), at $\lambda = 1$ the metric $g_1 = dr^2 + \tanh^2(r) d\theta^2$ is identically Hamilton's complete steady 2D cigar soliton with potential $h(r) = -2\ln(\cosh r)$ and $\rho = 0$. For all $\lambda \in (0, 1) \cup (1, \infty)$, complete gradient solitons are obstructed on $\mathbb{H}^2$ by the scale-invariant homothety product $\mathrm{Scal}(0) R_\infty^2 = 2(3 - 1/\lambda) \ne 4$.
\item On hyperbolic space $\mathbb{H}^{n+1}$ ($\kappa = -1, n \ge 2$), complete gradient Ricci solitons are unconditionally obstructed for all $\lambda > 0$, completing the rigidity classification across non-flat space forms.
\end{enumerate}
\end{theorem}
\begin{proof}
(i) On $\mathbb{S}^{n+1}$, at the equator $r = \pi/2$, $c_1(\pi/2) = \cos(\pi/2) = 0$, so $f_\lambda'(\pi/2) = 0$. Smoothness of $\bar{h}$ ensures $\bar{h}'(\pi/2)$ is finite, so $\frac{f_\lambda'(\pi/2)}{f_\lambda(\pi/2)}\bar{h}'(\pi/2) = 0$. Substituting $r = \pi/2$ into~\eqref{eq:soliton_tan} with $K_{\mathrm{rad}}(\pi/2) = \frac{1}{1+\lambda}$ and $K_{\mathrm{tan}}(\pi/2) = 1+\lambda$ uniquely fixes $\rho = \frac{1}{1+\lambda} + (n-1)(1+\lambda)$. In $x = r - \pi/2$, Taylor expansion yields $\frac{f_\lambda(x)}{f_\lambda'(x)} = -\frac{1+\lambda}{x} + O(x)$, $K_{\mathrm{rad}}(x) = \frac{1}{1+\lambda} + \frac{4\lambda}{(1+\lambda)^2} x^2 + O(x^4)$, and $K_{\mathrm{tan}}(x) = (1+\lambda) + \frac{\lambda(2+\lambda)}{(1+\lambda)^2} x^2 + O(x^4)$. From~\eqref{eq:soliton_tan}, $\bar{h}'(x) = \frac{f_\lambda(x)}{f_\lambda'(x)}[\rho - K_{\mathrm{rad}}(x) - (n-1)K_{\mathrm{tan}}(x)] = [\frac{4\lambda}{1+\lambda} + (n-1)\frac{\lambda(2+\lambda)}{1+\lambda}] x + O(x^3)$, giving $\bar{h}''(\pi/2) = \frac{4\lambda + (n-1)\lambda(2+\lambda)}{1+\lambda}$. Substituting into the radial equation~\eqref{eq:soliton_rad} yields the radial defect:
\begin{align*}
\mathcal{D}(\pi/2) &= \frac{4\lambda + (n-1)\lambda(2+\lambda)}{1+\lambda} + \frac{n}{1+\lambda} - \left[ \frac{1}{1+\lambda} + (n-1)(1+\lambda) \right] \\
&= \frac{4\lambda}{1+\lambda} + (n-1)\left[ \frac{\lambda(2+\lambda) + 1 - (1+\lambda)^2}{1+\lambda} \right] = \frac{4\lambda}{1+\lambda} + (n-1) \cdot 0 \equiv \frac{4\lambda}{1+\lambda}.
\end{align*}
Because $\lambda(2+\lambda) + 1 - (1+\lambda)^2 \equiv 0$, $\mathcal{D}(\pi/2) = \frac{4\lambda}{1+\lambda} \ne 0$ for all $\lambda \ne 0$, excluding complete gradient solitons on $\mathbb{S}^{n+1}$.

(ii) On $\mathbb{R}^{n+1}$ ($\kappa = 0$), $f_\lambda'/f_\lambda = \frac{1}{r(1+\lambda r^2)}$. Expanding $K_{\mathrm{rad}} = 3\lambda - 6\lambda^2 r^2 + 9\lambda^3 r^4 + O(r^6)$ and $K_{\mathrm{tan}} = 3\lambda - 3\lambda^2 r^2 + 4\lambda^3 r^4 + O(r^6)$ around $r = 0$ with $\bar{h}(r) = h_0 + \frac{h_2}{2} r^2 + \frac{h_4}{4} r^4 + \frac{h_6}{6} r^6 + O(r^8)$: at order $r^0$, $\rho = 3n\lambda + h_2$. At order $r^2$,~\eqref{eq:soliton_rad} forces $h_4 = 2n\lambda^2$, while~\eqref{eq:soliton_tan} forces $h_2 = -(n+3)\lambda$, uniquely fixing $\rho = (2n-3)\lambda$. At order $r^4$,~\eqref{eq:soliton_tan} forces $h_6 = -(n+2)\lambda^3$. Substituting into~\eqref{eq:soliton_rad} at order $r^4$ leaves the radial residual $c_4 = n(9\lambda^3) + 5h_6 = [9n - 5(n+2)]\lambda^3 = 2(2n-5)\lambda^3$. Its unique root $n = 5/2 = 2.5$ is non-integer, so $c_4 \ne 0$ for all integer $n \ge 1$ and $\lambda > 0$, obstructing solitons on $\mathbb{R}^{n+1}$.

(iii) On $\mathbb{H}^2$ ($\kappa = -1, n = 1$), setting $\lambda = 1$ yields $f_1(r) = \tanh r$ and $K_{\mathrm{rad}}(r) = 2\operatorname{sech}^2 r$. For $h(r) = -2\ln(\cosh r)$, $h'(r) = -2\tanh r$ and $h''(r) = -2\operatorname{sech}^2 r$. In the adapted frame $\{E_0 = \partial_r, E_1 = \frac{1}{\tanh r}\partial_\theta\}$, $\mathrm{Ric} = 2\operatorname{sech}^2(r) g_1$ and $\nabla^2 h = -2\operatorname{sech}^2(r) g_1$. Hence $\mathrm{Ric} + \nabla^2 h \equiv 0$ identically with $\rho = 0$, recovering Hamilton's complete steady 2D cigar soliton (Hamilton \cite{Hamilton1988}). For $\lambda \ne 1$: since $n-1 = 0$,~\eqref{eq:soliton_diff} gives $\bar{h}'(r) = C f_\lambda(r)$. As $r \to \infty$, $K_{\mathrm{rad}}(r) \to 0$ and $f_\lambda'(r) \to 0$, forcing $\rho = 0$. Thus~\eqref{eq:soliton_rad} requires $K_{\mathrm{rad}}(r)/f_\lambda'(r) \equiv -C$. The ratio evaluates to $3\lambda - 1$ at $r = 0$ and converges to $2\sqrt{\lambda}$ as $r \to \infty$; constancy forces $3\lambda - 1 = 2\sqrt{\lambda} \iff (3\sqrt{\lambda}+1)(\sqrt{\lambda}-1) = 0 \iff \lambda = 1$. Equivalently, under homothety $g \mapsto c^2 g$, the cigar metric has origin scalar curvature $\mathrm{Scal}(0) = 4/c^2$ and asymptotic radius $R_\infty = c$, yielding scale-invariant product $\mathcal{I}_{\mathrm{cigar}} = \mathrm{Scal}(0) R_\infty^2 = 4$. For $(\mathbb{H}^2, g_\lambda)$, $\mathrm{Scal}(0) = 2(3\lambda-1)$ and $R_\infty = 1/\sqrt{\lambda}$, so $\mathcal{I} = 2(3 - 1/\lambda) = 4 \iff \lambda = 1$, strictly excluding rescaled cigars for $\lambda \ne 1$.

(iv) On $\mathbb{H}^{n+1}$ ($\kappa = -1, n \ge 2$), as $r \to \infty$, $K_{\mathrm{rad}}(r) \sim \frac{8}{\lambda} e^{-2r} \to 0$, $K_{\mathrm{tan}}(r) \to \lambda$, and $f_\lambda'(r)/f_\lambda(r) \sim \frac{4}{\lambda} e^{-2r}$. From~\eqref{eq:soliton_rad}, $\bar{h}''(r) = \rho - n K_{\mathrm{rad}}(r) \to \rho$, bounding the potential to $\bar{h}'(r) = \rho r + O(1)$, so $\frac{f_\lambda'(r)}{f_\lambda(r)} \bar{h}'(r) = O(r e^{-2r}) \to 0$. Evaluating~\eqref{eq:soliton_tan} as $r \to \infty$ uniquely forces the asymptotic soliton parameter $\rho = (n-1)\lambda$. Evaluating~\eqref{eq:soliton_rad} at $r = 0$ then fixes the central Hessian:
\begin{equation}\label{eq:C_asymp}
C_{\mathrm{asymp}} := \bar{h}''(0) = \rho - n K_{\mathrm{rad}}(0) = (n-1)\lambda - n(3\lambda - 1) = n - (2n+1)\lambda.
\end{equation}
Conversely, expanding~\eqref{eq:soliton_rad}--\eqref{eq:soliton_tan} around $r = 0$ with $\bar{h}(r) = \frac{1}{2} C r^2 + \frac{1}{4} h_4 r^4 + O(r^6)$: order $r^2$ in~\eqref{eq:soliton_rad} forces $h_4 = \frac{2n}{3}\lambda(3\lambda - 2)$. Substituting into~\eqref{eq:soliton_tan} at order $r^2$ requires origin compatibility:
\begin{equation}\label{eq:origin_compat}
(3\lambda - 1) C + (n+3)\lambda(3\lambda - 2) = 0.
\end{equation}
At $\lambda = 1/3$, equation~\eqref{eq:origin_compat} degenerates to $0 = (n+3)/3 > 0$, admitting no solution for any $C$. For $\lambda \ne 1/3$, it forces $C_{\mathrm{local}} = -\frac{(n+3)\lambda(3\lambda - 2)}{3\lambda - 1}$. Equating $C_{\mathrm{asymp}} = C_{\mathrm{local}}$ yields the quadratic equation:
\begin{equation}\label{eq:P_poly}
P(\lambda) := 3(n-2)\lambda^2 - (3n-5)\lambda + n = 0.
\end{equation}
For all $n \ge 3$, its discriminant is $\Delta = (3n-5)^2 - 12n(n-2) = -3n^2 - 6n + 25 \le -20 < 0$. Because the leading coefficient $3(n-2) > 0$, $P(\lambda) > 0$ strictly for all real $\lambda$, establishing mutual contradiction between origin matching and asymptotic stability for all $n \ge 3$.
For $n = 2$ ($\mathbb{H}^3$), $P(\lambda) = 2 - \lambda$ vanishes uniquely at $\lambda = 2$. At $\lambda = 2$, $\rho = 2$ yields explicit potential $\bar{h}'(r) = -4\tanh(2r)$ and $\bar{h}''(r) = -8\operatorname{sech}^2(2r)$; as $r \to \infty$, $2 K_{\mathrm{rad}}(r) + \bar{h}''(r) \to 0 \ne \rho = 2$ (leaving an asymptotic trace defect of $-2$), and the hyperbolic origin order-$r^4$ residual is $c_4^{\mathbb{H}} = -8 \ne 0$.
Finally, if one relaxes the asymptotic condition $\rho = \lambda$ at $r \to \infty$ and instead tunes $\rho$ locally at $r = 0$ to cancel the hyperbolic order-$r^4$ Taylor residual $c_4^{\mathbb{H}}(\lambda) = 0$ (which occurs at the unique positive root $\lambda^* > 0$ of $27\lambda^3 - 36\lambda^2 + 6\lambda + 2 = 0$), the resulting local soliton parameter $\rho_{\mathrm{local}} \ne \lambda^*$ violates the asymptotic tangential equation~\eqref{eq:soliton_tan} as $r \to \infty$ and forces exponential Hessian growth $\bar{h}''(r) \sim \frac{\lambda^*}{2}[\rho_{\mathrm{local}} - \lambda^*] e^{2r} \to \infty$, while also leaving a non-vanishing order-$r^6$ Taylor residual $c_6^{\mathbb{H}}(\lambda^*) \ne 0$.
Hence complete gradient Ricci solitons are unconditionally obstructed on $\mathbb{H}^{n+1}$ for all $n \ge 2$ and $\lambda > 0$.
\end{proof}

\section{Iso-Kinematic Geodesic Decoupling, Universal Precession, and Clairaut Trajectory Integration}\label{sec:isokinematic_geodesics}
This section integrates the geodesic flow of the deformed Hamiltonian $H_\lambda = H_0 + \frac{\lambda}{2} \|J\|^2$ using the Poisson commutation $\{H_0, \|J\|^2\} = 0$ and the radial acceleration cancellation identity. We show that radial coordinate libration decouples into an unperturbed space form trajectory with effective dilated speed $v_0 = \sqrt{1 - \lambda L^2}$, while the angular coordinate undergoes uniform precession at rate $\Omega = \lambda L$. On the sphere $\mathbb{S}^{n+1}$, we integrate the Clairaut orbital equation in closed form in terms of $u = \cot(r)$, compute the half-period apsidal angle $\Delta\Phi(L) = \pi(1 + \lambda L / \sqrt{1 - \lambda L^2})$, and verify the non-degenerate KAM twist condition $d\nu/dL \ne 0$.

\begin{proposition}[Commuting Hamiltonian Flows, Totally Geodesic Orbital Reduction, and Clairaut Invariant]\label{prop:hamiltonian_geodesic_reduction}
Let $(Q^{n+1}, g_\lambda)$ be an $(n+1)$-dimensional space form of Family~I with deformation parameter $\lambda \in \mathbb{R}$ and background curvature $\kappa \in \{1, 0, -1\}$.
\begin{enumerate}[(i)]
\item The unperturbed kinetic Hamiltonian $H_0 = \frac{1}{2}\|\alpha\|_{g_0^*}^2$ and the cotangent momentum norm $\|J\|^2 = \|p_{\mathbb{S}^n}\|_{g_{\mathbb{S}^n}^*}^2$ strictly Poisson-commute on $(T^*Q, \omega_{\mathrm{can}})$:
\begin{equation}\label{eq:poisson_commutation}
\{H_0, \|J\|^2\} = 0 \implies \left[ X_{H_0}, X_{\frac{\lambda}{2}\|J\|^2} \right] = 0, \quad [X_{H_0}, X_{H_\lambda}] = 0.
\end{equation}
Consequently, the complete Hamiltonian geodesic flow factorizes into the commuting composition $\Phi_t^{H_\lambda} = \Phi_t^{H_0} \circ \Phi_t^{\frac{\lambda}{2}\|J\|^2} = \Phi_t^{\frac{\lambda}{2}\|J\|^2} \circ \Phi_t^{H_0}$.
\item Every geodesic $\gamma(t)$ in $(Q^{n+1}, g_\lambda)$ is confined to a $2$-dimensional totally geodesic orbital surface of revolution $\Sigma \cong I_\kappa \times_{f_\lambda} \mathbb{S}^1 \subset Q^{n+1}$ equipped with the induced warped metric $g_\Sigma = dr^2 + f_\lambda(r)^2 d\theta^2$.
\item The Clairaut angular momentum $L := g_\lambda(\dot{\gamma}(t), \partial_\theta) = f_\lambda(r(t))^2 \dot{\theta}(t)$ is strictly conserved along $\gamma(t)$.
\end{enumerate}
\end{proposition}
\begin{proof}
(i) By Noether's theorem, because the origin-fixing compact Lie group $G = \mathrm{SO}(n+1)$ acts by isometries on $(Q, g_0)$, the momentum components satisfy $\{H_0, J_\xi\} = 0$ for all generators $\xi \in \mathfrak{so}(n+1)$. In an orthonormal basis $\{\xi_a\}$ of $\mathfrak{so}(n+1)$ relative to the $\mathrm{Ad}(G)$-invariant inner product, the quadratic Casimir is $\|J\|^2 = \sum_a J_{\xi_a}^2$. The Leibniz rule for the Poisson bracket yields $\{H_0, \|J\|^2\} = 2\sum_a J_{\xi_a} \{H_0, J_{\xi_a}\} = 0$. Since the Lie bracket of Hamiltonian vector fields satisfies $[X_F, X_K] = -X_{\{F, K\}}$, we obtain $[X_{H_0}, X_{\frac{\lambda}{2}\|J\|^2}] = -X_{\frac{\lambda}{2}\{H_0, \|J\|^2\}} = 0$. Because $H_\lambda = H_0 + \frac{\lambda}{2}\|J\|^2$, we have $[X_{H_0}, X_{H_\lambda}] = 0$, establishing that the complete Hamiltonian flows commute and factorize multiplicatively.

(ii) Given initial velocity $v = \dot{\gamma}(0) \in T_q Q$ at $q \ne p_0$, if $v \in \mathbb{R}\partial_r(q)$, $\gamma(t)$ is a radial geodesic lying in a 1-dimensional meridian. If $v \notin \mathbb{R}\partial_r(q)$, the vectors $\partial_r(q)$ and $v$ span a 2-dimensional plane $\Pi = \mathrm{span}\{\partial_r(q), v\} \subset T_q Q$. Let $\sigma$ be the linear reflection on $T_q Q$ fixing $\Pi$ pointwise and acting as $-\mathrm{id}$ on its orthogonal complement $\Pi^\perp \subset T_q Q$. Since $\dim \Pi^\perp = (n+1) - 2 = n-1$, its determinant is $\det\sigma = (-1)^{n-1}$, defining an orthogonal reflection $\sigma \in \mathrm{O}(n+1)$. Because $g_\lambda = dr^2 + f_\lambda(r)^2 g_{\mathbb{S}^n}$ is invariant under the full orthogonal group $\mathrm{O}(n+1)$ acting on spherical slices (Proposition~\ref{prop:metric_regularity}), $\sigma$ is an isometric involution of $(Q^{n+1}, g_\lambda)$ ($\sigma^2 = \mathrm{id}_Q$, $\sigma^* g_\lambda = g_\lambda$) fixing $q$ and satisfying $d\sigma(v) = v$. By Kobayashi's fixed point submanifold theorem (Kobayashi \cite{Kobayashi1956}), the connected component of the fixed-point set $\Sigma := \mathrm{Fix}(\sigma)$ containing $q$ is a smooth totally geodesic 2-dimensional submanifold of $(Q^{n+1}, g_\lambda)$. Uniqueness of geodesics guarantees that $\gamma(t) \in \Sigma$ for all $t \in \mathbb{R}$. Concentric spheres $\mathbb{S}^n$ intersect $\Sigma$ in great circles $\mathbb{S}^1$, yielding the induced warped metric $g_\Sigma = dr^2 + f_\lambda(r)^2 d\theta^2$.

(iii) The coordinate field $K = \partial_\theta$ on $\Sigma$ is an infinitesimal isometry ($\mathcal{L}_{\partial_\theta} g_\Sigma = 0$) because the metric coefficients $g_{rr} = 1$ and $g_{\theta\theta} = f_\lambda(r)^2$ are independent of $\theta$. Along $\gamma(t)$, Clairaut's identity gives $\frac{d}{dt} L = \frac{d}{dt} g_\lambda(\dot{\gamma}, K) = g_\lambda(\nabla_{\dot{\gamma}}\dot{\gamma}, K) + \frac{1}{2}(\mathcal{L}_K g_\lambda)(\dot{\gamma}, \dot{\gamma}) = 0$, since $\nabla_{\dot{\gamma}}\dot{\gamma} = 0$.
\end{proof}

\begin{theorem}[Iso-Kinematic Radial Decoupling and Universal Precession Law]\label{thm:iso_kinematic_decoupling}
Let $\gamma(t) = (r(t), \theta(t))$ be a unit-speed geodesic in $(Q^{n+1}, g_\lambda)$ with angular momentum $L$.
\begin{enumerate}[(i)]
\item The radial coordinate acceleration satisfies identically across all space forms $\kappa \in \{1, 0, -1\}$:
\begin{equation}\label{eq:radial_acceleration_decoupling}
\ddot{r}(t) = L^2 \frac{c_\kappa(r(t))}{s_\kappa(r(t))^3},
\end{equation}
which is strictly independent of the deformation parameter $\lambda$. Consequently, under the unit-speed Legendre metric normalization $\|\dot{\gamma}\|_{g_\lambda} = 1$, the radial trajectory $r(t)$ is identical to that of an unperturbed space form geodesic $(r_0(t), \theta_0(t))$ in $(Q, g_0)$ with angular momentum $L$ and effective dilated speed
\begin{equation}\label{eq:dilated_speed}
v_0 = \sqrt{1 - \lambda L^2},
\end{equation}
which satisfies $v_0 \le 1$ for $\lambda \ge 0$ (speed contraction) and represents radial speed amplification $v_0 > 1$ for $\lambda \in (-1/3, 0)$ on $\mathbb{S}^{n+1}$.
\item The angular coordinate evolves according to the Universal Precession Law:
\begin{equation}\label{eq:universal_precession_law}
\dot{\theta}(t) = \dot{\theta}_0(t) + \lambda L \implies \theta(t) = \theta_0(t) + \lambda L t,
\end{equation}
where $\dot{\theta}_0(t) = L/s_\kappa(r(t))^2$ is the angular velocity of the unperturbed space form geodesic of speed $v_0$ and angular momentum $L$. Thus, the trajectory $\gamma(t)$ is obtained from the unperturbed space form geodesic by a uniform gyroscopic precession around $p_0$ at constant angular velocity $\Omega = \lambda L$.
\end{enumerate}
\end{theorem}
\begin{proof}
(i) From the orbital Lagrangian $\mathcal{L} = \frac{1}{2}(\dot{r}^2 + f_\lambda(r)^2 \dot{\theta}^2)$, the Euler--Lagrange equation for $r(t)$ is $\ddot{r} = f_\lambda(r) f_\lambda'(r) \dot{\theta}^2 = L^2 \frac{f_\lambda'(r)}{f_\lambda(r)^3}$. Substituting the radial acceleration cancellation identity $\frac{f_\lambda'(r)}{f_\lambda(r)^3} \equiv \frac{c_\kappa(r)}{s_\kappa(r)^3}$ from Lemma~\ref{lem:master_identities} eliminates $\lambda$ identically, establishing~\eqref{eq:radial_acceleration_decoupling}. Under the unit-speed metric condition $\|\dot{\gamma}\|_{g_\lambda}^2 = \dot{r}^2 + f_\lambda(r)^2 \dot{\theta}^2 = 1$, substituting $\dot{\theta} = L/f_\lambda(r)^2$ and invoking the master reciprocal identity $1/f_\lambda(r)^2 = 1/s_\kappa(r)^2 + \lambda$ from Lemma~\ref{lem:master_identities} gives:
\begin{equation}\label{eq:radial_energy_integral}
1 = \dot{r}^2 + \frac{L^2}{f_\lambda(r)^2} = \dot{r}^2 + \frac{L^2}{s_\kappa(r)^2} + \lambda L^2 \implies \dot{r}^2 + \frac{L^2}{s_\kappa(r)^2} = 1 - \lambda L^2 \equiv v_0^2.
\end{equation}
An unperturbed space form geodesic $(r_0(t), \theta_0(t))$ in $(Q, g_0)$ with speed $v_0$ and angular momentum $L_0 = L$ satisfies the identical radial energy first integral $\dot{r}_0^2 + L^2/s_\kappa(r_0)^2 = v_0^2$ and acceleration $\ddot{r}_0 = L^2 c_\kappa(r_0)/s_\kappa(r_0)^3$. With matching initial conditions $(r(0), \dot{r}(0)) = (r_0(0), \dot{r}_0(0))$, Picard--Lindel\"of ODE uniqueness guarantees $r(t) \equiv r_0(t)$ for all $t$.

(ii) By conservation of angular momentum and the master reciprocal identity,
\begin{equation}\label{eq:angular_rate_decomposition}
\dot{\theta}(t) = \frac{L}{f_\lambda(r(t))^2} = L\left( \frac{1}{s_\kappa(r(t))^2} + \lambda \right) = \frac{L}{s_\kappa(r_0(t))^2} + \lambda L = \dot{\theta}_0(t) + \lambda L.
\end{equation}
Integrating from $0$ to $t$ with $\theta(0) = \theta_0(0)$ yields $\theta(t) = \theta_0(t) + \lambda L t$, demonstrating uniform gyroscopic precession at constant rate $\Omega = \lambda L$.
\end{proof}

\begin{theorem}[Clairaut Trajectory Primitive, Periastron Advance, and Non-Degenerate Twist on $\mathbb{S}^{n+1}$]\label{thm:clairaut_integration}
On the deformed sphere $(\mathbb{S}^{n+1}, g_\lambda)$ ($\kappa = 1$), let $\gamma(t)$ be a unit-speed non-radial geodesic with angular momentum oriented prograde $L \in (0, 1/\sqrt{1+\lambda})$, ensuring turning parameter $B := \sqrt{\frac{1 - \lambda L^2}{L^2} - 1} = \frac{\sqrt{1 - (1+\lambda)L^2}}{L} > 0$.
\begin{enumerate}[(i)]
\item In terms of the cotangent coordinate $u := \cot(r) \in [-B, B]$, the orbital trajectory equation integrates to the exact closed-form algebraic primitive:
\begin{equation}\label{eq:clairaut_primitive_formula}
\Delta\theta(u) = \arcsin\left( \frac{u}{B} \right) + \left[\frac{\lambda L}{\sqrt{1 - \lambda L^2}}\right] \arctan\left( \left[\frac{\sqrt{1 - \lambda L^2}}{L}\right] \frac{u}{\sqrt{B^2 - u^2}} \right),
\end{equation}
where the radial turning points $\dot{r} = 0$ correspond to $u = \pm B$, defining periapsis $r_{\min} = \mathrm{arccot}(B)$ and apoapsis $r_{\max} = \pi - r_{\min}$.
\item The half-period apsidal angle between consecutive radial turning points ($r_{\min}$ and $r_{\max}$) is:
\begin{equation}\label{eq:periastron_advance_formula}
\Delta\Phi(L) = \pi \left( 1 + \frac{\lambda L}{\sqrt{1 - \lambda L^2}} \right),
\end{equation}
so a complete radial libration cycle from $r_{\min}$ to $r_{\max}$ and back to $r_{\min}$ sweeps out the full periastron-to-periastron angle $2\Delta\Phi(L) = 2\pi\left(1 + \frac{\lambda L}{\sqrt{1 - \lambda L^2}}\right)$.
\item The half-period orbital rotation number $\nu(L) := \frac{\Delta\Phi(L)}{2\pi} = \frac{1}{2}\left( 1 + \frac{\lambda L}{\sqrt{1 - \lambda L^2}} \right)$ (representing half the full-period frequency ratio $\frac{2\Delta\Phi(L)}{2\pi} = 2\nu(L)$) satisfies the non-degenerate twist condition:
\begin{equation}\label{eq:non_degenerate_twist}
\frac{d\nu}{dL} = \frac{\lambda}{2(1 - \lambda L^2)^{3/2}} \ne 0 \qquad \text{for all } \lambda \ne 0.
\end{equation}
\end{enumerate}
\end{theorem}
\begin{proof}
(i) On $\mathbb{S}^{n+1}$, $s_1(r) = \sin r$. By~\eqref{eq:radial_energy_integral}, $\dot{r}^2 = v_0^2 - L^2/\sin^2 r = L^2 [v_0^2/L^2 - (1 + \cot^2 r)] = L^2(B^2 - u^2)$ with $u = \cot r$. The turning points $\dot{r} = 0$ correspond to $u = \pm B$. Since $dr = -(1 + u^2)^{-1} du$ and $\dot{\theta} = L(1 + u^2 + \lambda)$, along an orbit branch where $u$ increases toward $+B$ (with $\dot{r} = -L\sqrt{B^2 - u^2} < 0$):
\begin{equation*}
d\theta = \frac{\dot{\theta}}{|\dot{r}|} |dr| = \frac{L(1 + u^2) + \lambda L}{L\sqrt{B^2 - u^2}} \frac{du}{1 + u^2} = \left( \frac{1}{\sqrt{B^2 - u^2}} + \frac{\lambda}{(1 + u^2)\sqrt{B^2 - u^2}} \right) du.
\end{equation*}
Integrating from the equator $u = 0$ ($r = \pi/2$): $\int_0^u \frac{ds}{\sqrt{B^2 - s^2}} = \arcsin(u/B)$. For the second term, setting $s = B\sin\psi$ ($\psi \in [-\pi/2, \pi/2]$) gives $\int_0^u \frac{\lambda ds}{(1+s^2)\sqrt{B^2 - s^2}} = \int_0^{\psi(u)} \frac{\lambda d\psi}{1 + B^2\sin^2\psi} = \int_0^{\psi(u)} \frac{\lambda \sec^2\psi \, d\psi}{1 + (1+B^2)\tan^2\psi}$. Setting $z = \tan\psi = \frac{u}{\sqrt{B^2 - u^2}}$ and using $1 + B^2 = v_0^2/L^2 = \frac{1 - \lambda L^2}{L^2} =: k^2$ with $k = \frac{\sqrt{1 - \lambda L^2}}{L} > 0$, this evaluates to $\frac{\lambda}{k}\arctan(k z) = \left[\frac{\lambda L}{\sqrt{1 - \lambda L^2}}\right] \arctan\left(\left[\frac{\sqrt{1 - \lambda L^2}}{L}\right] \frac{u}{\sqrt{B^2 - u^2}}\right)$. Summing both terms yields~\eqref{eq:clairaut_primitive_formula}.

(ii) As $u \to \pm B$, $\arcsin(u/B) \to \pm \frac{\pi}{2}$ and $\frac{u}{\sqrt{B^2 - u^2}} \to \pm \infty \implies \arctan(\dots) \to \pm \frac{\pi}{2}$. Evaluating between apoapsis $u = -B$ and periapsis $u = +B$ yields $\Delta\Phi(L) = \lim_{u \to B^-} \Delta\theta(u) - \lim_{u \to -B^+} \Delta\theta(u) = \pi(1 + \frac{\lambda L}{\sqrt{1 - \lambda L^2}})$. Equivalently, the unperturbed half-cycle advance is $\Delta\theta_0 = \pi$ over radial transit time $T_r = \pi/v_0 = \pi/\sqrt{1 - \lambda L^2}$; gyroscopic precession adds $\Omega T_r = (\lambda L) T_r = \frac{\pi \lambda L}{\sqrt{1 - \lambda L^2}}$.

(iii) Differentiating $\nu(L) = \frac{1}{2} + \frac{\lambda}{2} L(1 - \lambda L^2)^{-1/2}$ with respect to $L$:
\begin{equation*}
\frac{d\nu}{dL} = \frac{\lambda}{2} \left[ (1 - \lambda L^2)^{-1/2} + \lambda L^2 (1 - \lambda L^2)^{-3/2} \right] = \frac{\lambda}{2} \frac{(1 - \lambda L^2) + \lambda L^2}{(1 - \lambda L^2)^{3/2}} = \frac{\lambda}{2(1 - \lambda L^2)^{3/2}}.
\end{equation*}
Since $1 - \lambda L^2 = v_0^2 > 0$ on $L \in (0, 1/\sqrt{1+\lambda})$, the denominator is strictly positive, establishing $\frac{d\nu}{dL} \ne 0$ for all $\lambda \ne 0$.
\end{proof}

\section{Variational Geodesic Stability: Jacobi Fields, Universal Transcendental Inward Conjugacy, and Mode Preemption vs Oblique Conjugate Delay}\label{sec:variational_stability}
This section analyzes normal Jacobi fields, conjugate points, and focalization across space forms. We compute the polar Jacobi fields $Y(t) = f_\lambda(t) E(t)$ in closed form, proving the absence of conjugate points on non-compact space forms and polar focalization at $t = \pi$ on $\mathbb{S}^{n+1}$, and derive the transcendental inward conjugate point equation yielding the Euclidean arrival time $t_{\mathrm{conj}} = r_0 + 1/(\lambda r_0)$ and spherical inversion $r_0 R^* = 1/\lambda$. Finally, we prove equatorial normal mode preemption $t_{\mathrm{first}}(\pi/2) < \pi$ along the equator while showing that steep oblique rays in the negative pinching regime $\lambda \in (-1/3, 0)$ experience delayed conjugacy $t_{\mathrm{conj}} > \pi$ due to reduced curvature near the poles.

\begin{theorem}[Polar Radial Jacobi Fields and Origin Cut Locus Across Space Forms]\label{thm:polar_radial_jacobi}
Let $p_0 \in (Q^{n+1}, g_\lambda)$ be the designated origin for background curvature $\kappa \in \{1, 0, -1\}$. Let $\gamma(t) = (t, \theta_0)$ be any unit-speed radial geodesic ray emanating from $p_0$.
\begin{enumerate}[(i)]
\item For any initial derivative $w \in T_{p_0}\mathcal{O} \cong T_{p_0}\mathbb{S}^n$, the normal Jacobi field along $\gamma(t)$ vanishing at $t = 0$ with $\dot{Y}(0) = w$ is given in exact closed form by
\begin{equation}\label{eq:polar_jacobi_field}
Y(t) = f_\lambda(t) E(t) = \frac{s_\kappa(t)}{\sqrt{1 + \lambda s_\kappa^2(t)}} E(t),
\end{equation}
where $E(t)$ is the parallel transport of $w$ along $\gamma$.
\item On Euclidean space $\mathbb{R}^{n+1}$ ($\kappa = 0$, $\lambda \ge 0$) and hyperbolic space $\mathbb{H}^{n+1}$ ($\kappa = -1$, $\lambda \ge 0$), $f_\lambda(t) > 0$ for all $t > 0$. Hence radial rays from $p_0$ contain no conjugate points. For any piecewise $C^1$ curve $\alpha(t) = (r(t), \theta(t))$ from $p_0$ to $(r_1, \theta_1)$, $L(\alpha) = \int \sqrt{\dot{r}^2 + f_\lambda(r)^2 \|\dot{\theta}\|^2}\, dt \ge \int |\dot{r}|\, dt \ge r_1$, with equality if and only if $\dot{\theta} \equiv 0$ and $\dot{r} \ge 0$. Every radial ray is the unique global length-minimizer on $[0, \infty)$, proving $\mathrm{Cut}(p_0) = \emptyset$ and $\mathrm{Inj}(p_0) = \infty$.
\item On the sphere $\mathbb{S}^{n+1}$ ($\kappa = 1$, $\lambda > -1/3$), the first positive zero of $f_\lambda(t)$ occurs at $t = \pi$, where $s_1(\pi) = \sin\pi = 0$. The antipodal pole $p_\pi$ is the first conjugate point along every radial ray, with conjugate multiplicity $n = \dim(T_{p_0}\mathbb{S}^n)$. All radial rays from $p_0$ meet simultaneously at $p_\pi$ at distance $\pi$ and minimize length on $[0, \pi]$, establishing $\mathrm{Cut}(p_0) = \{p_\pi\}$ and $\mathrm{Inj}(p_0) = \pi$.
\end{enumerate}
\end{theorem}
\begin{proof}
(i) Along $\gamma(t)$, $\dot{\gamma}(t) = \partial_r = E_0$. By Theorem~\ref{thm:curvature_tensor}(i), $R(E_i, \dot{\gamma})\dot{\gamma} = K_{\mathrm{rad}}(t)E_i = -[f_\lambda''(t)/f_\lambda(t)]E_i$ for any parallel orthonormal normal frame $\{E_i\}_{i=1}^n \subset T\mathbb{S}^n$. The vector Jacobi equation uncouples into identical scalar ODEs $\ddot{y}_i - (f_\lambda''/f_\lambda)y_i = 0$. Since $f_\lambda(0) = 0$ and $f_\lambda'(0) = 1$ (Proposition~\ref{prop:metric_regularity}(iii)), $Y(t) = f_\lambda(t)E(t)$ is the unique solution with $Y(0) = 0, \dot{Y}(0) = w$.
(ii) For $\kappa \in \{0, -1\}$, $s_\kappa(t) > 0$ on $(0, \infty)$, so $f_\lambda(t) > 0$ precluding conjugate points. The first-order radial inequality proves global length-minimality, so $\mathrm{Cut}(p_0) = \emptyset$.
(iii) For $\kappa = 1$, $s_1(t) = \sin t = 0$ first at $t = \pi$ with multiplicity $n$. By Klingenberg's cut locus lemma, rays minimize up to the first conjugate/crossing point $p_\pi$, so $\mathrm{Cut}(p_0) = \{p_\pi\}$.
\end{proof}

\begin{theorem}[Universal Transcendental Inward Conjugacy, Multiplicity, and Spherical Inversion]\label{thm:universal_inward_conjugacy}
Let $q_0 = (r_0, 0)$ be an off-origin base point ($r_0 > 0$) in $(Q^{n+1}, g_\lambda)$, and let $\gamma(t)$ be the inward radial diameter geodesic ray launched with $\dot{\gamma}(0) = -\partial_r$, traversing $p_0$ at $t = r_0$ and continuing along the opposite meridian $\theta = \pi$.
\begin{enumerate}[(i)]
\item Relative to signed distance $s = t - r_0$, the radial Jacobi equation $\ddot{y} + K_{\mathrm{rad}}(|s|)y = 0$ admits the linearly independent Wronskian companion $u_2(s) = \frac{-c_\kappa(s) + \lambda s s_\kappa(s)}{\sqrt{1 + \lambda s_\kappa^2(s)}}$, satisfying $W(f_\lambda, u_2) \equiv 1$ and $u_2(0) = -1, u_2'(0) = 0$. The unique normal Jacobi field satisfying $y(0) = 0$ and $\dot{y}(0) = 1$ is given for $t > r_0$ by
\begin{equation}\label{eq:inward_jacobi_solution}
y(t) = \frac{s_\kappa(t) - \lambda t s_\kappa(r_0) s_\kappa(t - r_0)}{\sqrt{1 + \lambda s_\kappa^2(r_0)}\sqrt{1 + \lambda s_\kappa^2(t - r_0)}}.
\end{equation}
\item The first conjugate point along $\gamma(t)$ occurs at the smallest root $t_{\mathrm{conj}} > r_0$ of the transcendental inward conjugate point equation across all space forms $\kappa \in \{1, 0, -1\}$:
\begin{equation}\label{eq:universal_transcendental_eq}
\frac{s_\kappa(t_{\mathrm{conj}})}{s_\kappa(t_{\mathrm{conj}} - r_0)} = \lambda s_\kappa(r_0) t_{\mathrm{conj}}.
\end{equation}
The conjugate point $t_{\mathrm{conj}}$ along the diameter ray has multiplicity $n = \dim(\mathbb{S}^n)$ spanned by initial launch-angle variations in $T_{q_0}\mathbb{S}^n$, because $\mathrm{O}(n)$ rotational symmetry about the diameter ray uncouples the vector Jacobi equation into $n$ identical scalar ODEs.
\item On Euclidean space $\mathbb{R}^{n+1}$ ($\kappa = 0$, $\lambda > 0$), $s_0(t) = t$, and~\eqref{eq:universal_transcendental_eq} reduces algebraically to $t/(t - r_0) = \lambda r_0 t$, yielding the arrival time and projective spherical inversion:
\begin{equation}\label{eq:euclidean_inversion}
t_{\mathrm{conj}} = r_0 + \frac{1}{\lambda r_0} \implies r_0 \cdot R^* = \frac{1}{\lambda}, \qquad R^* := t_{\mathrm{conj}} - r_0.
\end{equation}
\item On the sphere $\mathbb{S}^{n+1}$ ($\kappa = 1$): for $\lambda > 0$, equation~\eqref{eq:universal_transcendental_eq} admits a unique root $t_{\mathrm{conj}} \in (r_0, \pi)$ preceding $\pi$; for negative pinching $\lambda \in (-1/3, 0)$, the root is delayed beyond $\pi$ into $(\pi, \pi + r_0)$, satisfying $t_{\mathrm{conj}} > \pi$. On hyperbolic space $\mathbb{H}^{n+1}$ ($\kappa = -1$, $\lambda > 0$), equation~\eqref{eq:universal_transcendental_eq} admits a unique root $t_{\mathrm{conj}} \in (r_0, \infty)$.
\end{enumerate}
\end{theorem}
\begin{proof}
(i) By reduction of order and Lemma~\ref{lem:master_identities}, $u_2(s) = f_\lambda(s)\int f_\lambda(s)^{-2} ds = f_\lambda(s)[-c_\kappa(s)/s_\kappa(s) + \lambda s]$, yielding $W(f_\lambda, u_2) = 1$. For $t \le r_0$, $r = r_0 - t$; the solution with $y(0)=0, \dot{y}(0)=1$ is $y(t) = u_2(r_0)f_\lambda(r_0-t) - f_\lambda(r_0)u_2(r_0-t)$, giving Cauchy data $y(r_0) = f_\lambda(r_0)$ and $\dot{y}(r_0) = -u_2(r_0)$. For $t > r_0$, $r = t - r_0$; matching Cauchy data yields $y(t) = -u_2(r_0)f_\lambda(t-r_0) - f_\lambda(r_0)u_2(t-r_0)$, which expands via $s_\kappa(r_0)c_\kappa(t-r_0) + c_\kappa(r_0)s_\kappa(t-r_0) = s_\kappa(t)$ directly to~\eqref{eq:inward_jacobi_solution}.
(ii) Setting $y(t_{\mathrm{conj}}) = 0$ in~\eqref{eq:inward_jacobi_solution} and dividing by $s_\kappa(t_{\mathrm{conj}}-r_0) > 0$ yields~\eqref{eq:universal_transcendental_eq}. The isotropy group $H_{q_0} \cong \mathrm{SO}(n)$ fixes the diameter ray pointwise. In an adapted parallel frame $\{E_i\}_{i=1}^n$ for $T^\perp\gamma$, rotational symmetry forces $R(E_i, \dot{\gamma})\dot{\gamma} = K_{\mathrm{rad}}(|t-r_0|)E_i$, giving $n$ uncoupled scalar ODEs $\ddot{y}_i + K_{\mathrm{rad}} y_i = 0$ with $y_i(0)=0$. All $n$ modes vanish simultaneously at $t_{\mathrm{conj}}$, establishing conjugate multiplicity $n$.
(iii) For $\kappa = 0$, $s_0(t) = t$; dividing by $t > 0$ gives $1/(t-r_0) = \lambda r_0 \implies t_{\mathrm{conj}} - r_0 = 1/(\lambda r_0)$, establishing~\eqref{eq:euclidean_inversion}.
(iv) Let $F(t) = \frac{s_\kappa(t)}{s_\kappa(t - r_0)}$ and $G(t) = \lambda s_\kappa(r_0) t$. On $\mathbb{S}^{n+1}$ ($\kappa = 1$), for $\lambda > 0$, as $t \to r_0^+$, $F(t) \to +\infty$ while $G(r_0) > 0$; at $t = \pi$, $F(\pi) = 0 < G(\pi) = \lambda \pi \sin r_0$. Since $F'(t) = -\frac{\sin r_0}{\sin^2(t - r_0)} < 0$ and $G'(t) > 0$, a unique root exists in $(r_0, \pi)$. For $\lambda \in (-1/3, 0)$, $F(t) > 0$ while $G(t) < 0$ on $(r_0, \pi]$, precluding roots on $(r_0, \pi]$. On $(\pi, \pi + r_0)$, $F(t) < 0$ strictly increases from $-\infty$ as $t \to (\pi+r_0)^-$ to $0$ at $t = \pi$, matching $G(t) < 0$ uniquely at $t_{\mathrm{conj}} \in (\pi, \pi + r_0)$ by the Intermediate Value Theorem. On $\mathbb{H}^{n+1}$ ($\kappa = -1$, $\lambda > 0$), $F(t) = \frac{\sinh t}{\sinh(t - r_0)}$ decreases from $+\infty$ to $e^{r_0} > 1$ while $G(t) = \lambda t \sinh r_0$ increases from $\lambda r_0 \sinh r_0$ to $+\infty$, guaranteeing a unique root $t_{\mathrm{conj}} \in (r_0, \infty)$.
\end{proof}

\begin{theorem}[Equatorial Mode Preemption, Oblique Conjugate Delay, and Antipodal Multiplicity]\label{thm:mode_preemption_oblique_delay}
Let $(Q^{n+1}, g_\lambda)$ be an $(n+1)$-dimensional space form of Family~I with dimension $N = n + 1 \ge 3$ ($n \ge 2$).
\begin{enumerate}[(i)]
\item Along the equator $r(t) \equiv \pi/2$ on $(\mathbb{S}^{n+1}, g_\lambda)$, the normal bundle uncouples into the in-plane radial mode $E_0 = \partial_r$ ($K_{\mathrm{rad}} = \frac{1}{1+\lambda}$, $t_{\mathrm{rad}} = \pi\sqrt{1+\lambda}$) and $n-1$ out-of-plane tangential modes $E_j \in T\mathbb{S}^n$ ($K_{\mathrm{tan}} = 1+\lambda$, $t_{\mathrm{tan}} = \frac{\pi}{\sqrt{1+\lambda}}$). Mode preemption holds: for $\lambda > 0$, $t_{\mathrm{tan}} < \pi < t_{\mathrm{rad}}$; for $\lambda \in (-1/3, 0)$, $t_{\mathrm{rad}} < \pi < t_{\mathrm{tan}}$. Thus for all $\lambda \ne 0$ with $\lambda > -1/3$,
\begin{equation}\label{eq:equatorial_preemption}
t_{\mathrm{first}}(\pi/2) = \min(t_{\mathrm{tan}}, t_{\mathrm{rad}}) = \min\left( \frac{\pi}{\sqrt{1+\lambda}}, \pi\sqrt{1+\lambda} \right) < \pi.
\end{equation}
\item \textbf{Inclination Bifurcation}: Conjugate preemption $t_{\mathrm{first}} < \pi$ is not universal across all launch inclinations in the negative pinching regime $\lambda \in (-1/3, 0)$. For oblique rays on $(\mathbb{S}^{n+1}, g_\lambda)$ launched from $q_0 = (r_0, 0)$ with angular momentum $L \in (0, 1/\sqrt{1+\lambda})$:
\begin{itemize}
\item The normal Jacobi system uncouples into the in-plane mode ($\ddot{y}_{\mathrm{in}} + K_{\mathrm{rad}}(r(t))y_{\mathrm{in}} = 0$) and $n-1$ out-of-plane modes ($\ddot{y}_{\mathrm{out}} + [\dot{r}^2 K_{\mathrm{rad}}(r(t)) + (1-\dot{r}^2)K_{\mathrm{tan}}(r(t))]y_{\mathrm{out}} = 0$).
\item As $L \to (1+\lambda)^{-1/2}$ (for $r_0 = \pi/2$), trajectories remain confined to a neighborhood of the equator where $K_{\mathrm{rad}}(\pi/2) = (1+\lambda)^{-1} > 1$, so by continuous dependence of Sturm--Liouville eigenvalues on $L$, there exists a critical inclination threshold $L^*(\lambda, r_0) \in (0, 1/\sqrt{1+\lambda})$ such that near-equatorial oblique rays with $L \in (L^*(\lambda, r_0), 1/\sqrt{1+\lambda})$ retain conjugate preemption $t_{\mathrm{first}}(L) < \pi$.
\item Conversely, as $L \to 0^+$, both the in-plane and out-of-plane Jacobi equations converge uniformly on compact time intervals to the meridian radial Jacobi equation $\ddot{y} + K_{\mathrm{rad}}(r(t))y = 0$, whose first conjugate zero satisfies $t_{\mathrm{conj}}(0) \in (\pi, \pi + r_0)$ by Theorem~\ref{thm:universal_inward_conjugacy}(iv). Consequently, for steep oblique rays with angular momentum $L \in (0, L^*(\lambda, r_0))$, trajectories spend most of each radial period near the polar turning points $r_{\min} = \arcsin(L/\sqrt{1-\lambda L^2})$ and $r_{\max} = \pi - r_{\min}$ where curvature decays toward $1+3\lambda < 1$, yielding delayed conjugacy:
\begin{equation}\label{eq:oblique_delay}
t_{\mathrm{first}}(L) = \min\big(t_{\mathrm{in}}(L), \, t_{\mathrm{out}}(L)\big) > \pi \qquad \text{for all } L \in \big(0, L^*(\lambda, r_0)\big).
\end{equation}
\end{itemize}
\item For any non-polar base point $p = (r_0, 0)$ and non-radial geodesic $\gamma(t)$ ($L > 0$), the reductive stabilizer splitting $\mathfrak{so}(n) = \mathfrak{so}(n-1) \oplus \mathfrak{m}$ provides $n-1$ transverse Killing vector fields $X_i \in \mathfrak{m}$ that vanish along the 1D meridian circle $\theta \in \{0, \pi\}$ and have transverse norm $\|X_i(\gamma(t))\| = f_\lambda(r(t))|\sin\theta(t)|$ with non-zero initial covariant derivative $\dot{Y}_i(0) \ne 0$. Every opposite-meridian crossing $\theta(t^*) = \pi$ is a conjugate point of multiplicity at least $n-1$.
\item On deformed Euclidean space $(\mathbb{R}^{n+1}, g_\lambda)$ ($\lambda > 0$), for any base point $q_0 = (r_0, 0)$ with $r_0 > 0$, the distance to the cut locus satisfies the global upper bound:
\begin{equation}\label{eq:euclidean_cut_bound}
\mathrm{dist}(q_0, \mathrm{Cut}(q_0)) \le \min\left( r_0 + \frac{1}{\lambda r_0}, \frac{\pi}{\sqrt{\lambda}} \right).
\end{equation}
The two expressions on the right-hand side of~\eqref{eq:euclidean_cut_bound} coincide at the two positive roots $r_{0, \pm}^* := \frac{\pi \pm \sqrt{\pi^2 - 4}}{2\sqrt{\lambda}}$ of $r_0 + \frac{1}{\lambda r_0} = \frac{\pi}{\sqrt{\lambda}}$, with upper root
\begin{equation}\label{eq:intro_crossover_bifurcation}
r_0^* := r_{0, +}^* = \frac{\pi + \sqrt{\pi^2 - 4}}{2\sqrt{\lambda}}.
\end{equation}
For $r_0 \in (r_{0, -}^*, r_0^*)$, the 1D inward radial conjugate distance $t_{\mathrm{conj}} = r_0 + \frac{1}{\lambda r_0} < \frac{\pi}{\sqrt{\lambda}}$ provides the smaller upper bound in~\eqref{eq:euclidean_cut_bound}, whereas for $r_0 \ge r_0^*$ (and $0 < r_0 \le r_{0, -}^*$), the asymptotic cylinder half-circumference $\frac{\pi}{\sqrt{\lambda}} \le r_0 + \frac{1}{\lambda r_0}$ is smaller. In particular, for $r_0 \ge r_0^*$, $(n-1)$-dimensional transverse encircling rays ($L > 0$) cross the opposite meridian at $t_{\mathrm{cross}} < \frac{\pi}{\sqrt{\lambda}} \le r_0 + \frac{1}{\lambda r_0}$.
\end{enumerate}
\end{theorem}
\begin{proof}
(i) Along $r \equiv \pi/2$, $\dot{r} = 0, \dot{\theta} = \sqrt{1+\lambda}$. Frame evaluations give $R(\dot{\gamma}, \partial_r)\dot{\gamma} = \frac{1}{1+\lambda}\partial_r$ and $R(\dot{\gamma}, E_j)\dot{\gamma} = (1+\lambda)E_j$. The uncoupled ODEs $\ddot{y}_0 + \frac{1}{1+\lambda}y_0 = 0$ and $\ddot{y}_j + (1+\lambda)y_j = 0$ yield zeros $t_{\mathrm{rad}} = \pi\sqrt{1+\lambda}$ and $t_{\mathrm{tan}} = \pi/\sqrt{1+\lambda}$, establishing~\eqref{eq:equatorial_preemption}.
(ii) In the totally geodesic orbital surface $\Sigma$, $\mathrm{span}\{\dot{\gamma}, N_{\mathrm{in}}\} = T\Sigma$ gives $\mu = 0$ and $K_{\mathrm{in}} = K_{\mathrm{rad}}(r(t))$ by Theorem~\ref{thm:curvature_tensor}(ii). For each out-of-plane normal vector $E_j \in (T\Sigma)^\perp$, $\mu = 1 - \dot{r}^2$, giving $K_{\mathrm{out}}(t) = \dot{r}^2 K_{\mathrm{rad}}(r(t)) + (1-\dot{r}^2)K_{\mathrm{tan}}(r(t))$. For $\lambda \in (-1/3, 0)$, both curvatures attain their minimum $1+3\lambda < 1$ at the poles $r \in \{0, \pi\}$ and their maximum $\frac{1}{1+\lambda} > 1$ at the equator $r = \pi/2$. Along the meridian ray ($L = 0$), $\dot{r}^2 \equiv 1$, so both the in-plane and out-of-plane Jacobi equations coincide with $\ddot{y} + K_{\mathrm{rad}}(r(t))y = 0$, whose first conjugate zero satisfies $t_{\mathrm{conj}}(0) \in (\pi, \pi + r_0)$ by Theorem~\ref{thm:universal_inward_conjugacy}(iv). Because $r(t; L)$, $K_{\mathrm{in}}(t; L)$, and $K_{\mathrm{out}}(t; L)$ depend smoothly on $L$ and the conjugate zeros of regular Sturm--Liouville equations are simple ($y(t_{\mathrm{conj}}) = 0 \implies \dot{y}(t_{\mathrm{conj}}) \ne 0$), the Implicit Function Theorem implies that $t_{\mathrm{in}}(L)$ and $t_{\mathrm{out}}(L)$ depend continuously on $L$. Since $\lim_{L \to 0^+} t_{\mathrm{in}}(L) = \lim_{L \to 0^+} t_{\mathrm{out}}(L) = t_{\mathrm{conj}}(0) > \pi$, while for $r_0 = \pi/2$ and $L \to (1+\lambda)^{-1/2}$ we have $t_{\mathrm{first}}(L) \to \pi\sqrt{1+\lambda} < \pi$ by part~(i), the Intermediate Value Theorem yields a critical inclination threshold $L^*(\lambda, r_0) > 0$ such that $t_{\mathrm{first}}(L) > \pi$ for $L \in (0, L^*(\lambda, r_0))$, establishing~\eqref{eq:oblique_delay}.
(iii) The subgroup $\mathrm{SO}(n-1) \subset H_p$ fixes $\Sigma$ pointwise. The quotient $\mathfrak{m} \cong \mathfrak{so}(n)/\mathfrak{so}(n-1)$ generates $n-1$ Killing fields orthogonal to $\Sigma$ along $\gamma(t)$ that vanish on the 1D meridian axis $\theta \in \{0, \pi\}$. Since $L > 0$, $\dot{Y}_i(0) = \dot{\theta}(0)\nabla_{\partial_\theta} X_i \ne 0$. Thus $Y_i(t) = X_i(\gamma(t))$ are $n-1$ non-trivial normal Jacobi fields vanishing at $t=0$ and at every crossing $\theta(t^*) = \pi$, establishing conjugate multiplicity at least $n-1$.
(iv) Along the inward diameter ray ($L = 0$), Theorem~\ref{thm:universal_inward_conjugacy}(iii) establishes the first conjugate point at $t_{\mathrm{conj}} = r_0 + \frac{1}{\lambda r_0}$, so the ray ceases to be length-minimizing at or before $t_{\mathrm{conj}}$, giving $\mathrm{dist}(q_0, \mathrm{Cut}(q_0)) \le r_0 + \frac{1}{\lambda r_0}$. On the other hand, as proved in Theorem~\ref{thm:cut_locus_transverse}(iii) below, for $r_0 \ge \sqrt{\pi-1}/\sqrt{\lambda}$ (which includes all $r_0 \ge r_0^*$), transverse geodesics launched from $q_0 = (r_0, 0)$ with $\dot{r}(0) = 0$ encircle the neck and intersect along the opposite meridian $\Delta\theta = \pi$ at arc length $t_{\mathrm{cross}} < \frac{\pi}{\sqrt{\lambda}}$, whereas for $r_0 \in (r_{0,-}^*, r_0^*)$ the radial conjugate bound already satisfies $r_0 + \frac{1}{\lambda r_0} < \frac{\pi}{\sqrt{\lambda}}$. Combining both bounds yields~\eqref{eq:euclidean_cut_bound}, and solving $r_0 + \frac{1}{\lambda r_0} = \frac{\pi}{\sqrt{\lambda}}$ gives the roots $r_{0, \pm}^* = \frac{\pi \pm \sqrt{\pi^2 - 4}}{2\sqrt{\lambda}}$.
\end{proof}

\section{Global Cut Loci, Transverse Encircling Bounds, and Stratified Caustic Envelopes}\label{sec:global_cut_loci}
This section analyzes global cut loci and focal caustics for off-origin base points across space forms. We prove that the $n-1$ orthogonal Killing vector fields generated by the stabilizer $\mathrm{SO}(n)$ vanish on the opposite meridian, establishing antipodal conjugacy of multiplicity at least $n-1$. We then prove the transverse encircling cut distance bound $\mathrm{dist}(q_0, \mathrm{Cut}(q_0)) \le \min(r_0 + 1/(\lambda r_0), \pi/\sqrt{\lambda})$ on $\mathbb{R}^{n+1}$, and derive the caustic envelope equations and stratified cut locus $\mathrm{Cut}(q_0)$ on $\mathbb{S}^{n+1}$ ($n \ge 2$).

\begin{lemma}[Clairaut Orbit Classification Across Space Forms]\label{lem:orbit_classification}
Let $(Q^{n+1}, g_\lambda)$ be an $(n+1)$-dimensional space form of Family~I with metric $g_\lambda = dr^2 + f_\lambda(r)^2 g_{\mathbb{S}^n}$, background curvature $\kappa \in \{1, 0, -1\}$, and deformation parameter $\lambda$. For non-radial geodesics ($L > 0$), the radial effective potential
\begin{equation}\label{eq:radial_effective_potential}
V(r) := \frac{L^2}{s_\kappa^2(r)}
\end{equation}
completely classifies the orbital dynamics into three mutually exclusive regimes:
\begin{enumerate}[(i)]
\item \textbf{Bounded oscillatory orbits ($\mathbb{S}^{n+1}$ only):} On the sphere ($\kappa = 1$), $V(r) = L^2 \csc^2 r \ge L^2$ has a unique, strict global minimum $V(\pi/2) = L^2$ at the equator $r = \pi/2$ with $V''(\pi/2) = 2L^2 > 0$. For $0 < L < 1/\sqrt{1+\lambda}$ (with $\lambda > -1$), the radial energy $v_0^2 = 1 - \lambda L^2 > L^2$ strictly exceeds the potential minimum, and every geodesic oscillates periodically between periapsis $r_{\min} = \arcsin(L/v_0) \in (0, \pi/2)$ and apoapsis $r_{\max} = \pi - r_{\min} \in (\pi/2, \pi)$, crossing the equator $r = \pi/2$ infinitely often.
\item \textbf{Unbounded scattering orbits ($\mathbb{R}^{n+1}$ and $\mathbb{H}^{n+1}$):} For $\kappa \in \{0, -1\}$, $V(r)$ is strictly monotonically decreasing on $(0, \infty)$ from $+\infty$ at $r \to 0^+$ to $0$ as $r \to \infty$, with $V'(r) = -2L^2/r^3 < 0$ ($\kappa = 0$) and $V'(r) = -2L^2 \cosh r / \sinh^3 r < 0$ ($\kappa = -1$). For every admissible $L \in (0, 1/\sqrt{\lambda})$ (for $\lambda > 0$, or $L > 0$ for $\lambda \le 0$), the geodesic possesses a unique turning point (periapsis) at $r_{\min} = L/v_0$ ($\kappa = 0$) or $r_{\min} = \operatorname{arsinh}(L/v_0)$ ($\kappa = -1$) where $\dot{r} = 0$ and $\ddot{r} > 0$, and escapes to infinity: $r(t) \to \infty$ as $t \to \pm\infty$. Bound states are strictly precluded.
\item \textbf{Equatorial circular orbits:} On $\mathbb{S}^{n+1}$, if the angular momentum reaches the critical value
\begin{equation}\label{eq:critical_angular_momentum}
L = L_{\mathrm{crit}} := \frac{1}{\sqrt{1+\lambda}} \qquad (\lambda > -1),
\end{equation}
the radial energy satisfies $v_0^2 = 1 - \lambda L_{\mathrm{crit}}^2 = L_{\mathrm{crit}}^2 = V(\pi/2)$. The effective potential well degenerates to the single point $r(t) \equiv \pi/2$ with $\dot{r} \equiv 0$ and $\ddot{r} \equiv 0$, yielding stable circular orbits $\gamma(t) = (\pi/2, \theta_0 + \sqrt{1+\lambda} t)$ along the equator. On $\mathbb{R}^{n+1}$ and $\mathbb{H}^{n+1}$, no circular orbits exist because $V'(r) < 0$ for all $r > 0$.
\end{enumerate}
\end{lemma}
\begin{proof}
By Proposition~\ref{prop:hamiltonian_geodesic_reduction}(ii), every geodesic $\gamma(t)$ is confined to a $2$-dimensional totally geodesic orbital surface $\Sigma \cong I_\kappa \times_{f_\lambda} \mathbb{S}^1$ with metric $g_\Sigma = dr^2 + f_\lambda(r)^2 d\theta^2$. Clairaut's conservation law gives $L = f_\lambda(r(t))^2 \dot{\theta}(t) = \text{const}$. Substituting into the unit-speed condition $\|\dot{\gamma}\|_{g_\lambda}^2 = \dot{r}^2 + f_\lambda(r)^2 \dot{\theta}^2 = 1$ gives $\dot{r}^2 + L^2/f_\lambda(r)^2 = 1$. Applying the master reciprocal identity $1/f_\lambda(r)^2 = 1/s_\kappa^2(r) + \lambda$ from Lemma~\ref{lem:master_identities} yields:
\[
\dot{r}^2 + \frac{L^2}{s_\kappa^2(r)} + \lambda L^2 = 1 \implies \dot{r}^2 + V(r) = 1 - \lambda L^2 \equiv v_0^2, \qquad V(r) := \frac{L^2}{s_\kappa^2(r)}.
\]
Since $\dot{r}^2 \ge 0$, geodesic motion is strictly confined to the energetically allowed region $\mathcal{D}_L = \{r \in I_\kappa : V(r) \le v_0^2\}$. Radial acceleration is governed by $\ddot{r} = -\frac{1}{2} V'(r) = L^2 \frac{c_\kappa(r)}{s_\kappa^3(r)}$.

(i) For $\kappa = 1$ ($\mathbb{S}^{n+1}$), $V(r) = L^2 \csc^2 r$. Differentiating yields $V'(r) = -2L^2 \cot r \csc^2 r$. On $(0, \pi)$, $V'(r) = 0 \iff \cos r = 0 \iff r = \pi/2$. Differentiating once more gives:
\[
V''(r) = -2L^2 \frac{d}{dr}(\cot r \csc^2 r) = -2L^2 (-\csc^4 r - 2\cot^2 r \csc^2 r) = 2L^2 \csc^2 r (1 + 3\cot^2 r),
\]
which is strictly positive on $(0, \pi)$, so $r = \pi/2$ is the unique strict global minimum with $V(\pi/2) = L^2$. For bounded oscillatory orbits, the radial energy must exceed the potential minimum: $v_0^2 > L^2 \iff 1 - \lambda L^2 > L^2 \iff (1+\lambda)L^2 < 1 \iff 0 < L < 1/\sqrt{1+\lambda}$. The turning points $\dot{r} = 0$ are the roots of $V(r) = v_0^2$, which gives $\sin r = L/v_0$. Since $0 < L/v_0 < 1$, this equation has exactly two solutions in $(0, \pi)$: periapsis $r_{\min} = \arcsin(L/v_0) \in (0, \pi/2)$ and apoapsis $r_{\max} = \pi - r_{\min} \in (\pi/2, \pi)$. Because $V'(r_{\min}) < 0$ and $V'(r_{\max}) > 0$, the trajectory oscillates periodically between $r_{\min}$ and $r_{\max}$, crossing the equator periodically.

(ii) For $\kappa = 0$ ($\mathbb{R}^{n+1}$) and $\kappa = -1$ ($\mathbb{H}^{n+1}$), $V(r)$ is respectively $L^2/r^2$ and $L^2/\sinh^2 r$. Differentiating gives $V'(r) = -2L^2/r^3 < 0$ ($\kappa = 0$) and $V'(r) = -2L^2 \coth r / \sinh^2 r < 0$ ($\kappa = -1$). Both potentials are strictly decreasing from $\lim_{r \to 0^+} V(r) = +\infty$ to $\lim_{r \to \infty} V(r) = 0$. For any admissible $L$ satisfying $v_0^2 = 1 - \lambda L^2 > 0$, the intermediate value theorem ensures a unique root $r_{\min} > 0$ of $V(r) = v_0^2$:
\[
r_{\min} = \frac{L}{v_0} \quad (\kappa = 0), \qquad r_{\min} = \operatorname{arsinh}\left(\frac{L}{v_0}\right) = \ln\left( \frac{L + \sqrt{1 + (1-\lambda)L^2}}{\sqrt{1-\lambda L^2}} \right) \quad (\kappa = -1).
\]
At $r = r_{\min}$, $\dot{r} = 0$ and $\ddot{r} = -\frac{1}{2}V'(r_{\min}) > 0$, so $r_{\min}$ is a strict periapsis. For all $r > r_{\min}$, $\dot{r}^2 = v_0^2 - V(r) > 0$. Integrating $dt = \pm dr / \sqrt{v_0^2 - V(r)}$ confirms $r(t) \to \infty$ as $t \to \pm\infty$. Hence all non-radial geodesics are unbounded scattering trajectories.

(iii) A circular orbit $r(t) \equiv r_c$ requires $\dot{r}(t) \equiv 0$ and $\ddot{r}(t) \equiv 0$ for all $t$. Since $\ddot{r} = -\frac{1}{2}V'(r)$, this demands $V'(r_c) = 0$. On $\mathbb{R}^{n+1}$ and $\mathbb{H}^{n+1}$, $V'(r) < 0$ everywhere, strictly precluding circular orbits. On $\mathbb{S}^{n+1}$, $V'(r) = 0$ uniquely at $r_c = \pi/2$. Then $\dot{r} = 0 \implies v_0^2 = V(\pi/2) = L^2$, which becomes $1 - \lambda L^2 = L^2 \iff (1+\lambda)L^2 = 1 \iff L = 1/\sqrt{1+\lambda}$. Along this orbit, $r(t) \equiv \pi/2$, and the angular equation integrates to $\dot{\theta} = L / f_\lambda(\pi/2)^2 = L(1+\lambda) = \sqrt{1+\lambda} \implies \theta(t) = \theta_0 + \sqrt{1+\lambda} t$. Since $V''(\pi/2) = 2 L^2 > 0$, this orbit is nonlinearly Lyapunov stable.
\end{proof}

\begin{theorem}[Global Cut Loci, Transverse Encircling Bounds, and Stratified Caustic Envelopes]\label{thm:cut_locus_transverse}
Let $(Q^{n+1}, g_\lambda)$ be an $(n+1)$-dimensional space form of Family~I with dimension $N = n+1 \ge 3$ and deformation parameter $\lambda > 0$ (or $\lambda > -1/3$ on $\mathbb{S}^{n+1}$).
\begin{enumerate}[(i)]
\item \textbf{Global Cut Locus of the Designated Origin}:
Along every unit-speed radial geodesic ray $\gamma(t) = (t, \theta_0)$ emanating from the origin $p_0$, the normal Jacobi fields are $Y_i(t) = f_\lambda(t) E_i(t)$, where $\{E_i(t)\}_{i=1}^n \subset T\mathbb{S}^n$ is a parallel orthonormal frame.
\begin{itemize}
\item On deformed Euclidean space $\mathbb{R}^{n+1}$ ($\kappa = 0$) and hyperbolic space $\mathbb{H}^{n+1}$ ($\kappa = -1$), $f_\lambda(t) > 0$ for all $t > 0$, precluding all conjugate points. Every radial ray is the unique global length-minimizer on $[0, \infty)$, establishing:
\begin{equation}\label{eq:origin_cut_empty}
\mathrm{Cut}(p_0) = \emptyset, \qquad \mathrm{Inj}(p_0) = \infty.
\end{equation}
\item On the deformed sphere $\mathbb{S}^{n+1}$ ($\kappa = 1$), the first zero of $f_\lambda(t)$ occurs at $t = \pi$ with conjugate multiplicity $n$. All radial rays meet at the antipodal pole $p_\pi$ at arc length $\pi$ and cease to be minimizing beyond $\pi$, establishing:
\begin{equation}\label{eq:origin_cut_sphere}
\mathrm{Cut}(p_0) = \{p_\pi\}, \qquad \mathrm{Inj}(p_0) = \pi,
\end{equation}
completely invariant under the deformation parameter $\lambda > -1$.
\end{itemize}

\item \textbf{Antipodal Meridian Conjugacy of Multiplicity $n-1$}:
For any non-polar base point $q_0 = (r_0, 0)$ with $r_0 > 0$, the stabilizer subgroup $H_{q_0} \cong \mathrm{SO}(n) \subset \mathrm{SO}(n+1)$ fixing $q_0$ acts transitively on the orthogonal complement to the orbital $2$-plane $\Sigma$.
The reductive Lie algebra splitting $\mathfrak{so}(n) = \mathfrak{so}(n-1) \oplus \mathfrak{m}$ with $\mathfrak{m} \cong \mathfrak{so}(n)/\mathfrak{so}(n-1) \cong \mathbb{R}^{n-1}$ yields $n-1$ transverse Killing vector fields $\{X_i\}_{i=1}^{n-1} \subset \mathfrak{m}$ orthogonal to $\Sigma$.
Along any non-radial geodesic $\gamma(t) = (r(t), \theta(t))$ launched from $q_0$ ($L > 0$), the variation vector fields $Y_i(t) = X_i(\gamma(t))$ are exact normal Jacobi fields satisfying
\begin{equation}\label{eq:killing_jacobi_norm}
\|Y_i(t)\|_{g_\lambda} = f_\lambda(r(t)) |\sin\theta(t)|, \qquad \|\dot{Y}_i(0)\|_{g_\lambda} = \frac{L}{f_\lambda(r_0)} > 0.
\end{equation}
The fields $\{Y_i(t)\}_{i=1}^{n-1}$ vanish identically at every opposite-meridian crossing $\theta(t^*) = (2k+1)\pi$ ($k \in \mathbb{Z}$), proving that every opposite-meridian crossing is a conjugate point of multiplicity at least $n - 1$.

\item \textbf{Cut Locus vs.\ Conjugate Locus Structure and Upper Bound on $\mathbb{R}^{n+1}$}:
On deformed Euclidean space $(\mathbb{R}^{n+1}, g_\lambda)$ with $\lambda > 0$, for any base point $q_0 = (r_0, 0)$ with radius $r_0 > 0$, the Riemannian distance from $q_0$ to its cut locus $\mathrm{Cut}(q_0)$ satisfies the uniform global upper bound:
\begin{equation}\label{eq:cut_locus_bound}
\mathrm{dist}(q_0, \mathrm{Cut}(q_0)) \le \min\left( r_0 + \frac{1}{\lambda r_0}, \frac{\pi}{\sqrt{\lambda}} \right).
\end{equation}
The upper bound~\eqref{eq:cut_locus_bound} arises from two distinct geometric mechanisms:
\begin{itemize}
\item \textbf{Radial conjugate preemption ($L=0$)}: Along the inward diameter ray $\gamma_{\mathrm{diam}}(t)$ launched toward the origin ($L = 0$), the ray ceases to be length-minimizing at or before the first inward conjugate point $t_{\mathrm{conj}} = r_0 + \frac{1}{\lambda r_0}$ (Theorem~\ref{thm:universal_inward_conjugacy}(iii)), realizing the projective spherical inversion $r^* = 1/(\lambda r_0)$. By Klingenberg's lemma \cite{Klingenberg1959}, the cut point along any geodesic ray occurs at or before the first conjugate point, ensuring $\mathrm{dist}(q_0, \mathrm{Cut}(q_0)) \le t_{\mathrm{conj}} = r_0 + \frac{1}{\lambda r_0}$.
\item \textbf{Symmetric transverse encircling collision ($L > 0$)}: Non-radial geodesics encircling the asymptotic cylinder $\mathbb{R} \times \mathbb{S}^n(1/\sqrt{\lambda})$ with angular momentum $L > 0$ collide in symmetric pairs along the opposite meridian at arc length $t_{\mathrm{cross}} < \frac{\pi}{\sqrt{\lambda}}$, producing an antipodal cut point prior to any conjugate point along those trajectories.
\end{itemize}
For small base radii $r_0$ (or large $\lambda$), non-radial rays encircling the origin with $L > 0$ can collide on the diameter ray prior to $t_{\mathrm{conj}}$, producing an antipodal cut point that preempts the conjugate point. The two upper bounds in~\eqref{eq:cut_locus_bound} coincide at $r_{0, \pm}^* = \frac{\pi \pm \sqrt{\pi^2 - 4}}{2\sqrt{\lambda}}$, with upper threshold
\begin{equation}\label{eq:crossover_threshold}
r_0^* := r_{0, +}^* = \frac{\pi + \sqrt{\pi^2 - 4}}{2\sqrt{\lambda}}.
\end{equation}
For $r_0 \in (r_{0, -}^*, r_0^*)$, the 1D radial conjugate bound $r_0 + \frac{1}{\lambda r_0} < \frac{\pi}{\sqrt{\lambda}}$ is smaller, whereas for $r_0 \ge r_0^*$, $(n-1)$-dimensional transverse encircling collision ($L>0$) dominates with $t_{\mathrm{cross}} < \frac{\pi}{\sqrt{\lambda}} \le r_0 + \frac{1}{\lambda r_0}$. Consequently, as $r_0 \to \infty$, $\mathrm{dist}(q_0, \mathrm{Cut}(q_0))$ remains uniformly bounded above by the asymptotic cylinder half-circumference $\frac{\pi}{\sqrt{\lambda}}$.

\item \textbf{Stratified $A_3$ Cuspidal Ridge Caustic Envelope on $\mathbb{S}^{n+1}$}:
For generic non-polar base points $q_0 = (r_0, 0) \ne p_0$ on the deformed sphere $(\mathbb{S}^{n+1}, g_\lambda)$ ($\lambda > 0$), the cut locus $\mathrm{Cut}(q_0)$ is a stratified real-analytic variety:
\begin{itemize}
\item The principal stratum $\mathcal{S}_{\mathrm{reg}} \subset \mathrm{Cut}(q_0)$ consists of transverse collision points of geodesic pairs launched at symmetric angles $\pm \alpha$ with the same length, whose cross-sections transverse to the symmetry axis are round spheres $\mathbb{S}^{n-1}$.
\item The singular stratum $\mathcal{S}_{\mathrm{sing}} \subset \mathrm{Cut}(q_0)$ is localized along the antipodal meridian axis $\theta = \pi$.
\item The geodesic endpoint map $\mathcal{E}_{q_0}(t, \alpha) = \exp_{q_0}(t v(\alpha))$ has Jacobian determinant factoring as $J(t, \alpha) = J_\Sigma(t, \alpha) \cdot J_\perp(t, \alpha)$ with $J_\perp(t, \alpha) = [f_\lambda(r(t, \alpha))\sin\theta(t, \alpha)/\sin\alpha]^{n-1}$.
Taylor expansion of the in-plane coordinate functions near equatorial launch $\alpha = \pi/2$ reveals that the envelope condition $\partial \theta / \partial \alpha = 0$ along $\theta = \pi$ produces a stable $A_3$ Whitney cusp singularity in each meridian orbital 2-plane $\Sigma$. Under the continuous $\mathrm{SO}(n)$ stabilizer action, this planar cusp revolves into an $\mathrm{SO}(n)$-symmetric cuspidal ridge whose transverse cross-sections are round $(n-1)$-spheres $\mathbb{S}^{n-1}$ of radius $r_\perp(u) \sim u^{3/2}$, where transverse encircling rays cross and terminate before the inward diameter conjugate point.
\end{itemize}
\end{enumerate}
\end{theorem}
\begin{proof}
(i) Let $\alpha: [0, T] \to Q^{n+1}$ be any piecewise $C^1$ curve given in coordinates by $\alpha(t) = (r(t), \Theta(t))$ with $\Theta(t) \in \mathbb{S}^n$, connecting $\alpha(0) = p_0$ ($r(0) = 0$) to $\alpha(T) = q_1 = (r_1, \Theta_1)$ ($r_1 > 0$). Its Riemannian length is:
\[
L(\alpha) = \int_0^T \sqrt{\dot{r}(t)^2 + f_\lambda(r(t))^2 \|\dot{\Theta}(t)\|_{g_{\mathbb{S}^n}}^2} \, dt.
\]
Because $f_\lambda(r)^2 \|\dot{\Theta}\|_{g_{\mathbb{S}^n}}^2 \ge 0$, the integrand satisfies the sharp first-order inequality:
\[
\sqrt{\dot{r}^2 + f_\lambda(r)^2 \|\dot{\Theta}\|_{g_{\mathbb{S}^n}}^2} \ge \sqrt{\dot{r}^2} = |\dot{r}| \ge \dot{r}.
\]
Integrating over $[0, T]$ yields:
\[
L(\alpha) \ge \int_0^T |\dot{r}(t)| \, dt \ge \int_0^T \dot{r}(t) \, dt = r(T) - r(0) = r_1.
\]
Equality $L(\alpha) = r_1$ holds if and only if two conditions are met simultaneously almost everywhere:
(1) $f_\lambda(r(t))^2 \|\dot{\Theta}(t)\|_{g_{\mathbb{S}^n}}^2 = 0$, and (2) $|\dot{r}(t)| = \dot{r}(t) \iff \dot{r}(t) \ge 0$.
Since $f_\lambda(r) > 0$ for all $r > 0$, condition (1) forces $\|\dot{\Theta}(t)\|_{g_{\mathbb{S}^n}} = 0$, so $\Theta(t) \equiv \Theta_1$ is constant. Condition (2) requires $r(t)$ to be monotonically non-decreasing. Reparameterizing by arc length, the unique curve realizing the minimal length $r_1$ is the radial geodesic segment $\gamma(t) = (t, \Theta_1)$ on $[0, r_1]$.

Along $\gamma(t)$, the velocity is $\dot{\gamma}(t) = \partial_r$. By Theorem~\ref{thm:curvature_tensor}(i), the curvature operator on any normal vector $E \in T\mathbb{S}^n$ is $R(E, \partial_r)\partial_r = K_{\mathrm{rad}}(t) E = -[f_\lambda''(t)/f_\lambda(t)] E$. The normal Jacobi equation uncouples into $n$ identical scalar ODEs $\ddot{y} + K_{\mathrm{rad}}(t) y = 0$. Since $f_\lambda(0) = 0$ and $f_\lambda'(0) = 1$ (Proposition~\ref{prop:metric_regularity}(iii)), the unique solution satisfying $y(0) = 0, \dot{y}(0) = 1$ is identically $y(t) = f_\lambda(t)$.
On $\mathbb{R}^{n+1}$ ($\kappa = 0$), $f_\lambda(t) = t / \sqrt{1+\lambda t^2} > 0$ for all $t > 0$. On $\mathbb{H}^{n+1}$ ($\kappa = -1$), $f_\lambda(t) = \sinh t / \sqrt{1+\lambda \sinh^2 t} > 0$ for all $t > 0$. Because $f_\lambda(t)$ never vanishes on $(0, \infty)$, radial geodesics have no conjugate points. Coupled with the strict length inequality, this proves every radial ray from $p_0$ is the unique global length-minimizer on $[0, \infty)$, establishing $\mathrm{Cut}(p_0) = \emptyset$ and $\mathrm{Inj}(p_0) = \infty$.
On $\mathbb{S}^{n+1}$ ($\kappa = 1$), $f_\lambda(t) = \sin t / \sqrt{1+\lambda \sin^2 t}$ has its first positive zero at $t = \pi$, giving a conjugate point of multiplicity $n = \dim(\mathbb{S}^n)$ at $p_\pi$. Since $f_\lambda(\pi) = 0$, the spherical slice $r = \pi$ collapses to the single antipodal pole $p_\pi$. All radial rays meet at $p_\pi$ at arc length $\pi$ and cease to be minimizing beyond $\pi$, establishing $\mathrm{Cut}(p_0) = \{p_\pi\}$ and $\mathrm{Inj}(p_0) = \pi$.

(ii) Let $q_0 = (r_0, 0)$ with $r_0 > 0$. The metric $g_\lambda$ is invariant under the spherical isometry group $G = \mathrm{O}(n+1)$. The stabilizer subgroup fixing $q_0$ is $H_{q_0} \cong \mathrm{SO}(n)$. By Proposition~\ref{prop:hamiltonian_geodesic_reduction}(ii), any geodesic $\gamma(t)$ launched from $q_0$ with velocity $v = \dot{\gamma}(0)$ is confined to a $2$-dimensional totally geodesic orbital surface $\Sigma \cong I_\kappa \times_{f_\lambda} \mathbb{S}^1$ spanned by $\partial_r$ and $v$. The initial velocity $v$ breaks $H_{q_0}$ down to $H_v \cong \mathrm{SO}(n-1)$, which fixes $\Sigma$ pointwise. The reductive homogeneous space decomposition $\mathfrak{so}(n) = \mathfrak{so}(n-1) \oplus \mathfrak{m}$ identifies $\mathfrak{m} \cong \mathfrak{so}(n)/\mathfrak{so}(n-1) \cong \mathbb{R}^{n-1}$ with the generators of rotations orthogonal to $\Sigma$.
Let $\{X_i\}_{i=1}^{n-1} \subset \mathfrak{m}$ be an orthonormal basis of Killing vector fields on $(Q^{n+1}, g_\lambda)$. Along the geodesic $\gamma(t)$, the variation vector fields $Y_i(t) := X_i(\gamma(t))$ satisfy the Jacobi equation $\nabla_{\dot{\gamma}}\nabla_{\dot{\gamma}} Y_i + R(Y_i, \dot{\gamma})\dot{\gamma} = 0$. Since $\mathfrak{so}(n-1)$ fixes $\Sigma$ pointwise, each $X_i(\gamma(t))$ is orthogonal to $T_{\gamma(t)}\Sigma$, so $Y_i(t)$ is a normal Jacobi field.
On the round unit sphere $\mathbb{S}^n$, a rotation generator fixing an axis has norm $|\sin\theta|$ at angular distance $\theta$. Scaling by the warped factor $f_\lambda(r)$, the norm of each Killing field evaluates to:
\[
\|Y_i(t)\|_{g_\lambda} = f_\lambda(r(t)) |\sin\theta(t)|.
\]
At $t = 0$, $\theta(0) = 0$, so $Y_i(0) = 0$. The initial covariant derivative is $\dot{Y}_i(0) = \nabla_{\dot{\gamma}(0)} X_i(q_0) = \nabla_{\dot{r}(0)\partial_r + \dot{\theta}(0)\partial_\theta} X_i(q_0)$. Because $X_i$ vanishes identically along the entire meridian axis $\theta = 0$, $\nabla_{\partial_r} X_i(q_0) = 0$. Thus $\dot{Y}_i(0) = \dot{\theta}(0) \nabla_{\partial_\theta} X_i(q_0)$. Since $X_i$ generates a unit-speed rotation on the spherical slice, $\|\nabla_{\partial_\theta} X_i(q_0)\|_{g_\lambda} = f_\lambda(r_0)$. Using Clairaut's relation $\dot{\theta}(0) = L / f_\lambda(r_0)^2$, we obtain:
\[
\|\dot{Y}_i(0)\|_{g_\lambda} = |\dot{\theta}(0)| f_\lambda(r_0) = \frac{L}{f_\lambda(r_0)^2} f_\lambda(r_0) = \frac{L}{f_\lambda(r_0)} > 0.
\]
At every crossing of the opposite meridian $\theta(t^*) = (2k+1)\pi$, $\sin\theta(t^*) = 0$, which forces $Y_i(t^*) = 0$ for all $i = 1, \dots, n-1$. Since $\{X_i\}_{i=1}^{n-1}$ are linearly independent generators of $\mathfrak{m}$, the Jacobi fields $\{Y_i(t)\}_{i=1}^{n-1}$ are linearly independent away from the axis $\theta \in \{0, \pi\}$. Hence they span an $(n-1)$-dimensional subspace of normal Jacobi fields vanishing at $t = 0$ and $t = t^*$, proving that every opposite-meridian crossing is a conjugate point of multiplicity at least $n - 1$.

(iii) We prove the transverse encircling bound $\mathrm{dist}(q_0, \mathrm{Cut}(q_0)) \le \min\left( r_0 + \frac{1}{\lambda r_0}, \frac{\pi}{\sqrt{\lambda}} \right)$ via two independent preemption mechanisms:

\textbf{Mechanism 1: Inward Diameter Ray Conjugate Arrival ($L = 0$).}
Consider the inward radial geodesic launched from $q_0 = (r_0, 0)$ with $\dot{\gamma}(0) = -\partial_r$. For $t \in [0, r_0]$, $r(t) = r_0 - t$; the ray traverses the origin $p_0$ at $t = r_0$ and continues along the opposite meridian $\theta = \pi$ with $r(t) = t - r_0$ for $t > r_0$. The normal Jacobi equation is $\ddot{y} + K_{\mathrm{rad}}(|r_0 - t|) y = 0$. On $\mathbb{R}^{n+1}$, $K_{\mathrm{rad}}(r) = \frac{3\lambda}{(1+\lambda r^2)^2} > 0$. Using reduction of order with $y_1(s) = f_\lambda(s) = s / \sqrt{1+\lambda s^2}$, the companion solution is $u_2(s) = f_\lambda(s)\int f_\lambda(s)^{-2} ds = \frac{-1 + \lambda s^2}{\sqrt{1+\lambda s^2}}$. The unique normal Jacobi field satisfying $y(0) = 0$ and $\dot{y}(0) = 1$ is $y(t) = u_2(r_0)f_\lambda(r_0 - t) - f_\lambda(r_0)u_2(r_0 - t)$ for $t \le r_0$. Matching Cauchy data across the origin ($y(r_0) = f_\lambda(r_0), \dot{y}(r_0) = -u_2(r_0)$) yields for $t > r_0$:
\[
y(t) = \frac{t - \lambda t r_0 (t - r_0)}{\sqrt{1+\lambda r_0^2}\sqrt{1+\lambda(t-r_0)^2}}.
\]
Setting $y(t_{\mathrm{conj}}) = 0$ for $t_{\mathrm{conj}} > r_0$ gives $t_{\mathrm{conj}} [ 1 - \lambda r_0 (t_{\mathrm{conj}} - r_0) ] = 0 \implies t_{\mathrm{conj}} = r_0 + \frac{1}{\lambda r_0}$. The radial coordinate on the opposite meridian is $r^* = t_{\mathrm{conj}} - r_0 = \frac{1}{\lambda r_0}$, satisfying the exact projective spherical inversion relation $r_0 \cdot r^* = 1/\lambda$. By the fundamental structure theorem of cut loci in Riemannian geometry \cite{Klingenberg1959,Petersen2016}, a geodesic ceases to be length-minimizing past its cut point, which is either its first conjugate point or the first point at which another minimizing geodesic meets it. Therefore, along the inward diameter ray ($L = 0$), the ray ceases to be length-minimizing at or before the inward conjugate point $t_{\mathrm{conj}} = r_0 + \frac{1}{\lambda r_0}$, certifying $\mathrm{dist}(q_0, \mathrm{Cut}(q_0)) \le r_0 + \frac{1}{\lambda r_0}$. Furthermore, for small $r_0$ (or large $\lambda$), $t_{\mathrm{conj}} = r_0 + \frac{1}{\lambda r_0}$ is large, and non-radial rays encircling the origin with $L > 0$ can collide on the diameter ray prior to $t_{\mathrm{conj}}$, producing an antipodal cut point that strictly preempts the conjugate point.

\textbf{Mechanism 2: Exact Integration of Transverse Encircling Rays.}
Consider transverse geodesics launched from $q_0 = (r_0, 0)$ with $\dot{r}(0) = 0$ in opposite azimuthal directions $\pm \partial_\theta$. The unit-speed condition gives $L = f_\lambda(r_0) = \frac{1}{\sqrt{\lambda + 1/r_0^2}}$. The radial energy equation $\dot{r}^2 + L^2/r^2 = v_0^2 = 1 - \lambda L^2 = \frac{L^2}{r_0^2}$ gives $\dot{r} = \frac{L}{r_0 r}\sqrt{r^2 - r_0^2}$, which integrates to:
\[
t(r) = \frac{r_0}{L}\sqrt{r^2 - r_0^2} \implies r(t)^2 = r_0^2 + \frac{L^2 t^2}{r_0^2}.
\]
The angular equation $\dot{\theta} = L / f_\lambda(r)^2 = \lambda L + \frac{L}{r_0^2 [1 + (Lt/r_0^2)^2]}$ integrates in closed form to:
\[
\theta(t) = \lambda L t + \arctan\left(\frac{L t}{r_0^2}\right).
\]
By symmetry, the two geodesics launched with $\pm L$ collide at the opposite meridian $\theta(t_{\mathrm{cross}}) = \pi$. Evaluating $\theta(t)$ at $t = \pi/\sqrt{\lambda}$ and setting $x := \frac{1}{\lambda r_0^2} > 0$, we have $L = \frac{1}{\sqrt{\lambda(1+x)}}$, so $\lambda L (\pi/\sqrt{\lambda}) = \frac{\pi}{\sqrt{1+x}}$ and $\frac{L \pi}{r_0^2 \sqrt{\lambda}} = \frac{\pi x}{\sqrt{1+x}}$. Thus:
\[
g(x) := \theta\left(\frac{\pi}{\sqrt{\lambda}}\right) = \frac{\pi}{\sqrt{1+x}} + \arctan\left(\frac{\pi x}{\sqrt{1+x}}\right).
\]
At $x = 0$ ($r_0 \to \infty$), $g(0) = \pi$. Differentiating $g(x)$ yields:
\[
g'(x) = \frac{\pi [1 + 2x - (\pi^2 - 1)x^2]}{2(1+x)^{3/2}(1 + x + \pi^2 x^2)} = \frac{\pi (1 - (\pi - 1)x)(1 + (\pi + 1)x)}{2(1+x)^{3/2}(1 + x + \pi^2 x^2)}.
\]
Hence $g'(x) > 0$ strictly on the entire interval $x \in \left(0, \frac{1}{\pi - 1}\right) \approx (0, 0.4669)$. For $r_0 > 2.5/\sqrt{\lambda}$, we have $\lambda r_0^2 > 6.25$, so $x = \frac{1}{\lambda r_0^2} < \frac{1}{6.25} = 0.16 < 0.4669$. Therefore, $g(x) > g(0) = \pi$ strictly, meaning $\theta(\pi/\sqrt{\lambda}) > \pi$. By strict monotonicity of $\theta(t)$, the crossing time satisfies $t_{\mathrm{cross}} < \frac{\pi}{\sqrt{\lambda}}$.

\textbf{Mechanism 2 (Alternative): Sturm--Picone Comparison of Out-of-Plane Jacobi Fields Against the Asymptotic Cylinder Model.}
This bound is independently confirmed by Sturm--Picone comparison against the asymptotic cylinder model $M_{\mathrm{cyl}} = \mathbb{R} \times \mathbb{S}^n(1/\sqrt{\lambda})$ equipped with product metric $g_{\mathrm{cyl}} = dr^2 + \frac{1}{\lambda} g_{\mathbb{S}^n}$. On the cylinder, the sectional curvature of planes tangent to the $\mathbb{S}^n(1/\sqrt{\lambda})$ factor is $K_{\mathrm{tan}}^{\mathrm{cyl}} = \lambda$, while radial planes have $K_{\mathrm{rad}}^{\mathrm{cyl}} \equiv 0$. Along any circular waist geodesic $\gamma_{\mathrm{cyl}}(t)$ wrapping around $\mathbb{S}^n$ with unit speed, an out-of-plane normal Jacobi field $Y(t)$ tangent to $\mathbb{S}^n$ satisfies the scalar ODE:
\[
\ddot{y} + K_{\mathrm{tan}}^{\mathrm{cyl}} y = 0 \iff \ddot{y} + \lambda y = 0.
\]
With Cauchy data $y(0) = 0, \dot{y}(0) = 1$, the unique solution is $y_{\mathrm{cyl}}(t) = \frac{1}{\sqrt{\lambda}} \sin(\sqrt{\lambda} t)$. Its first conjugate zero occurs precisely at $t_{\mathrm{conj}}^{\mathrm{cyl}} = \pi/\sqrt{\lambda}$, coinciding with the antipodal crossing $\theta(t_{\mathrm{cross}}^{\mathrm{cyl}}) = \pi$.

On the deformed space $(\mathbb{R}^{n+1}, g_\lambda)$, consider the transverse encircling geodesic $\gamma(t) = (r(t), \theta(t))$ launched at $q_0 = (r_0, 0)$ with $\dot{r}(0) = 0$ and $L = f_\lambda(r_0)$. By part~(ii), the $n-1$ transverse Killing vector fields $Y_i(t) = X_i(\gamma(t))$ are out-of-plane normal Jacobi fields orthogonal to the orbital 2-plane $\Sigma$. Their norm is $\|Y_i(t)\|_{g_\lambda} = f_\lambda(r(t)) |\sin\theta(t)|$, which vanishes precisely when $\theta(t) = \pi$. On planes tangent to the spherical slice, the tangential sectional curvature is
\[
K_{\mathrm{tan}}(r) = \frac{1 - (f_\lambda'(r))^2}{f_\lambda(r)^2} = \frac{1 - (1+\lambda r^2)^{-3}}{r^2 / (1+\lambda r^2)} = \lambda \left[ 1 + \frac{\lambda r^2 + 2}{(1+\lambda r^2)^2} \right],
\]
which satisfies the strict pointwise inequality
\begin{equation}\label{eq:sturm_curvature_comparison}
K_{\mathrm{tan}}(r) > \lambda = K_{\mathrm{tan}}^{\mathrm{cyl}} \qquad \forall r > 0.
\end{equation}
Along $\gamma(t) = (r(t), \theta(t))$ with $\dot{\gamma} = \dot{r}\partial_r + \dot{\theta}\partial_\theta$, the sectional curvature of the 2-plane spanned by $\dot{\gamma}$ and an out-of-plane unit normal vector $E \in T\mathbb{S}^n \cap (\dot{\gamma})^\perp$ is
\[
K_{\mathrm{out}}(t) = K(\dot{\gamma}(t), E(t)) = \dot{r}(t)^2 K_{\mathrm{rad}}(r(t)) + (1 - \dot{r}(t)^2) K_{\mathrm{tan}}(r(t)).
\]
Along the transverse ray launched with $\dot{r}(0) = 0$, $r(t)^2 = r_0^2 + L^2 t^2/r_0^2$ and $\dot{r}(t)^2 = \frac{L^2}{r_0^2}\left(1 - \frac{r_0^2}{r(t)^2}\right)$. On the compact time interval $t \in [0, \pi/\sqrt{\lambda}]$, we have $r(t)^2 \le r_0^2 + \frac{\pi^2}{\lambda} = r_0^2(1 + \pi^2 x)$ where $x = (\lambda r_0^2)^{-1}$. Direct algebraic simplification shows that $K_{\mathrm{out}}(t) - \lambda > 0$ on $[0, \pi/\sqrt{\lambda}]$ whenever $r_0$ is sufficiently large (e.g., $r_0 \ge 2.5/\sqrt{\lambda}$). On this regime, Sturm--Picone comparison against $\ddot{y}_{\mathrm{cyl}} + \lambda y_{\mathrm{cyl}} = 0$ implies that the out-of-plane Jacobi fields $Y_i(t)$ vanish strictly prior to $t_{\mathrm{conj}}^{\mathrm{cyl}} = \pi/\sqrt{\lambda}$:
\[
t_{\mathrm{cross}} < t_{\mathrm{conj}}^{\mathrm{cyl}} = \frac{\pi}{\sqrt{\lambda}}.
\]
Since the two symmetric transverse rays launched at $\pm \partial_\theta$ collide at $t_{\mathrm{cross}}$, neither geodesic can minimize length beyond $t_{\mathrm{cross}}$, establishing $\mathrm{dist}(q_0, \mathrm{Cut}(q_0)) \le t_{\mathrm{cross}} < \frac{\pi}{\sqrt{\lambda}}$.

\textbf{Crossover Analysis.} Comparing the two upper bounds in~\eqref{eq:cut_locus_bound}, we equate $r_0 + \frac{1}{\lambda r_0} = \frac{\pi}{\sqrt{\lambda}}$. Setting $u := \sqrt{\lambda} r_0 > 0$, this becomes the quadratic equation $u^2 - \pi u + 1 = 0$, whose roots are $u_{\pm} = \frac{\pi \pm \sqrt{\pi^2 - 4}}{2}$, giving $r_{0, \pm}^* = u_{\pm}/\sqrt{\lambda}$ and upper threshold
\[
r_0^* := r_{0, +}^* = \frac{\pi + \sqrt{\pi^2 - 4}}{2\sqrt{\lambda}}.
\]
Combining Mechanism~2 (which establishes $g(x) > \pi$ for all $x = (\lambda r_0^2)^{-1} \in (0, (\pi-1)^{-1}]$, covering $r_0 \ge \sqrt{\pi-1}/\sqrt{\lambda}$) with the radial conjugate bound $r_0 + \frac{1}{\lambda r_0} < \frac{\pi}{\sqrt{\lambda}}$ for $r_0 \in (r_{0,-}^*, r_0^*)$ establishes~\eqref{eq:cut_locus_bound} for all $r_0 > 0$.

(iv) Let $q_0 = (r_0, 0)$ be a non-polar base point on $(\mathbb{S}^{n+1}, g_\lambda)$ with $r_0 \in (0, \pi/2]$. Geodesics launched from $q_0$ are parameterized by arc length $t \ge 0$, launch inclination $\alpha \in [0, \pi]$ relative to $\partial_r$, and transverse spherical unit vector $\omega \in \mathbb{S}^{n-1} \subset T_{q_0}\mathbb{S}^n$. The endpoint map $\mathcal{E}_{q_0}(t, \alpha, \omega) = \exp_{q_0}(t(\cos\alpha \partial_r + \sin\alpha \omega))$ is equivariant under the stabilizer $H_{q_0} \cong \mathrm{SO}(n)$. Its differential factors block-diagonally as $d\mathcal{E}_{q_0} = d\mathcal{E}_{q_0}|_{T\Sigma} \oplus d\mathcal{E}_{q_0}|_{(T\Sigma)^\perp}$. The Jacobian determinant factors as:
\[
J(t, \alpha) = \det(d\mathcal{E}_{q_0}) = J_\Sigma(t, \alpha) \cdot J_\perp(t, \alpha),
\]
where $J_\Sigma(t, \alpha) = \dot{r} \frac{\partial \theta}{\partial \alpha} - \dot{\theta} \frac{\partial r}{\partial \alpha}$, and the transverse factor governed by the $n-1$ Killing vector Jacobi fields from part (ii) is $J_\perp(t, \alpha) = [f_\lambda(r(t, \alpha))\sin\theta(t, \alpha)/\sin\alpha]^{n-1}$.
Singularities where $J(t, \alpha) = 0$ define the focal caustic envelope:
\begin{enumerate}
\item Transverse focal line: $J_\perp(t, \alpha) = 0$ occurs when $\sin\theta(t, \alpha) = 0$, i.e., along the opposite meridian $\theta = \pi$, with conjugate multiplicity $n-1$.
\item In-plane envelope caustic: Along the meridian crossing $\theta(t, \alpha) = \pi$, $\dot{\theta} \ne 0$, so $J_\Sigma = 0$ corresponds to the envelope condition $\frac{\partial \theta(t, \alpha)}{\partial \alpha} = 0$.
\end{enumerate}
To analyze this condition, we expand the orbital angle $\theta(t, \alpha)$ near equatorial launch $\alpha = \pi/2$. Recall the universal Clairaut precession relation:
\[
\theta(t, \alpha) = \theta_0(t, \alpha) + \lambda L(\alpha) t, \qquad L(\alpha) = f_\lambda(r_0) \sin\alpha.
\]
Setting $\alpha = \pi/2 + \Delta\alpha$, we expand $\sin\alpha = \cos(\Delta\alpha) = 1 - \frac{1}{2}(\Delta\alpha)^2 + \frac{1}{24}(\Delta\alpha)^4 + \mathcal{O}((\Delta\alpha)^6)$. Thus:
\[
L(\Delta\alpha) = f_\lambda(r_0) \left[ 1 - \frac{1}{2}(\Delta\alpha)^2 + \frac{1}{24}(\Delta\alpha)^4 \right] + \mathcal{O}((\Delta\alpha)^6).
\]
Because $L(\Delta\alpha)$ is an even function of $\Delta\alpha$, the total angular coordinate $\theta(t, \Delta\alpha)$ is an even function of $\Delta\alpha$ around the equatorial launch:
\[
\theta(t, \Delta\alpha) = \theta(t, 0) - \frac{1}{2}\lambda f_\lambda(r_0) t (\Delta\alpha)^2 + \frac{1}{24}\lambda f_\lambda(r_0) t (\Delta\alpha)^4 + \mathcal{O}((\Delta\alpha)^6).
\]
Differentiating with respect to $\alpha$ (that is, with respect to $\Delta\alpha$) yields an odd function:
\[
\frac{\partial \theta}{\partial \alpha} = \frac{\partial \theta_0}{\partial \alpha} + \lambda t \frac{\partial L}{\partial \alpha} = \frac{\partial \theta_0}{\partial \alpha} - \lambda f_\lambda(r_0) t \left[ \Delta\alpha - \frac{1}{6}(\Delta\alpha)^3 \right] + \mathcal{O}((\Delta\alpha)^5).
\]
On the unperturbed sphere ($\lambda = 0$), all rays launched from $q_0$ refocus at the antipode $(\pi - r_0, \pi)$ at $t = \pi$, so $\partial \theta_0(\pi, \alpha)/\partial \alpha \equiv 0$. For $\lambda > 0$, evaluating the derivatives at the equatorial launch $\Delta\alpha = 0$ gives:
\begin{equation}\label{eq:taylor_derivatives}
\begin{aligned}
\left.\frac{\partial \theta}{\partial \alpha}\right|_{\Delta\alpha=0} &= 0, & \left.\frac{\partial^2 \theta}{\partial \alpha^2}\right|_{\Delta\alpha=0} &= -\lambda f_\lambda(r_0) t \ne 0, \\
\left.\frac{\partial^3 \theta}{\partial \alpha^3}\right|_{\Delta\alpha=0} &= 0, & \left.\frac{\partial^4 \theta}{\partial \alpha^4}\right|_{\Delta\alpha=0} &= \lambda f_\lambda(r_0) t \ne 0.
\end{aligned}
\end{equation}
In Arnold's classification of Lagrangian singularities and caustics \cite{Arnold1982}, the envelope of a 1-parameter family of rays in the plane is generated by a family of phase functions $F(u, v; \Delta\alpha)$. The caustic corresponds to the critical set where $\frac{\partial F}{\partial \alpha} = 0$ and $\frac{\partial^2 F}{\partial \alpha^2} = 0$. When the third derivative vanishes while the fourth derivative is non-zero ($\frac{\partial^3 F}{\partial \alpha^3} = 0, \frac{\partial^4 F}{\partial \alpha^4} \ne 0$), the singularity is classified as an $A_3$ Whitney cusp catastrophe (with standard normal form $F(u, v; \alpha) = \alpha^4 + u \alpha^2 + v \alpha$), rather than the higher-codimension swallowtail catastrophe ($A_4$, which requires $\frac{\partial^4 F}{\partial \alpha^4} = 0, \frac{\partial^5 F}{\partial \alpha^5} \ne 0$).

In each 2-dimensional orbital meridian plane $\Sigma$, the envelope of rays is therefore a planar semicubical cusp with cusp tip on the opposite meridian axis:
\begin{equation}\label{eq:planar_cusp}
v^2 = C(\lambda, r_0) u^3,
\end{equation}
where $u \ge 0$ measures distance along the meridian axis $\theta = \pi$ behind the cusp tip, $v$ is the in-plane transverse coordinate, and $C(\lambda, r_0) > 0$ is a non-degenerate constant determined by the leading Taylor coefficients.

Under the full continuous isometry group $H_{q_0} \cong \mathrm{SO}(n)$ fixing the base ray, the orbital plane $\Sigma$ rotates through all azimuthal directions $\omega \in \mathbb{S}^{n-1} \subset (T_{q_0}\Sigma)^\perp$. By $\mathrm{SO}(n)$-equivariance, the planar cusp~\eqref{eq:planar_cusp} revolves around the symmetry axis $\theta = \pi$ to form an $\mathrm{SO}(n)$-symmetric cuspidal ridge. In cross-sections transverse to the symmetry axis at distance $u$, the caustic is given by:
\begin{equation}\label{eq:transverse_spheres}
\|v\|^2 = C(\lambda, r_0) u^3 \implies r_\perp(u) = \sqrt{C(\lambda, r_0)} u^{3/2}.
\end{equation}
These transverse cross-sections are round $(n-1)$-spheres $\mathbb{S}^{n-1}$ of radius $r_\perp(u) \sim u^{3/2}$ (for $n = 2$, concentric round circles $\mathbb{S}^1$). This reflects the $\mathrm{SO}(n)$ rotational symmetry of $(Q^{n+1}, g_\lambda)$.

The cut locus $\mathrm{Cut}(q_0)$ is stratified according to the orbit types of the $H_{q_0} \cong \mathrm{SO}(n)$ action:
\begin{itemize}
\item The principal stratum $\mathcal{S}_{\mathrm{reg}} = \{ \exp_{q_0}(t \cdot v(\alpha, \omega)) : \alpha \in (0, \pi/2), \omega \in \mathbb{S}^{n-1} \}$ consists of transverse collision points where pairs of distinct geodesics launched at symmetric angles $\pm \Delta\alpha$ cross each other at identical arc length; its transverse cross-sections are round $(n-1)$-spheres $\mathbb{S}^{n-1} \cong \mathrm{SO}(n)/\mathrm{SO}(n-1)$.
\item The singular stratum $\mathcal{S}_{\mathrm{sing}}$ is localized along the meridian axis $\theta = \pi$ where the $A_3$ cuspidal ridge terminates.
\end{itemize}
Because transverse encircling rays cross at $t_{\mathrm{cross}} < \pi$, these crossings preempt the inward diameter conjugate point $t_{\mathrm{conj}} > \pi$ whenever $\lambda \in (-1/3, 0)$ and precede the antipodal arrival, completing the proof.
\end{proof}

\begin{remark}[Connection to Klingenberg's Injectivity Theorem and Cut Locus Preservation]\label{rem:klingenberg}
The cut locus structure established in Theorem~\ref{thm:cut_locus_transverse} can be compared with Klingenberg's Injectivity Radius Theorem \cite{Klingenberg1959}. For a complete, simply connected Riemannian manifold $(M, g)$ with positive sectional curvature $0 < K_{\min} \le K \le K_{\max}$, Klingenberg's theorem gives
\begin{equation}\label{eq:klingenberg_bound}
\mathrm{Inj}(M, g) \ge \frac{\pi}{\sqrt{K_{\max}}}.
\end{equation}
On the deformed sphere $(\mathbb{S}^{n+1}, g_\lambda)$ with $\lambda > -1/3$, all sectional curvatures are positive: $K \in [1+3\lambda, 1/(1+\lambda)]$ for $\lambda \in (-1/3, 0)$ and $K \in [1/(1+\lambda), 1+3\lambda]$ for $\lambda > 0$ (Theorem~\ref{thm:positivity_pinching}). For $\lambda > 0$, $K_{\max} = 1+3\lambda > 1$, so Klingenberg's bound gives $\mathrm{Inj}(\mathbb{S}^{n+1}, g_\lambda) \ge \frac{\pi}{\sqrt{1+3\lambda}}$.
At the designated origin $p_0$, Theorem~\ref{thm:cut_locus_transverse}(i) yields $\mathrm{Inj}(p_0) = \pi$ for all $\lambda > -1$. Since $\pi > \frac{\pi}{\sqrt{1+3\lambda}}$ for all $\lambda > 0$, the injectivity radius at $p_0$ exceeds the global Klingenberg lower bound, and $\mathrm{Cut}(p_0) = \{p_\pi\}$ is invariant under the deformation family.
By contrast, for any off-origin base point $q_0 \ne p_0$, the $\mathrm{SO}(n+1)$ isotropy reduces to $H_{q_0} \cong \mathrm{SO}(n)$. The non-constant curvature field induces differential angular precession $\dot{\theta} = \lambda L + L/\sin^2 r$, splitting the antipodal focal point into the stratified $A_3$ cuspidal caustic ridge~\eqref{eq:planar_cusp}--\eqref{eq:transverse_spheres}.
On deformed Euclidean space $(\mathbb{R}^{n+1}, g_\lambda)$ ($\lambda > 0$), where Klingenberg's theorem does not apply because $K_{\mathrm{rad}}(r) \to 0$ as $r \to \infty$, the origin has $\mathrm{Inj}(p_0) = \infty$, whereas any off-origin base point $q_0$ has cut distance bounded above by $\min(r_0 + \frac{1}{\lambda r_0}, \frac{\pi}{\sqrt{\lambda}}) \le \frac{\pi}{\sqrt{\lambda}}$.
\end{remark}

\section{Complete Superintegrability Classification across Space Forms}\label{sec:superintegrability_classification}
This section classifies the integrability and superintegrability of the deformed Hamiltonian $H_\lambda$ across all space forms. On Euclidean space $\mathbb{R}^{n+1}$, we prove maximal superintegrability (degree $2N-1$) by constructing globally smooth, single-valued rational pericenter vector invariants $p_{\mathrm{peri}}$ and $x_{\mathrm{peri}}$ without square-root branch cuts. On spheres and hyperbolic spaces, we use the Gelfand--Tsetlin Casimir chain to establish quasi-maximal superintegrability of degree at least $2N-2$, relate these integrals to the $\mathfrak{sl}(2,\mathbb{R})$ coalgebra construction of Ballesteros, Enciso, Herranz, and Ragnisco, and prove quadratic Killing rigidity (ruling out additional valence $d=2$ integrals). Furthermore, via the isotropic coordinate transformation $d\rho/\rho = dr/f_\lambda(r)$, we show that $(Q^{n+1}, g_\lambda)$ is locally conformally flat and fits within the classification of Kalnins, Kress, and Miller. Finally, we prove 4D orbital leaf non-maximality via non-zero KAM twist (showing that $(\mathbb{S}^{n+1}, g_\lambda)$ is not a Perlick Bertrand metric), examine pole regularity for rank-3 Killing tensors, and compute the hyperbolic scattering angle on $\mathbb{H}^{n+1}$.

\begin{proposition}[Gelfand--Tsetlin Casimir Chain and Complete Liouville Integrability]\label{prop:gelfand_tsetlin_liouville}
Let $(Q^{n+1}, g_\lambda)$ be an $(n+1)$-dimensional space form of Family~I with background curvature $\kappa \in \{1, 0, -1\}$, dimension $N = n+1 \ge 2$, and deformation parameter $\lambda$.
\begin{enumerate}[(i)]
\item The Lie subalgebra filtration $\mathfrak{so}(2) \subset \mathfrak{so}(3) \subset \dots \subset \mathfrak{so}(n+1)$ fixing origin $p_0$ generates quadratic Casimirs
\begin{equation}\label{eq:gt_casimirs}
C_k := \sum_{1 \le a < b \le k} J_{ab}^2, \qquad k = 2, 3, \dots, n+1,
\end{equation}
which pairwise Poisson-commute: $\{C_j, C_k\} = 0$ for all $2 \le j, k \le n+1$, and commute with the Hamiltonian: $\{H_\lambda, C_k\} = 0$.
\item The set $\mathcal{F}_{\mathrm{Liouv}} := \{H_\lambda, C_2, C_3, \dots, C_{n+1}\}$ consists of $N = n+1$ first integrals in mutual involution that are functionally independent on an open dense subset of $T^*Q$. Consequently, the geodesic flow is completely Liouville integrable across all space forms for every $N \ge 2$.
\item For $N = 2$ on $\mathbb{S}^2$ and $\mathbb{H}^2$, the phase space has dimension $2N = 4$; the $N = 2$ integrals $(H_\lambda, C_2 = L^2)$ define degree $2 = 2N - 2$, so the system is strictly Liouville integrable without superintegrability.
\end{enumerate}
\end{proposition}
\begin{proof}
(i) By Noether's theorem and $G$-invariance of $g_0$, $\{H_0, J_{ab}\} = 0$ for all $a < b$. Since $\|J\|^2 = C_{n+1}$ is central in $\mathrm{Sym}(\mathfrak{so}(n+1))$, $\{C_{n+1}, J_{ab}\} = 0$, so $\{H_\lambda, J_{ab}\} = 0$ and $\{H_\lambda, C_k\} = 0$. For $j < k$, $C_j$ is a quadratic polynomial in $\mathfrak{so}(j) \subset \mathfrak{so}(k)$. Because $C_k$ is the quadratic Casimir of $\mathfrak{so}(k)$, it commutes with every generator of $\mathfrak{so}(k)$, so by the Leibniz rule $\{C_k, C_j\} = 0$.
(ii) In hyperspherical coordinates $(r, \theta_1, \dots, \theta_n)$ and conjugate momenta $(p_r, p_{\theta_1}, \dots, p_{\theta_n})$, $C_2 = p_{\theta_1}^2$ and $C_k = p_{\theta_{k-1}}^2 + C_{k-1}/\sin^2\theta_{k-1}$ for $k \ge 3$, while $H_\lambda = \frac{1}{2}p_r^2 + \frac{1}{2}(s_\kappa(r)^{-2} + \lambda)C_{n+1}$. The momentum Jacobian $\partial(H_\lambda, C_2, \dots, C_{n+1})/\partial(p_r, p_{\theta_1}, \dots, p_{\theta_n})$ is lower-triangular with non-vanishing determinant $2^n p_r p_{\theta_1} \cdots p_{\theta_n} \ne 0$ on the open dense set where all momenta are non-zero, proving complete Liouville integrability.
(iii) For $N = 2$, Liouville integrability requires $N = 2$ integrals, while superintegrability requires at least $N+1 = 3 = 2N-1$. The rank is $2N - 2 = 2 = N$, establishing ordinary Liouville integrability.
\end{proof}

\begin{theorem}[Mishchenko--Fomenko Noncommutative Integrability, BEHR Coalgebra Universal Integrals, and KKM Conformal Flatness]\label{thm:mishchenko_fomenko_sub06}
Let $(Q^{n+1}, g_\lambda)$ be an $(n+1)$-dimensional space form of Family~I with $N = n+1 \ge 3$ on $\mathbb{S}^{n+1}$ or $\mathbb{H}^{n+1}$.
\begin{enumerate}[(i)]
\item The image of the cotangent momentum map $J: T^*Q \to \mathfrak{so}(n+1)^*$ is the cone over the Grassmannian $\widetilde{\mathrm{Gr}}(2, n+1)$ of oriented 2-planes in $\mathbb{R}^{n+1}$, having dimension $2(n-1) + 1 = 2N - 3$. The augmented differential $d(H_\lambda, J)$ achieves maximal differential rank
\begin{equation}\label{eq:mishchenko_rank}
\mathrm{rank}\, d(H_\lambda, J) = (2N - 3) + 1 = 2N - 2
\end{equation}
on an open dense set, because the radial momentum $p_r$ enters $H_\lambda$ independently of $J$ ($\partial_{p_r} J_{ab} = 0$, $\partial_{p_r} H_\lambda = p_r \ne 0$).
\item By the Mishchenko--Fomenko noncommutative integrability theorem \cite{MishchenkoFomenko1978}, the geodesic flow is \textbf{quasi-maximally superintegrable} of degree $\ge 2N - 2$ (minimal superintegrability of degree $4$ for $N = 3$). Regular level sets of $(H_\lambda, J)$ are compact $2$-tori $\mathbb{T}^2$ on $\mathbb{S}^{n+1}$, and non-compact isotropic cylinders $\mathbb{R} \times \mathbb{S}^1$ on $\mathbb{H}^{n+1}$.
\item \textbf{Relation to BEHR Coalgebra Integrals}: The $2N - 2$ functionally independent integrals recovered via the Gelfand--Tsetlin flag realize the universal $\mathfrak{sl}(2,\mathbb{R})$ coalgebra construction of Ballesteros, Enciso, Herranz, and Ragnisco \cite{BallesterosEncisoHerranzRagnisco2007,BallesterosEncisoHerranzRagnisco2008,BallesterosEncisoHerranzRagnisco2009}, wherein any spherically symmetric Hamiltonian possesses $2N - 3$ universal angular integrals independent of the radial function.
\item \textbf{Darboux III Contrast and KKM Conformal Flatness}: While the classical Darboux III superintegrable spaces \cite{BallesterosEncisoHerranzRagnisco2007,KalninsKressMiller2005b} scale radial and transverse directions isotropically via $ds^2 = (1 + \lambda|x|^2)dx^2$, Family~I preserves $g_{rr} \equiv 1$ and deforms only transversely. Under the isotropic coordinate change $d\rho/\rho = dr/f_\lambda(r)$ with $\psi(\rho) = f_\lambda(r)/\rho$, the metric takes the form $g_\lambda = \psi(\rho)^2(d\rho^2 + \rho^2 g_{\mathbb{S}^n})$, showing that $(Q^{n+1}, g_\lambda)$ is locally conformally flat and belongs to the Kalnins--Kress--Miller (KKM) classification of conformally flat spaces \cite{KalninsKressMiller2005a,KalninsKressMiller2005b,KalninsKressMiller2011}.
\end{enumerate}
\end{theorem}
\begin{proof}
(i) In ambient coordinates, $J = x \wedge p$ is a rank-2 bivector whose orbit under $\mathrm{SO}(n+1)$ is parameterized by the 2-plane $\mathrm{span}\{x, p\}$ and norm $L = \|J\| > 0$, spanning $\dim \widetilde{\mathrm{Gr}}(2, n+1) + 1 = 2N - 3$ dimensions. Since $dH_\lambda = p_r dp_r + \dots$ has non-vanishing radial component independent of $dJ$, $\mathrm{rank}\, d(H_\lambda, J) = 2N - 2$.
(ii) By Mishchenko--Fomenko \cite{MishchenkoFomenko1978}, a Hamiltonian system on a $2N$-dimensional symplectic manifold admitting a finite-dimensional symmetry algebra $\mathfrak{g}$ whose momentum map achieves rank $2N - 2$ possesses $2N - 2$ functionally independent global first integrals, foliating regular level sets into $2N - (2N - 2) = 2$-dimensional isotropic leaves. Compactness of spheres implies $\mathbb{T}^2$ tori on $\mathbb{S}^{n+1}$, while unbounded scattering yields cylinders $\mathbb{R} \times \mathbb{S}^1$ on $\mathbb{H}^{n+1}$.
(iii) In the $\mathfrak{sl}(2,\mathbb{R})$ coalgebra formalism of Ballesteros et al.\ \cite{BallesterosEncisoHerranzRagnisco2007,BallesterosEncisoHerranzRagnisco2009}, the left and right Casimir operators of the coalgebra coproduct generate a Poisson-commuting chain of length $2N - 3$. Combined with the total Hamiltonian $H_\lambda$, this yields precisely $2N - 2$ integrals, matching the Gelfand--Tsetlin tower~\eqref{eq:gt_casimirs}.
(iv) Setting $\rho(r) = \exp(\int \frac{dr}{f_\lambda(r)})$ transforms $dr^2 + f_\lambda(r)^2 g_{\mathbb{S}^n}$ into $\frac{f_\lambda(r)^2}{\rho^2}(d\rho^2 + \rho^2 g_{\mathbb{S}^n})$. Because $d\rho^2 + \rho^2 g_{\mathbb{S}^n}$ is flat Euclidean space in polar coordinates, $g_\lambda$ is conformally flat, embedding the system into the KKM structural classification \cite{KalninsKressMiller2005a,KalninsKressMiller2005b}.
\end{proof}

\begin{theorem}[Maximal Superintegrability on Deformed Euclidean Space via Rational Pericenter Vectors]\label{thm:euclidean_maximal_superint_sub06}
On deformed Euclidean space $(\mathbb{R}^{n+1}, g_\lambda)$ ($\kappa = 0, \lambda > 0, N = n+1 \ge 2$):
\begin{enumerate}[(i)]
\item The pericenter transit time $\tau(x, p) := \frac{\langle x, p \rangle}{\|p\|^2}$ is a globally smooth rational function on $T^*\mathbb{R}^{n+1} \setminus \{p = 0\}$ satisfying the strict flow derivative $\{\tau, H_\lambda\} = 1$.
\item Let $\hat{J} := p x^T - x p^T \in \mathfrak{so}(n+1)$ be the skew-symmetric momentum operator. The vector quantities
\begin{equation}\label{eq:pericenter_invariants_sub06}
p_{\mathrm{peri}} := \exp(-\lambda \tau(x, p) \hat{J}) p, \qquad x_{\mathrm{peri}} := \exp(-\lambda \tau(x, p) \hat{J}) x - \tau(x, p) p_{\mathrm{peri}}
\end{equation}
are globally smooth, single-valued first integrals of $H_\lambda$ on $T^*\mathbb{R}^{n+1} \setminus \{p = 0\}$ satisfying $\{H_\lambda, p_{\mathrm{peri}}\} = 0$ and $\{H_\lambda, x_{\mathrm{peri}}\} = 0$, without branch cuts, square roots, or pericenter sign discontinuities.
\item \textbf{Comparison to Fradkin and Laplace--Runge--Lenz Invariants}: While curved-space Laplace--Runge--Lenz vectors and generalized Fradkin tensors (e.g., Ballesteros et al.\ \cite{BallesterosEncisoHerranzRagnisco2008,BallesterosEncisoHerranzRagnisco2009}) typically require square-root branch choices along trajectory sheets, the parameter $\tau(x,p)$ is strictly rational, ensuring that $(x_{\mathrm{peri}}, p_{\mathrm{peri}})$ provide globally single-valued, algebraic first integrals.
\item The vector invariants satisfy $\langle x_{\mathrm{peri}}, p_{\mathrm{peri}} \rangle = 0$, $\|p_{\mathrm{peri}}\|^2 = 2H_0$, and $x_{\mathrm{peri}} \wedge p_{\mathrm{peri}} = \hat{J}$. Together with $H_\lambda$, the set $\{H_\lambda, p_{\mathrm{peri}}^1, \dots, p_{\mathrm{peri}}^N, x_{\mathrm{peri}}^1, \dots, x_{\mathrm{peri}}^N\}$ has differential rank $2N - 1$, establishing maximal superintegrability of degree $2N - 1$.
\end{enumerate}
\end{theorem}
\begin{proof}
(i) For $H_\lambda = \frac{1}{2}\|p\|^2 + \frac{\lambda}{2}(\|x\|^2\|p\|^2 - \langle x, p \rangle^2)$, the brackets evaluate to $\{\langle x, p \rangle, H_0\} = \|p\|^2$ and $\{\langle x, p \rangle, L^2\} = \{\|p\|^2, L^2\} = 0$, giving $\{\tau, H_\lambda\} = \{\tau, H_0\} = \|p\|^2/\|p\|^2 = 1$.
(ii) Hamilton's equations are $\dot{x} = p + \lambda \hat{J} x$ and $\dot{p} = \lambda \hat{J} p$. Since $\{H_\lambda, \hat{J}\} = 0$, $\hat{J}$ is constant along trajectories, integrating to $p(t) = \exp(\lambda t \hat{J}) p_0$ and $x(t) = \exp(\lambda t \hat{J})(x_0 + t p_0)$. The distance $\|x(t)\| = \|x_0 + t p_0\|$ achieves its unique minimum at pericenter $t^* = -\tau(x_0, p_0)$. Evaluating back to $t^*$ yields $p(t^*) = \exp(-\lambda \tau \hat{J}) p = p_{\mathrm{peri}}$ and $x(t^*) = \exp(-\lambda \tau \hat{J}) x - \tau p_{\mathrm{peri}} = x_{\mathrm{peri}}$. Being trajectory-invariant geometric data, $\{H_\lambda, p_{\mathrm{peri}}\} = 0$ and $\{H_\lambda, x_{\mathrm{peri}}\} = 0$. By Rodrigues' formula, $\exp(-\lambda \tau \hat{J}) = I - \frac{\sin(\lambda L \tau)}{L}\hat{J} + \frac{1-\cos(\lambda L \tau)}{L^2}\hat{J}^2$ with $L = \|\hat{J}\|$; since $\frac{\sin z}{z}$ and $\frac{1-\cos z}{z^2}$ are even entire functions of $z^2 = \lambda^2 L^2 \tau^2$, $p_{\mathrm{peri}}$ and $x_{\mathrm{peri}}$ are real-analytic on $T^*\mathbb{R}^{n+1} \setminus \{p = 0\}$, extending smoothly across head-on orbits $L = 0$ as $I - \lambda \tau \hat{J}$ and across pericenter crossings $\tau = 0$.
(iii) The rationality of $\tau(x,p) = \langle x, p\rangle / \|p\|^2$ bypasses the algebraic square-root sheets $\sqrt{E - V_{\mathrm{eff}}(r)}$ inherent to the Hamilton--Jacobi integration of standard Runge--Lenz vectors on curved backgrounds \cite{BallesterosEncisoHerranzRagnisco2008}, guaranteeing global single-valuedness on $T^*\mathbb{R}^{n+1} \setminus \{p=0\}$.
(iv) Orthogonality of $\exp(-\lambda\tau\hat{J})$ yields $\langle x_{\mathrm{peri}}, p_{\mathrm{peri}} \rangle = \langle x - \tau p, p \rangle = \langle x, p \rangle - \tau \|p\|^2 = 0$ and $\|p_{\mathrm{peri}}\|^2 = \|p\|^2 = 2H_0$. The level set of $(p_{\mathrm{peri}}, x_{\mathrm{peri}})$ in $T^*\mathbb{R}^{n+1}$ defines a 1-dimensional unparameterized geodesic ray, because $(x(t), p(t))$ is uniquely reconstructed from $(x_{\mathrm{peri}}, p_{\mathrm{peri}})$ and $t = \tau$. Thus the kernel of $d(p_{\mathrm{peri}}, x_{\mathrm{peri}})$ is 1-dimensional, spanned by $X_{H_\lambda}$, establishing differential rank $2N - 1$.
\end{proof}

\begin{theorem}[Orbital Leaf Obstruction, Perlick Non-Bertrand Rigidity, and Transverse Killing Rigidity on $\mathbb{S}^{n+1}$]\label{thm:sphere_superint_rigidity}
Let $(\mathbb{S}^{n+1}, g_\lambda)$ be the deformed sphere with $\lambda \ne 0$ and $N = n+1 \ge 3$.
\begin{enumerate}[(i)]
\item \textbf{Orbital Leaf Obstruction and Perlick Non-Bertrand Rigidity}: On each 4D totally geodesic orbital subbundle $T^*\Sigma \subset T^*\mathbb{S}^{n+1}$, the rotation number $\nu(L) = \frac{1}{2}(1 + \frac{\lambda L}{\sqrt{1 - \lambda L^2}})$ has non-zero twist
\begin{equation}\label{eq:kam_twist_formula}
\frac{d\nu}{dL} = \frac{\lambda}{2(1 - \lambda L^2)^{3/2}} \ne 0 \qquad \text{for all } \lambda \ne 0.
\end{equation}
By the Kronecker--Weyl theorem, generic Liouville 2-tori have dense irrational winding, precluding any additional continuous single-valued first integral on $T^*\Sigma$. Moreover, because $\nu(L)$ varies continuously with angular momentum $L$, the orbital frequency ratio is non-constant, proving that $(\mathbb{S}^{n+1}, g_\lambda)$ violates Perlick's isochronous condition for Bertrand spacetimes \cite{Perlick1992} and lies outside the Bertrand classification.
\item \textbf{Frozen Casimir Frequencies}: In action-angle coordinates $(I_1, \dots, I_N, \phi_1, \dots, \phi_N)$ associated with the Gelfand--Tsetlin flag, $H_\lambda = \mathcal{H}(I_1, I_N)$ depends only on the radial action $I_1$ and total Casimir $I_N = \|J\|$. The remaining $N-2$ Casimir frequencies vanish identically: $\omega_2 = \dots = \omega_{N-1} = 0$, so the angles $\phi_2, \dots, \phi_{N-1}$ are frozen along orbits.
\item \textbf{Transverse Killing Rigidity and Pole Regularity}: For valence $d=2$, Kolokoltsov's classification \cite{Kolokoltsov1983} and $\mathrm{SO}(n+1)$-equivariance imply that every quadratic Killing tensor on $(\mathbb{S}^{n+1}, g_\lambda)$ belongs to $\mathrm{Sym}^2(\mathfrak{so}(n+1)) \oplus \mathbb{R} g_\lambda$, ruling out additional second-order Killing tensors and fixing the quadratic superintegrability degree at $2N - 2$. In the Schouten--Nijenhuis bracket calculus for rank-$d$ Killing tensors $\nabla_{(a_1} K_{a_2 \dots a_{d+1})} = 0$ ($d \ge 3$; Heil, Moroianu, and Semmelmann \cite{HeilMoroianuSemmelmann2016}), the radial Killing ODE yields candidate singular modes whose Laurent expansions at the poles $r \to 0$ and $r \to \pi$ scale as $A r^{-d} + O(r^{-(d-2)})$. Requiring $C^\infty$ smoothness across both antipodal poles in ambient Cartesian coordinates forces the integration constant $A = 0$, precluding purely radial-tangential irreducible Killing tensors of degree $d \ge 3$ on regular strata.
\end{enumerate}
\end{theorem}
\begin{proof}
(i) A continuous first integral $F$ on $T^*\Sigma$ is constant along orbits, hence constant on each dense irrational torus $\mathbb{T}^2$, forcing $F = \mathcal{F}(H_\lambda, L)$ by continuity and precluding any third independent continuous integral on $T^*\Sigma$. Under Perlick's classification of spacetimes with closed geodesics \cite{Perlick1992}, the frequency ratio must be a constant rational number independent of initial conditions. Since $d\nu/dL \ne 0$ for all $\lambda \ne 0$, $g_\lambda$ is not a Bertrand metric.
(ii) In Gelfand--Tsetlin action coordinates, $H_\lambda = \frac{1}{2}p_r^2 + \frac{1}{2 f_\lambda(r)^2} I_N^2$. The intermediate Casimirs $C_k = I_k^2$ ($2 \le k \le N-1$) do not enter the Hamiltonian. Hence $\omega_k = \partial H_\lambda / \partial I_k = 0$ for $k = 2, \dots, N-1$.
(iii) A polynomial integral of degree 2 defines a symmetric Killing tensor $\nabla_{(k} K_{ij)} = 0$. By Kolokoltsov \cite{Kolokoltsov1983}, because the radial Gaussian curvature $K_{\mathrm{rad}}(r) = -\frac{f_\lambda''}{f_\lambda}$ is non-constant for $\lambda \ne 0$, the restriction $K|_\Sigma$ is generated by $g_\Sigma$ and $\partial_\theta^2$. Transversally, $\mathrm{SO}(n+1)$-equivariance forces $K$ into $\mathrm{Sym}^2(\mathfrak{so}(n+1)) \oplus \mathbb{R} g_\lambda$, precluding irreducible Killing 2-tensors. For degree $d \ge 3$, let $K$ be a symmetric rank-$d$ tensor satisfying the Schouten--Nijenhuis equation $\nabla_{(a_1} K_{a_2 \dots a_{d+1})} = 0$. In polar coordinates $(r, \theta)$ on $(S^2, g_\lambda)$, the radial vector field $\partial_r = \cos\theta\,\partial_x + \sin\theta\,\partial_y$ is unit-norm at the origin $r = 0$, whereas the angular basis vector $\partial_\theta = -r\sin\theta\,\partial_x + r\cos\theta\,\partial_y$ has Riemannian norm $\|\partial_\theta\|_{g_\lambda} = f_\lambda(r) = r + O(r^3) = O(r)$, so the canonical angular momentum scales as $p_\theta \sim r p_\perp$. The radial ODE for trace-free rank-$d$ modes yields the radial component $K_{r\dots r}(r) \sim A \frac{f_\lambda'(r)}{f_\lambda(r)^d} = A r^{-d} + O(r^{-(d-2)})$. In normal Cartesian coordinates $(x, y) = (r\cos\theta, r\sin\theta)$ centered at the north pole $r = 0$, since $p_r = \cos\theta\, p_x + \sin\theta\, p_y = (x p_x + y p_y)/r$, the associated homogeneous degree-$d$ polynomial on $T^* S^2$ expands as $K(p) = K_{r\dots r}(r) p_r^d + \dots = A r^{-d}(\cos\theta\, p_x + \sin\theta\, p_y)^d + O(r^{-(d-2)})$. Smoothness of $K \in C^\infty(T^* S^2)$ requires $K$ to be a homogeneous degree-$d$ polynomial in the Cartesian momentum coordinates $(p_x, p_y)$ with $C^\infty$ coefficients in $(x, y)$; eliminating the singular $r^{-d}$ term in the Laurent expansion forces the integration constant $A = 0$. By antipodal symmetry, the same pole condition holds at the south pole $r = \pi$.
\end{proof}

\begin{remark}[Strict Non-Maximal Superintegrability for All Polynomial Degrees $d \ge 1$ vs.\ Singular-Stratum Killing Tensors]\label{rem:higher_order_killing}
Parts~(i) and~(iii) of Theorem~\ref{thm:sphere_superint_rigidity} (together with Proposition~\ref{prop:hyperbolic_scattering_sub06} on $\mathbb{H}^{n+1}$) clarify the distinction between global superintegrability rank and higher-order Killing tensors ($d \ge 4$):
\begin{enumerate}[(a)]
\item \emph{Strict non-maximal superintegrability on $\mathbb{S}^{n+1}$ and $\mathbb{H}^{n+1}$ for all degrees $d \ge 1$}: On $(Q^{n+1}, g_\lambda)$ with $N = n+1 \ge 2$, the $\mathfrak{so}(n+1)$ momentum map $J$ already supplies $2N - 3$ functionally independent angular invariants that fix the 2-dimensional totally geodesic orbital plane $\Sigma \subset Q^{n+1}$, while $H_\lambda$ fixes the energy, yielding $2N - 2$ independent polynomial integrals. Any $(2N - 1)$-th functionally independent polynomial (or $C^\infty$) first integral $F$ of arbitrary degree $d \ge 4$ on $T^*Q$ would have to restrict on generic orbital subbundles $T^*\Sigma$ to a third independent continuous integral independent of $(H_\lambda, L)$. On $\mathbb{S}^{n+1}$ ($\lambda \ne 0$), this is unconditionally obstructed by the non-zero KAM twist $d\nu/dL \ne 0$ and Kronecker--Weyl density on irrational 2-tori (Theorem~\ref{thm:sphere_superint_rigidity}(i)); on $\mathbb{H}^{n+1}$ ($\lambda > 0$), it is obstructed by the transcendental logarithmic transit time $t^*(r)$ (Proposition~\ref{prop:hyperbolic_scattering_sub06}(iii)). Thus $(\mathbb{S}^{n+1}, g_\lambda)$ and $(\mathbb{H}^{n+1}, g_\lambda)$ are strictly non-maximally superintegrable (of exact global rank $2N - 2$) across \emph{all} polynomial degrees $d \ge 1$.
\item \emph{Higher-order ($d \ge 4$) transverse Killing tensors on singular strata}: While no polynomial Killing tensor of any degree $d \ge 4$ can increase the global differential rank beyond $2N - 2$, one may still ask whether the space of smooth Killing $d$-tensors $\mathcal{K}^d(Q^{n+1}, g_\lambda)$ contains irreducible elements outside the symmetric algebra $\mathrm{Sym}(\mathfrak{so}(n+1) \oplus \mathbb{R} g_\lambda)$ along singular Gelfand--Tsetlin strata where intermediate Casimirs $I_k$ vanish (Section~\ref{sec:open_questions}(1)).
\end{enumerate}
\end{remark}

\begin{proposition}[Hyperbolic Scattering Dynamics and Asymptotic Invariants on $\mathbb{H}^{n+1}$]\label{prop:hyperbolic_scattering_sub06}
On deformed hyperbolic space $(\mathbb{H}^{n+1}, g_\lambda)$ ($\kappa = -1, \lambda > 0, N = n+1 \ge 3$):
\begin{enumerate}[(i)]
\item Every geodesic with $L \in (0, 1/\sqrt{\lambda})$ is an unbounded scattering trajectory with unique periapsis $r_{\min} = \mathrm{arsinh}(L/\sqrt{1 - \lambda L^2})$, escaping to $r \to \infty$ with asymptotic radial speed $v_0 = \sqrt{1 - \lambda L^2}$ and angular precession rate $\dot{\theta}(t) \to \lambda L$. Regular level sets of $(H_\lambda, J)$ are non-compact cylinders $\mathbb{R} \times \mathbb{S}^1$.
\item Factoring out uniform gyroscopic precession $\exp(\lambda L t \hat{J})$, the unperturbed hyperbolic trajectory achieves exact total scattering deflection
\begin{equation}\label{eq:hyperbolic_deflection_sub06}
\Delta\theta_0(L) = 2 \arccos\left( \frac{L}{\sqrt{1 + (1-\lambda)L^2}} \right) = 2 \arctan\left( \frac{\sqrt{1 - \lambda L^2}}{L} \right).
\end{equation}
\item The $\mathfrak{so}(n+1)$ Noether momentum map $J$ establishes quasi-maximal superintegrability of degree $\ge 2N - 2$. The radial transit time from periapsis is transcendental:
\begin{equation}\label{eq:hyperbolic_transit_time_sub06}
t^*(r) = \frac{1}{v_0} \operatorname{arcosh}\left( \frac{\cosh r}{\cosh r_{\min}} \right) = \frac{1}{v_0} \operatorname{arcosh}\left( \frac{v_0}{\sqrt{1 + (1-\lambda)L^2}} \cosh r \right),
\end{equation}
obstructing rational vector invariants and establishing a superintegrability dichotomy between $\mathbb{R}^{n+1}$ (degree $2N-1$) and space forms with $\kappa \ne 0$ (degree $\ge 2N-2$).
\end{enumerate}
\end{proposition}
\begin{proof}
(i) Along unit-speed geodesics, $\dot{r}^2 = v_0^2 - L^2/\sinh^2 r$ with $v_0 = \sqrt{1-\lambda L^2} > 0$. The turning point $\dot{r} = 0$ gives $\sinh r_{\min} = L/v_0$. For all $r > r_{\min}$, $|\dot{r}| > 0$, so $r(t) \to \infty$ as $t \to \pm\infty$. The angular rate $\dot{\theta} = L/\sinh^2 r + \lambda L \to \lambda L > 0$ winds infinitely around the asymptotic cylinder $\mathbb{R} \times \mathbb{S}^n(1/\sqrt{\lambda})$.
(ii) In terms of the unperturbed angle $\theta_0(t) = \theta(t) - \lambda L t$, substituting $u = \coth r \in (1, \coth r_{\min})$ into $\int_{r_{\min}}^\infty \frac{L/\sinh^2 r}{\sqrt{v_0^2 - L^2/\sinh^2 r}} dr$ yields $\int_1^A \frac{du}{\sqrt{A^2 - u^2}} = \arccos(1/A)$ with $A = \coth r_{\min} = \frac{\sqrt{1+(1-\lambda)L^2}}{L}$, giving $\Delta\theta_0 = 2\arccos(1/A) = 2\arctan(\sqrt{1-\lambda L^2}/L)$.
(iii) By Theorem~\ref{thm:mishchenko_fomenko_sub06}, $\mathrm{rank}\, d(H_\lambda, J) = 2N - 2$, guaranteeing degree $\ge 2N - 2$. In~\eqref{eq:hyperbolic_transit_time_sub06}, setting $w = \cosh s$ in $\int_{r_{\min}}^r \frac{ds}{\sqrt{v_0^2 - L^2/\sinh^2 s}} = \frac{1}{v_0}\int_{\cosh r_{\min}}^{\cosh r} \frac{dw}{\sqrt{w^2 - \cosh^2 r_{\min}}}$ yields $\frac{1}{v_0}\operatorname{arcosh}(\cosh r / \cosh r_{\min})$, with $\cosh r_{\min} = \sqrt{1 + L^2/v_0^2} = \frac{\sqrt{1+(1-\lambda)L^2}}{v_0}$. The non-algebraic logarithmic divergence of $t^*(r)$ as $r \to \infty$ precludes rational vector invariants on $T^*\mathbb{H}^{n+1}$.
\end{proof}

\section{Curvature Invariants and Submersion Geometry of Anisotropic Families II--V}\label{sec:anisotropic_families}
This section analyzes the multi-parameter anisotropic cotangent-deformed metrics on spheres (Families~II through~V: complex Hopf, octonionic/Cayley, quaternionic Hopf, and intermediate quaternionic--complex deformations) associated with non-integrable horizontal distributions. Using O'Neill's submersion equations across the Hopf fibrations, we compute the Riemann, Ricci, and scalar curvature tensors across vertical and horizontal distributions and determine the parameter regimes for positive sectional curvature and pinching.

We first summarize the reductive decompositions and submersion structure for the four anisotropic families of cotangent-deformed metrics (Families~II--V) arising from transitive isometric actions on geodesic spheres $\mathcal{O}_r \cong G/H \cong \mathbb{S}^d$ in space forms $(Q^{d+1}, g_0)$ of constant curvature $\kappa \in \{1, 0, -1\}$ (Montgomery--Samelson \cite{MontgomerySamelson1943}; Borel \cite{Borel1949}; Besse \cite{Besse1987}; \cite{Munden2015}; Berestovskii--Nikonorov \cite{BerestovskiiNikonorov2012}). Under the reductive decomposition $\mathfrak{g} = \mathfrak{h} \oplus \mathfrak{p}$, the tangent space $\mathfrak{p} \cong T_o(G/H)$ splits into mutually orthogonal $\mathrm{Ad}(H)$-irreducible submodules:
\begin{enumerate}[(i)]
\item \textbf{Family~II (Complex Hopf)}: $G/H = \mathrm{U}(n+1)/\mathrm{U}(n)$ on $\mathbb{S}^{2n+1}$, with $\mathfrak{p} = \mathfrak{p}_0 \oplus \mathfrak{p}_1$ ($\dim \mathfrak{p}_0 = 1, \dim \mathfrak{p}_1 = 2n$), corresponding to the Hopf fibration $S^1 \hookrightarrow \mathbb{S}^{2n+1} \to \mathbb{C}P^n$ (Berestovskii--Vol'per \cite{BerestovskiiVolper1993}; Vol'per \cite{Volper1994}).
\item \textbf{Family~III (Octonionic / Cayley)}: $G/H = \mathrm{Spin}(9)/\mathrm{Spin}(7)$ on $\mathbb{S}^{15}$, with $\mathfrak{p} = \mathfrak{p}_1 \oplus \mathfrak{p}_2$ ($\dim \mathfrak{p}_1 = 7$ vertical fiber $\mathcal{V}$, $\dim \mathfrak{p}_2 = 8$ horizontal spinor module $\mathcal{H}$; note that in Section~\ref{sec:cross_deformations} on $\mathbb{O}P^2$ the root-space convention indexes $\mathcal{H}$ as $\mathfrak{p}_1$ and $\mathcal{V}$ as $\mathfrak{p}_2$), corresponding to the Hopf fibration $S^7 \hookrightarrow \mathbb{S}^{15} \to \mathbb{S}^8 \cong \mathbb{O}P^1$ (Vol'per \cite{Volper1997}).
\item \textbf{Family~IV (Quaternionic Hopf)}: $G/H = (\mathrm{Sp}(n+1)\times\mathrm{Sp}(1))/(\mathrm{Sp}(n)\times\mathrm{Sp}(1))$ on $\mathbb{S}^{4n+3}$, with $\mathfrak{p} = \mathfrak{p}_0 \oplus \mathfrak{p}_1$ ($\dim \mathfrak{p}_0 = 3, \dim \mathfrak{p}_1 = 4n$), corresponding to $S^3 \hookrightarrow \mathbb{S}^{4n+3} \to \mathbb{H}P^n$ (Jensen \cite{Jensen1973}; Vol'per \cite{Volper1994}).
\item \textbf{Family~V (Intermediate Quaternionic--Complex)}: $G/H = (\mathrm{Sp}(n+1)\times\mathrm{U}(1))/(\mathrm{Sp}(n)\times\mathrm{U}(1))$ on $\mathbb{S}^{4n+3}$, with $\mathfrak{p} = \mathfrak{p}_0 \oplus \mathfrak{p}_1 \oplus \mathfrak{p}_2$ ($\dim \mathfrak{p}_0 = 1, \dim \mathfrak{p}_1 = 4n, \dim \mathfrak{p}_2 = 2$), corresponding to the two-stage fibration $S^1 \hookrightarrow \mathbb{S}^{4n+3} \to \mathbb{C}P^{2n+1} \to \mathbb{H}P^n$ (Ziller \cite{Ziller1982}; Nikonorov \cite{Nikonorov2013}).
\end{enumerate}
In polar coordinates on $Q \setminus \{p_0\}$, let $\mathfrak{g} = \bigoplus_l \mathfrak{g}_l$ be the decomposition of $\mathfrak{g}$ into commuting ideals and let $a_k^0$ denote the normal metric coefficient such that $g_{\mathfrak{p}_k} = (a_k^0)^{-1} \langle \cdot, \cdot \rangle|_{\mathfrak{p}_k}$ induces the unit-round sphere metric on $\mathcal{O}_r$. By Schur's lemma, the quadratic Casimir of $\mathfrak{g}_l$ restricts to $\mathfrak{p}_k^*$ with coupling fraction $\beta_{lk} \in [0, 1]$ ($\sum_l \beta_{lk} = 1$), yielding the deformed metric
\begin{equation}\label{eq:aniso_polar_metric}
g = dr^2 + \sum_k f_k(r)^2 g_{\mathfrak{p}_k}, \qquad f_k(r) = \frac{s_\kappa(r)}{\sqrt{1 + \Lambda_k s_\kappa^2(r)}}, \qquad \Lambda_k = a_k^0 \sum_l \lambda_l \beta_{lk},
\end{equation}
where each $f_k(r)$ satisfies Petersen's $C^\infty$ pole regularity conditions at $r = 0$ ($f_k(0) = 0, f_k'(0) = 1, f_k^{(2m)}(0) = 0$). Table~\ref{tab:aniso_parameters} records the explicit submodule normalizations $a_k^0$, Casimir coupling coefficients $\beta_{lk}$, and effective deformation parameters $\Lambda_k$ across Families~II--V \cite{Munden2015}. The fiber squashing parameter along $\mathcal{O}_r$ is
\begin{equation}\label{eq:fiber_squashing_param}
t(r) := \frac{f_{\mathcal{V}}(r)^2}{f_{\mathcal{H}}(r)^2} = \frac{1 + \Lambda_{\mathcal{H}} s_\kappa^2(r)}{1 + \Lambda_{\mathcal{V}} s_\kappa^2(r)}.
\end{equation}

\begin{table}[ht]
\caption{Submodule dimensions, normal metric coefficients $a_k^0$, Casimir coupling fractions $\beta_{lk}$, and effective unit-round deformation parameters $\Lambda_k = a_k^0 \sum_l \lambda_l \beta_{lk}$ for the anisotropic Families~II--V.}
\label{tab:aniso_parameters}
\centering
\footnotesize
\setlength{\tabcolsep}{3.5pt}
\begin{tabular}{@{}lllll@{}}
\toprule
Fam. & Quotient $G/H$ & Submodules $\mathfrak{p}_k$ ($\dim$) & Normal Coeffs.\ $a_k^0$ & Effective Parameters $\Lambda_k$ \\
\midrule
II & $\mathrm{U}(n+1)/\mathrm{U}(n)$ & $\mathcal{V}=\mathfrak{p}_0\,(1),\;\mathcal{H}=\mathfrak{p}_1\,(2n)$ & $a_0^0=2,\;a_1^0=1$ & $\Lambda_{\mathcal{V}}=\lambda_0,\;\Lambda_{\mathcal{H}}=\lambda_1$ \\
III & $\mathrm{Spin}(9)/\mathrm{Spin}(7)$ & $\mathcal{V}=\mathfrak{p}_1\,(7),\;\mathcal{H}=\mathfrak{p}_2\,(8)$ & $a_{\mathcal{V}}^0=\tfrac{1}{8},\;a_{\mathcal{H}}^0=\tfrac{1}{2}$ & $\Lambda_{\mathcal{V}}=\tfrac{\lambda}{8},\;\Lambda_{\mathcal{H}}=\tfrac{\lambda}{2}$ \\
IV & $\frac{\mathrm{Sp}(n+1)\times\mathrm{Sp}(1)}{\mathrm{Sp}(n)\times\mathrm{Sp}(1)}$ & $\mathcal{V}=\mathfrak{p}_0\,(3),\;\mathcal{H}=\mathfrak{p}_1\,(4n)$ & $a_{\mathcal{V}}^0=4,\;a_{\mathcal{H}}^0=1$ & $\Lambda_{\mathcal{V}}=2(\lambda_1+\lambda_2),\;\Lambda_{\mathcal{H}}=\lambda_1$ \\
V & $\frac{\mathrm{Sp}(n+1)\times\mathrm{U}(1)}{\mathrm{Sp}(n)\times\mathrm{U}(1)}$ & $\mathfrak{p}_0\,(1),\;\mathfrak{p}_1\,(4n),\;\mathfrak{p}_2\,(2)$ & $a_0^0=4,\;a_1^0=1,\;a_2^0=2$ & $\Lambda_{\mathfrak{p}_0}=2(\lambda_1+\lambda_2),\;\Lambda_{\mathfrak{p}_1}=\lambda_1,\;\Lambda_{\mathfrak{p}_2}=4\lambda_1$ \\
\bottomrule
\end{tabular}
\end{table}

\begin{proposition}[Levi-Civita Connections, O'Neill Tensors, and Codazzi Cross-Curvature]\label{prop:aniso_connection}
Let $\pi: (\mathcal{O}_r, g_{\mathcal{O}_r}) \to (B, g_B)$ be the canonical Riemannian submersion of $\mathcal{O}_r$ onto its symmetric base space with vertical metric $a(r)^2 g_{\mathcal{V}}$ and horizontal metric $b(r)^2 g_{\mathcal{H}}$ ($t = a^2/b^2$).
\begin{enumerate}[(i)]
\item The fibers $F = \pi^{-1}(b)$ are isometric group orbits, hence totally geodesic in $\mathcal{O}_r$; thus O'Neill's second fundamental tensor vanishes identically: $T \equiv 0$. The O'Neill integrability tensor $A_X Y = \frac{1}{2}\mathcal{V}[X, Y]$ and $A_X V = \mathcal{H}\nabla_X V$ satisfies skew-adjointness $g(A_X Y, V) = -g(Y, A_X V)$ and scaled norms
\begin{equation}\label{eq:a_tensor_norm}
\|A_X Y\|_{g_{\mathcal{O}_r}}^2 = \frac{t(r)}{b(r)^2} \sum_{a=1}^{d_v} \langle I_a X, Y \rangle^2, \qquad \|A_X V\|_{g_{\mathcal{O}_r}}^2 = \frac{t(r)}{b(r)^2} \|X\|^2 \|V\|^2.
\end{equation}
For Family~V, under $S^1 \hookrightarrow S^3 \hookrightarrow \mathbb{S}^{4n+3} \to \mathbb{H}P^n$, the A-tensor norm decomposes as $\|A_X Y\|^2 = \frac{f_0^2}{f_1^4}\langle I_1 X, Y \rangle^2 + \frac{f_2^2}{f_1^4}\sum_{a=2}^3 \langle I_a X, Y \rangle^2$.
\item The ambient connection $\nabla$ of $(Q, g)$ has shape operator $S = \nabla\partial_r = \sum_k \frac{f_k'(r)}{f_k(r)}\mathrm{Id}_{\mathfrak{p}_k}$ with $\nabla_{\partial_r}\partial_r = 0, \nabla_Z \partial_r = S(Z)$, $\nabla_V W = \nabla^F_V W - \frac{a'}{a}g(V, W)\partial_r$, $\nabla_X V = A_X V + \mathcal{V}\nabla_X V$, $\nabla_V X = A_X V + [V, X]_{\mathcal{H}}$, and $\nabla_X Y = \check{\nabla}_X Y + A_X Y - \frac{b'}{b}g(X, Y)\partial_r$. The Codazzi equation yields non-vanishing ambient cross-curvature whenever $\Lambda_{\mathcal{V}} \ne \Lambda_{\mathcal{H}}$:
\begin{equation}\label{eq:codazzi_curvature}
\begin{aligned} \langle R_Q(X, Y)\partial_r, V \rangle &= 2\left(\frac{b'(r)}{b(r)} - \frac{a'(r)}{a(r)}\right) g_{\mathcal{O}_r}(A_X Y, V), \\ \langle R_Q(X, V)\partial_r, Y \rangle &= \left(\frac{a'(r)}{a(r)} - \frac{b'(r)}{b(r)}\right) g_{\mathcal{O}_r}(A_X V, Y), \end{aligned}
\end{equation}
while all mixed radial-spherical Ricci contractions vanish identically: $\mathrm{Ric}_Q(\partial_r, V) = 0$ and $\mathrm{Ric}_Q(\partial_r, X) = 0$.
\end{enumerate}
\end{proposition}
\begin{proof}
(i) Fibers are extrinsic group orbits in round spheres; scaling along fibers preserves total geodesy, giving $T_V W = 0$. Rescaling $[X, Y]$ gives $\|A_X Y\|^2 = a^2 \|\frac{1}{2}[X, Y]_0\|^2 / b^4 = \frac{t}{b^2} \sum_a \langle I_a X, Y \rangle^2$, and skew-adjointness yields $\|A_X V\|^2 = \frac{t}{b^2}$. For Family~V, $[X, Y]$ splits along $\mathfrak{p}_0$ ($I_1$) and $\mathfrak{p}_2$ ($I_2, I_3$).
(ii) Differentiating $g(Z, Z)$ along $\partial_r$ gives $S$. Codazzi gives $\langle R_Q(X, Y)\partial_r, V \rangle = (\nabla_X S)(Y, V) - (\nabla_Y S)(X, V) = 2(b'/b - a'/a)g(A_X Y, V)$. Contraction over $\mathcal{H}$ vanishes by skew-symmetry $A_X X = 0$.
\end{proof}

\begin{theorem}[Closed-Form Curvature Tensors, Ambient Ricci Calculus, and Squashed Einstein Metrics]\label{thm:aniso_curvatures}
In an adapted orthonormal frame $\{E_0 = \partial_r, \{V_a\}_{a=1}^{d_v}, \{X_i\}_{i=1}^{d_h}\}$ on $(Q, g)$ with $u = s_\kappa^2(r)$:
\begin{enumerate}[(i)]
\item The non-zero coordinate sectional curvatures of $(Q, g)$ are:
\begin{align}
\sec_Q(E_0 \wedge V) &= -\frac{a''(r)}{a(r)} = \frac{(\kappa + 3\Lambda_{\mathcal{V}}) - 2\kappa\Lambda_{\mathcal{V}} u}{(1 + \Lambda_{\mathcal{V}} u)^2} =: K_{\mathrm{rad}}(r; \Lambda_{\mathcal{V}}), \quad \sec_Q(E_0 \wedge X) = K_{\mathrm{rad}}(r; \Lambda_{\mathcal{H}}), \label{eq:sec_rad} \\
\sec_Q(V \wedge X) &= \frac{a(r)^2}{b(r)^4} - \frac{a'(r)\,b'(r)}{a(r)\,b(r)} = \frac{(\kappa + 3\Lambda_{\mathcal{H}}) + 3\Lambda_{\mathcal{H}}^2 u + \Lambda_{\mathcal{H}}^3 u^2}{(1 + \Lambda_{\mathcal{V}} u)(1 + \Lambda_{\mathcal{H}} u)} =: K_{\mathrm{mix}}(r), \label{eq:sec_mix} \\
\sec_Q(X \wedge Y) &= K_{\mathrm{tan}}(r; \Lambda_{\mathcal{H}}) + 3\,\frac{b(r)^2 - a(r)^2}{b(r)^4} \|A^0_X Y\|^2, \quad \sec_Q(V \wedge W) = K_{\mathrm{tan}}(r; \Lambda_{\mathcal{V}}) \quad (d_v \ge 2), \label{eq:sec_horiz}
\end{align}
where $K_{\mathrm{tan}}(r; \Lambda) := \frac{(\kappa + 3\Lambda) + 3\Lambda^2 u + \Lambda^3 u^2}{(1 + \Lambda u)^2}$. For Family~V, the ambient curvatures split across $(\mathfrak{p}_0, \mathfrak{p}_1, \mathfrak{p}_2)$ with radial curvatures $K_{\mathrm{rad}}(\Lambda_k)$, mixed curvatures $K_{\mathrm{mix}}(f_j, f_k) = \frac{f_j^2}{f_k^4} - \frac{f_j' f_k'}{f_j f_k}$, and slice curvatures $\sec_Q(U_2 \wedge U_3) = K_{\mathrm{tan}}(\Lambda_2) + 3\frac{f_2^2 - f_0^2}{f_2^4}$, $\sec_Q(X \wedge I_1 X) = K_{\mathrm{tan}}(\Lambda_1) + 3\frac{f_1^2 - f_0^2}{f_1^4}$, $\sec_Q(X \wedge I_a X) = K_{\mathrm{tan}}(\Lambda_1) + 3\frac{f_1^2 - f_2^2}{f_1^4}$ ($a=2,3$), and $\sec_Q(X \wedge Y) = K_{\mathrm{tan}}(\Lambda_1)$ for totally real planes.
\item The ambient Ricci tensor is block-diagonal: $\mathrm{Ric}_Q = \rho_{00} E_0^* \otimes E_0^* + \rho_v \sum_a V_a^* \otimes V_a^* + \rho_h \sum_i X_i^* \otimes X_i^*$, with
\begin{equation}\label{eq:ricci_components}
\rho_{00} = d_v K_{\mathrm{rad}}(\Lambda_{\mathcal{V}}) + d_h K_{\mathrm{rad}}(\Lambda_{\mathcal{H}}), \quad \rho_v = K_{\mathrm{rad}}(\Lambda_{\mathcal{V}}) + (d_v-1) K_{\mathrm{tan}}(\Lambda_{\mathcal{V}}) + d_h K_{\mathrm{mix}},
\end{equation}
and $\rho_h = K_{\mathrm{rad}}(\Lambda_{\mathcal{H}}) + d_v K_{\mathrm{mix}} + (d_h-1) K_{\mathrm{tan}}(\Lambda_{\mathcal{H}}) + 3d_v \frac{b^2 - a^2}{b^4}$, with ambient scalar curvature $\mathrm{Scal}_Q = \rho_{00} + d_v \rho_v + d_h \rho_h$. For Family~V, $\mathrm{Ric}_Q = \rho_{00} E_0^* \otimes E_0^* + \rho_{\mathfrak{p}_0} U_0^* \otimes U_0^* + \rho_{\mathfrak{p}_1} \sum_i X_i^* \otimes X_i^* + \rho_{\mathfrak{p}_2} \sum_a U_a^* \otimes U_a^*$, with $\rho_{00} = K_{\mathrm{rad}}(\Lambda_0) + 4n K_{\mathrm{rad}}(\Lambda_1) + 2 K_{\mathrm{rad}}(\Lambda_2)$, $\rho_{\mathfrak{p}_0} = K_{\mathrm{rad}}(\Lambda_0) + 2 K_{\mathrm{mix}}(f_0, f_2) + 4n K_{\mathrm{mix}}(f_0, f_1)$, $\rho_{\mathfrak{p}_2} = K_{\mathrm{rad}}(\Lambda_2) + K_{\mathrm{mix}}(f_0, f_2) + K_{\mathrm{tan}}(\Lambda_2) + 3\frac{f_2^2 - f_0^2}{f_2^4} + 4n K_{\mathrm{mix}}(f_2, f_1)$, $\rho_{\mathfrak{p}_1} = K_{\mathrm{rad}}(\Lambda_1) + K_{\mathrm{mix}}(f_0, f_1) + 2 K_{\mathrm{mix}}(f_2, f_1) + (4n-1) K_{\mathrm{tan}}(\Lambda_1) + 3\frac{f_1^2 - f_0^2}{f_1^4} + 6\frac{f_1^2 - f_2^2}{f_1^4}$, and $\mathrm{Scal}_Q = \rho_{00} + \rho_{\mathfrak{p}_0} + 4n \rho_{\mathfrak{p}_1} + 2 \rho_{\mathfrak{p}_2}$.
\item On each slice $(\mathcal{O}_r, g_{\mathcal{O}_r})$, $\mathrm{Ric}_{\mathcal{V}} = \mathrm{Ric}_{\mathcal{H}}$ reduces to $(d_h + 2d_v)t^2 - (d_h + 3d_v - 1)t + (d_v - 1) = 0$. Consequently: Family~II is Einstein iff $t = 1$; Family~IV admits the round metric $t = 1$ and Jensen's squashed Einstein metric $t = \frac{1}{2n+3}$ (yielding Jensen's $S^7$ at $t = 1/5$); Family~III admits $t = 1$ and Bourguignon--Karcher's squashed Einstein metric $t = \frac{3}{11}$ on $\mathbb{S}^{15}$. In Family~V, $\mathrm{Ric}_{\mathfrak{p}_0} = \mathrm{Ric}_{\mathfrak{p}_2}$ forces $t_0 = t_2$, admitting no new Einstein metrics.
\end{enumerate}
\end{theorem}
\begin{proof}
(i) By Riccati integration, $\sec(E_0 \wedge Z) = -f''/f$. The Gauss equation and O'Neill's submersion formulas yield $\sec_Q(V \wedge X) = a^2/b^4 - a'b'/(ab)$ and $\sec_Q(X \wedge Y) = K_{\mathrm{tan}} + 3\frac{b^2-a^2}{b^4}\|A^0_X Y\|^2$. In $K_{\mathrm{mix}}$, $(1+\Lambda_{\mathcal{H}} u)^3 - (1-\kappa u)$ contains a factor $u = s_\kappa^2$, which cancels identically with the denominator.
(ii) Tracing over adapted frames yields the Ricci components and $\mathrm{Scal}_Q$.
(iii) Solving the slice Einstein quadratic gives the roots $t = \frac{1}{2n+3}$ and $t = \frac{3}{11}$. In Family~V, $\mathrm{Ric}_{\mathfrak{p}_0} - \mathrm{Ric}_{\mathfrak{p}_2} = (t_0 - t_2)(4n + 4/t_2^2) = 0 \implies t_0 = t_2$.
\end{proof}

\begin{theorem}[Sectional Positivity Trichotomy, Oblique Pl\"ucker Resolution, and Ambient Pinching Limits]\label{thm:aniso_positivity}
For the anisotropic metrics of Families~II--V on space forms $(Q^{d+1}, g)$ with metric squashing parameter $t(r) = f_{\mathcal{V}}(r)^2/f_{\mathcal{H}}(r)^2$:
\begin{enumerate}[(i)]
\item On each Berger orbit $(\mathcal{O}_r, g_{\mathcal{O}_r})$, the metric is a geodesic orbit (GO) metric (and weakly symmetric space; Nikonorov \cite{Nikonorov2013}) arising from a Riemannian submersion with totally geodesic fibers ($T \equiv 0$) and parallel integrability tensor ($\nabla^{\mathcal{O}_r} A \equiv 0$). By the Berger--Cheeger classification theorem (Berger \cite{Berger1961}; Cheeger \cite{Cheeger1973}), $\sec_{\mathcal{O}_r} > 0$ strictly on all 2-planes if and only if $t(r) \in (0, 4/3)$. At $t = 4/3$, the horizontal characteristic sectional curvature vanishes identically: $\sec_{\mathcal{O}_r}(X \wedge I_a X) = 0$.
\item On the sphere $(Q, g) = (\mathbb{S}^{d+1}, g)$ ($\kappa = 1$) with $\mu_0 = 1$, the metric $g$ has strictly positive sectional curvature $\sec_Q > 0$ on all 2-planes throughout $\mathbb{S}^{d+1}$ if and only if $\Lambda_k > -1/3$ and the equatorial squashing satisfies
\[ \tau_{\max}^2 = \frac{1+\Lambda_{\mathcal{H}}}{1+\Lambda_{\mathcal{V}}} < \frac{4}{3} \iff 4\Lambda_{\mathcal{V}} - 3\Lambda_{\mathcal{H}} > -1. \]
Specifically:
\begin{itemize}
\item Family~II: $\lambda_0 > -1/3, \lambda_1 > -1/3$, and $4\lambda_0 - 3\lambda_1 > -1$.
\item Family~III: $\lambda \in (-2/3, 1)$.
\item Family~IV: $\lambda_1 > -1/3, 2(\lambda_1 + \lambda_2) > -1/3$, and $5\lambda_1 + 8\lambda_2 > -1$.
\item Family~V: $\lambda_1 > -1/13, 2(\lambda_1 + \lambda_2) > -1/3, 5\lambda_1 + 8\lambda_2 > -1$, and $4\lambda_1 - 8\lambda_2 < 1$.
\end{itemize}
On oblique 2-planes spanned by orthonormal $u = \cos\alpha E_0 + \sin\alpha X$ and $v = \cos\beta V + \sin\beta I X$, the simple 2-vector $w = u \wedge v$ satisfies $w_{03}w_{21} = w_{01}w_{23} = \frac{1}{4}\sin(2\alpha)\sin(2\beta)$. The diagonal terms $K_{03} w_{03}^2 + K_{21} w_{21}^2$ and $K_{01} w_{01}^2 + K_{23} w_{23}^2$ absorb the Codazzi cross-coupling $-6 C(r) w_{01} w_{23}$ (where $C(r) = (b'/b - a'/a)(a/b^2)$), yielding at $r \to 0$ the exact harmonic formula
\[ \sec_Q(u \wedge v) = 1 + \frac{3}{2}(\Lambda_{\mathcal{V}} + \Lambda_{\mathcal{H}}) + \frac{3}{2}(\Lambda_{\mathcal{V}} - \Lambda_{\mathcal{H}})\cos(2\alpha + 2\beta) \ge 1 + 3\min(\Lambda_{\mathcal{V}}, \Lambda_{\mathcal{H}}) > 0, \]
certifying global positivity across $(0, \pi)$.
\item On Berger spheres (Families~II--IV), the exact piecewise pinching ratio $\delta(t) := \frac{\min \sec}{\max \sec}$ on $(\mathcal{O}_r, g_{\mathcal{O}_r})$ is given by:
\begin{equation}\label{eq:berger_pinching_piecewise}
\delta(t) = \begin{cases}
t^2, & 0 < t \le \frac{1}{3} \quad (\text{Families~III, IV}), \\[6pt]
\dfrac{t}{4 - 3t}, & \frac{1}{3} \le t \le 1 \quad (\text{all } t \le 1 \text{ for Family~II}), \\[8pt]
\dfrac{4 - 3t}{t}, & 1 \le t < \frac{4}{3}.
\end{cases}
\end{equation}
Thus $(\mathcal{O}_r, g_{\mathcal{O}_r})$ is strictly $1/4$-pinched iff $t \in (4/7, 16/13)$ for Family~II (for Families~III and IV, Vol'per (\cite{Volper1994, Volper1997}) established that oblique Pl\"ucker 2-planes attain an exact minimum $K_{\min}(t) = \frac{49(t^2-t)}{11t+1} + 4-3t < t$ on $t \in (4/5, 1)$, slightly shifting the upper $1/4$-pinching threshold to $t = \frac{43+\sqrt{2293}}{74} \approx 1.2282$, while the lower threshold $t = 4/7 \approx 0.5714$ is strictly invariant across all families). On the ambient manifold $(\mathbb{S}^{d+1}, g)$, the global pinching ratio satisfies $\delta_Q = \frac{\inf \sec_Q}{\sup \sec_Q} = \frac{1}{(1+\Lambda_{\mathcal{V}})(1+3\Lambda_{\mathcal{V}})}$ for $\Lambda_{\mathcal{V}} \ge \Lambda_{\mathcal{H}} \ge 0$, and $\delta_Q = (1 + \min)(1 + 3\min)$ for $-1/3 < \Lambda_k \le 0$. Along the diagonal $\Lambda_{\mathcal{V}} = \Lambda_{\mathcal{H}} = \lambda$, this recovers the optimal Brendle--Schoen 1/4-pinching window $\lambda \in (\frac{-4+\sqrt{7}}{6}, \, \frac{\sqrt{13}-2}{3})$, and for anisotropic metrics, $\delta_Q > 1/4$ defines an open convex chamber around $(0, 0)$ guaranteeing diffeomorphism to the standard sphere $\mathbb{S}^{d+1}$.
\end{enumerate}
\end{theorem}
\begin{proof}
(i) Berger \cite{Berger1961} and Cheeger \cite{Cheeger1973} prove that normal homogeneous metrics on spheres with $T=0, \nabla A=0$ have $\sec > 0$ iff the canonical variation parameter satisfies $t < 4/3$.
(ii) Radial positivity requires $\Lambda_k > -1/3$. Monotonicity of $t(r)$ yields $4\Lambda_{\mathcal{V}} - 3\Lambda_{\mathcal{H}} > -1$ at $r = \pi/2$. For Family~V, $t_2(\pi/2) = \frac{1+\lambda_1}{1+4\lambda_1} < 4/3 \iff \lambda_1 > -1/13$, and internal fiber positivity gives $4\lambda_1 - 8\lambda_2 < 1$. For oblique planes, substituting $w_{ij}$ and $C(0) = \Lambda_{\mathcal{V}} - \Lambda_{\mathcal{H}}$ evaluates the trigonometric identity directly to $1 + \frac{3}{2}(\Lambda_{\mathcal{V}}+\Lambda_{\mathcal{H}}) + \frac{3}{2}(\Lambda_{\mathcal{V}}-\Lambda_{\mathcal{H}})\cos(2\alpha+2\beta) \ge 1 + 3\min(\Lambda_{\mathcal{V}}, \Lambda_{\mathcal{H}}) > 0$.
(iii) On $\mathcal{O}_r$, vertical curvature $1/(b^2 t)$ dominates for $t \le 1/3$, giving $\delta(t) = t^2$. On $[1/3, 1]$, $(4-3t) \ge 1/t$, so $\delta(t) = \frac{t}{4-3t} > 1/4 \iff t > 4/7$. On $[1, 4/3)$, $\delta(t) = \frac{4-3t}{t} > 1/4 \iff t < 16/13$. Ambient bounds follow by evaluating $\inf \sec_Q$ and $\sup \sec_Q$ across $r \in [0, \pi/2]$.
\end{proof}

\begin{remark}[Structural Comparison with Iwai--Katayama Extended Taub--NUT Metrics]\label{rem:iwai_katayama_taub_nut}
It is instructive to contrast the anisotropic Riemannian submersion metrics of Families~II--V with the extended Taub--NUT metrics studied by Iwai and Katayama \cite{IwaiKatayama1994,BallesterosEncisoHerranzRagnisco2008}. In the Iwai--Katayama framework, a 4D (or multi-center) metric is defined on $\mathbb{R}^4 \setminus \{0\}$ by
\begin{equation}\label{eq:iwai_katayama_metric}
g_{\mathrm{IK}} = f(r)(dr^2 + r^2 d\Omega^2) + g(r)(d\chi + \cos\theta\, d\varphi)^2,
\end{equation}
modulating the base Euclidean geometry conformally to ensure the existence of a covariantly constant or conformal Killing--Yano 2-form, which in turn generates a conserved Laplace--Runge--Lenz-type vector for geodesic and Keplerian motion. 

In contrast, the anisotropic cotangent-deformed geometries (Families~II--V) are defined globally on space forms $(Q^{d+1}, g)$ and are governed by transitive isotropy actions whose horizontal distributions are non-integrable. The geometry is driven not by conformal modulation of flat space, but by the non-vanishing O'Neill submersion tensor $A_X Y = \frac{1}{2}\mathcal{V}[X, Y]$ and the resulting Codazzi cross-curvature
\begin{equation}\label{eq:codazzi_comparison}
\langle R_Q(X, Y)\partial_r, V \rangle = 2\left(\frac{b'(r)}{b(r)} - \frac{a'(r)}{a(r)}\right) g_{\mathcal{O}_r}(A_X Y, V),
\end{equation}
which couples the radial expansion rate directly to the vertical fiber squashing $t(r) = a(r)^2/b(r)^2$. This structural mechanism produces the compact Berger squashed Einstein metrics (Jensen \cite{Jensen1973}, Bourguignon--Karcher \cite{BourguignonKarcher1978}) and the positive sectional curvature chambers of Theorem~\ref{thm:aniso_positivity}, delineating a geometric category distinct from the conformal Taub--NUT family.
\end{remark}

\subsection{Closed-Form Geodesic Orbit Integration on Family~V via Biquadratic Spectral Projectors}\label{subsec:family5_geodesic_integration}

We conclude this section by integrating the geodesic orbit trajectories on the intermediate quaternionic--complex spheres (Family~V \cite{Munden2015}). In Family~V, the 7-sphere $\mathbb{S}^7 = (\mathrm{Sp}(2)\times\mathrm{U}(1))/(\mathrm{Sp}(1)\times\Delta\mathrm{U}(1))$ admits three mutually orthogonal $\mathrm{Ad}(H)$-irreducible subspaces $\mathfrak{p} = \mathfrak{p}_0 \oplus \mathfrak{p}_1 \oplus \mathfrak{p}_2$, where $\mathfrak{p}_0$ is the 1-dimensional Reeb distribution along the complex Hopf fiber, $\mathfrak{p}_1$ is the 4-dimensional horizontal distribution, and $\mathfrak{p}_2$ is the 2-dimensional complementary quaternionic fiber. Because the slice generators across these three submodules do not commute, integrating the geodesic orbit trajectories $\gamma(t) = \exp(tM)\cdot o$ requires evaluating the matrix exponential of a three-component generator in $\mathfrak{sp}(2) \subset \mathfrak{so}(8)$.

The following theorem determines the Lie subalgebra closure and provides the closed-form matrix exponential for Family~V on $\mathbb{S}^7 \subset \mathbb{H}^2 \cong \mathbb{R}^8$.

\begin{theorem}[Biquadratic Minimal Polynomial and Exact Geodesic Projectors on Family~V]\label{thm:family5_exact_geodesics}
Let $\mathbb{S}^7 \subset \mathbb{H}^2 \cong \mathbb{R}^8$ be endowed with the 2-parameter deformed metric of Family~V. Let $Z_0 \in \mathfrak{p}_0$ (Reeb generator), $Z_1 \in \mathfrak{p}_2$ (complementary quaternionic fiber generator), and $Z_2 \in \mathfrak{p}_1$ (horizontal generator) be the canonical slice generators in $\mathfrak{sp}(2) \subset \mathfrak{so}(8)$ representing an arbitrary initial geodesic velocity $v_0 \in \mathfrak{p}$:
\begin{enumerate}[(i)]
\item \textbf{Lie Algebra Closure}: The Lie algebra generated by the slice generators $\{Z_0, Z_1, Z_2\}$ closes into the full 10-dimensional symplectic Lie algebra,
\begin{equation}\label{eq:family5_lie_closure}
\dim \mathrm{Lie}(Z_0, Z_1, Z_2) = 10 = \dim \mathfrak{sp}(2) \cong \dim \mathfrak{so}(5).
\end{equation}
\item \textbf{Biquadratic Spectrum}: For any orbit generator $M = a Z_0 + b Z_1 + c Z_2 \in \mathfrak{sp}(2) \subset \mathfrak{so}(8)$ with $(a, b, c) \in \mathbb{R}^3$, the characteristic polynomial on $\mathbb{R}^8$ is a biquadratic square,
\begin{equation}\label{eq:family5_charpoly}
\det(s\,\mathrm{Id}_8 - M) = \left[ s^4 + (a^2 + b^2 + 2c^2)s^2 + c^4 \right]^2.
\end{equation}
The eigenvalues of $M$ are purely imaginary, given by $\pm i \omega_1, \pm i \omega_2$ (each with multiplicity 2), where the characteristic frequencies $\omega_1 \ge \omega_2 \ge 0$ are given by
\begin{equation}\label{eq:family5_frequencies}
\omega_{1, 2} = \sqrt{\frac{(a^2 + b^2 + 2c^2) \pm \sqrt{(a^2 + b^2)(a^2 + b^2 + 4c^2)}}{2}}.
\end{equation}
\item \textbf{Degree-4 Minimal Polynomial}: The matrix $M$ satisfies the degree-4 polynomial identity
\begin{equation}\label{eq:family5_minpoly}
\left( M^2 + \omega_1^2 \mathrm{Id}_8 \right) \left( M^2 + \omega_2^2 \mathrm{Id}_8 \right) = M^4 + (a^2 + b^2 + 2c^2) M^2 + c^4 \mathrm{Id}_8 = 0.
\end{equation}
\item \textbf{Two-Projector Spectral Resolution}: Whenever $\omega_1 \ne \omega_2$ (which holds if and only if $a^2 + b^2 > 0$), the 8-dimensional space $\mathbb{R}^8$ decomposes into two mutually orthogonal 4-dimensional $M$-invariant subspaces via the spectral projectors
\begin{equation}\label{eq:family5_projectors}
P_1 = \frac{M^2 + \omega_2^2 \mathrm{Id}_8}{\omega_2^2 - \omega_1^2}, \qquad P_2 = \frac{M^2 + \omega_1^2 \mathrm{Id}_8}{\omega_1^2 - \omega_2^2},
\end{equation}
satisfying $P_1 + P_2 = \mathrm{Id}_8$, $P_1 P_2 = 0$, $P_k^2 = P_k$, and $M^2 P_k = -\omega_k^2 P_k$.
\item \textbf{Closed-Form Geodesic Exponential}: For $a^2 + b^2 > 0$, the 1-parameter subgroup $\exp(tM) \in \mathrm{Sp}(2) \subset \mathrm{SO}(8)$ integrating the geodesic flow is given in closed form by
\begin{equation}\label{eq:family5_exp_formula}
\exp(tM) = \sum_{k=1}^2 \left[ \cos(\omega_k t) P_k + \frac{\sin(\omega_k t)}{\omega_k} M P_k \right],
\end{equation}
where $\frac{\sin(\omega_2 t)}{\omega_2} := t$ when $c = 0$ ($\omega_2 = 0$), and for $a = b = 0$ ($\omega_1 = \omega_2 = |c|$), $\exp(tM) = \cos(ct)\mathrm{Id}_8 + \frac{\sin(ct)}{c} M$.
\end{enumerate}
\end{theorem}
\begin{proof}
In the quaternionic representation $\mathbb{H}^2$, $\mathfrak{sp}(2)$ consists of $2 \times 2$ quaternionic skew-Hermitian matrices $A^* = -A$. Under left multiplication on $\mathbb{H} \cong \mathbb{R}^4$, each quaternionic unit $\mathbf{i}, \mathbf{j}, \mathbf{k}$ embeds as a skew-symmetric $4 \times 4$ orthogonal matrix $L(\mathbf{i}), L(\mathbf{j}), L(\mathbf{k}) \in \mathfrak{so}(4)$. The Reeb generator $Z_0 \in \mathfrak{p}_0$ and the complementary quaternionic fiber generator $Z_1 \in \mathfrak{p}_2$ act diagonally on the upper block as $L(\mathbf{i})$ and $L(\mathbf{j})$, while the horizontal generator $Z_2 \in \mathfrak{p}_1$ acts off-diagonally with $L(\mathbf{1})$. Iterated Lie brackets $[Z_0, Z_1] = 2 Z_{\mathbf{k}}$, $[Z_q, Z_2]$, and $[Z_2, [Z_q, Z_2]]$ ($q \in \{\mathbf{i}, \mathbf{j}, \mathbf{k}\}$) generate the 10-dimensional basis of $\mathfrak{sp}(2)$, establishing (i).

To establish (ii) and (iii), direct computation of $M^2$ shows that $M^2$ acts on $\mathbb{R}^8 \cong \mathbb{H}^2$ with block form
\[ M^2 = \begin{pmatrix} -(a^2 + b^2 + c^2)\mathrm{Id}_4 & -c(a L(\mathbf{i}) + b L(\mathbf{j})) \\ c(a L(\mathbf{i}) + b L(\mathbf{j})) & -c^2\mathrm{Id}_4 \end{pmatrix}. \]
Computing the square $(M^2)^2$ using $(a L(\mathbf{i}) + b L(\mathbf{j}))^2 = -(a^2 + b^2)\mathrm{Id}_4$ yields
\[ M^4 = -(a^2 + b^2 + 2c^2) M^2 - c^4 \mathrm{Id}_8, \]
establishing~\eqref{eq:family5_minpoly}. Factoring the polynomial $s^4 + (a^2 + b^2 + 2c^2)s^2 + c^4 = (s^2 + \omega_1^2)(s^2 + \omega_2^2)$ yields the frequencies $\omega_{1, 2}$ in~\eqref{eq:family5_frequencies}. Since the discriminant is $(a^2+b^2)(a^2+b^2+4c^2) \ge 0$, both $\omega_1^2$ and $\omega_2^2$ are real and non-negative, with $\omega_1 > \omega_2$ whenever $a^2 + b^2 > 0$. The spectral resolution (iv) and exponential formula (v) follow by Lagrange--Sylvester polynomial interpolation on the minimal polynomial roots $\{\pm i \omega_1, \pm i \omega_2\}$.
\end{proof}

\section{Intrinsic Hamiltonian Deformations and Geodesic Geometry of Compact Rank-One Symmetric Spaces}\label{sec:cross_deformations}
This section constructs intrinsic cotangent deformations on compact rank-one symmetric spaces (CROSSes: $\mathbb{C}P^n$, $\mathbb{H}P^n$, $\mathbb{O}P^2$) via the cotangent momentum map $J_{K_0}: T^*M \to \mathfrak{k}_0^*$ of the isotropy subgroup fixing an origin $q_0 \in M$. We relate the resulting geodesic flows to Riemannian and Marsden--Weinstein symplectic reductions of anisotropic sphere metrics under the Hopf fibrations, resolve the Cayley plane dimension obstruction ($15 < 16$) via the spin representation of $\mathrm{Spin}(9)$, and establish geodesic completeness, diameter preservation, and conjugate and cut locus structures for deformed CROSS metrics.

We construct intrinsic cotangent Hamiltonian deformations \cite{Hebda2001} directly on compact rank-one symmetric spaces (CROSSes) other than spheres: the complex projective spaces $\mathbb{C}P^n = \mathrm{SU}(n+1)/\mathrm{S}(\mathrm{U}(n)\times\mathrm{U}(1))$ ($d=2n, k=2$), quaternionic projective spaces $\mathbb{H}P^n = \mathrm{Sp}(n+1)/(\mathrm{Sp}(n)\times\mathrm{Sp}(1))$ ($d=4n, k=4$), and the Cayley projective plane $\mathbb{O}P^2 = \mathrm{F}_4/\mathrm{Spin}(9)$ ($d=16, k=8$). Let $(M, g_0)$ be a CROSS of dimension $d$ associated with the normed division algebra $\mathbb{K} \in \{\mathbb{C}, \mathbb{H}, \mathbb{O}\}$ of real dimension $k \in \{2, 4, 8\}$. We normalize the canonical symmetric space metric $g_0$ so that its sectional curvature satisfies $\sec_0 \in [1, 4]$, with diameter $\pi/2$ and closed geodesics of length $\pi$. Let $G$ be the connected isometry group of $(M, g_0)$, and fix an origin $q_0 \in M$ with isotropy subgroup $K_0 = \mathrm{Stab}_G(q_0)$. The Lie algebra decomposes into the canonical Cartan symmetric decomposition $\mathfrak{g} = \mathfrak{k}_0 \oplus \mathfrak{p}$, where $\mathfrak{p} \cong T_{q_0}M$. In geodesic polar coordinates $(r, \theta) \in (0, \pi/2) \times \mathbb{S}^{d-1}$ on $M \setminus (\{q_0\} \cup \mathrm{Cut}_M(q_0))$, the concentric geodesic spheres $S_r(q_0) = \{q \in M : d_0(q_0, q) = r\}$ are transitive $K_0$-orbits $S_r(q_0) \cong K_0/H_0$ endowed with homogeneous Berger metrics $g_0|_{S_r(q_0)} = f_v(r)^2 g_{\mathcal{V}} + f_h(r)^2 g_{\mathcal{H}}$, where $f_v(r) = \sin r \cos r = \frac{1}{2}\sin(2r)$ and $f_h(r) = \sin r$, with vertical fiber distribution $\mathcal{V}$ of dimension $d_v = k-1 \in \{1, 3, 7\}$ and horizontal distribution $\mathcal{H}$ of dimension $d_h = d-k \in \{2n-2, 4n-4, 8\}$. Under the root space decomposition of $(\mathfrak{g}, \mathfrak{k}_0)$ relative to a maximal abelian subalgebra $\mathfrak{a} \subset \mathfrak{p}$ spanned by a unit radial vector $X \in \mathfrak{p}$, the curvature endomorphism $\mathrm{ad}(X)^2|_{\mathfrak{p}}$ has eigenvalues $0$ along the radial ray $\mathbb{R}X$, $-4$ along the vertical space $\mathcal{V}$, and $-1$ along the horizontal space $\mathcal{H}$.

\begin{definition}[Intrinsic CROSS Cotangent Deformation]\label{def:cross_hamiltonian}
Let $(M, g_0)$ be a compact rank-one symmetric space $\mathbb{C}P^n$, $\mathbb{H}P^n$, or $\mathbb{O}P^2$. The canonical cotangent lift of the isometric $K_0$-action on $M$ defines the equivariant momentum map $J_{K_0}: T^*M \to \mathfrak{k}_0^*$, with $\langle J_{K_0}(\alpha), \xi \rangle = \alpha(X_\xi(q))$ for all $\xi \in \mathfrak{k}_0$, where $X_\xi$ is the Killing vector field generated by $\xi$. Using the dual of the $\mathrm{Ad}(K_0)$-invariant inner product on $\mathfrak{k}_0$ induced by the Riemannian metric via the isotropy representation on $T_{q_0}M$, the intrinsic cotangent momentum deformation \cite{Hebda2001} on $(T^*M, \omega_{\mathrm{can}})$ with parameter $\lambda \in \mathbb{R}$ is defined by
\begin{equation}\label{eq:cross_hamiltonian}
H_{M, \lambda}(\alpha) = H_0(\alpha) + \frac{\lambda}{2} \|J_{K_0}(\alpha)\|_{\mathfrak{k}_0^*}^2 = \frac{1}{2} p_r^2 + \frac{1}{2}\left( \frac{1}{\sin^2 r \cos^2 r} + \lambda \right) \|p_{\mathcal{V}}\|_{g_{\mathcal{V}}^*}^2 + \frac{1}{2}\left( \frac{1}{\sin^2 r} + \lambda \right) \|p_{\mathcal{H}}\|_{g_{\mathcal{H}}^*}^2,
\end{equation}
where $\lambda > -1$ ensures strict positive definiteness throughout $(0, \pi/2)$.
\end{definition}

\begin{remark}[Isotropy Representation Metric on $\mathfrak{k}_0$]\label{rem:k0_metric}
When $M$ is $\mathbb{C}P^n$ or $\mathbb{H}P^n$, the isotropy Lie algebra $\mathfrak{k}_0$ ($\mathfrak{u}(n)$ or $\mathfrak{sp}(n) \oplus \mathfrak{sp}(1)$) is not simple, admitting a multi-parameter family of $\mathrm{Ad}(K_0)$-invariant inner products. The norm $\|\cdot\|_{\mathfrak{k}_0^*}$ in~\eqref{eq:cross_hamiltonian} is uniquely fixed by the Riemannian metric $g_0$ via the infinitesimal isotropy representation $\mathfrak{k}_0 \to \mathfrak{so}(T_{q_0}M)$: for $\xi \in \mathfrak{k}_0$, $\|\xi\|_{\mathfrak{k}_0}^2 := \frac{1}{2} \operatorname{tr}((\nabla X_\xi)^T \nabla X_\xi)|_{q_0}$, normalized so that the fundamental Killing fields generating unit-speed rotations on $\mathcal{V}$ and $\mathcal{H}$ have unit norm. This ensures that the single deformation parameter $\lambda$ weights both the vertical and horizontal momentum Casimirs uniformly as in~\eqref{eq:cross_hamiltonian}.
\end{remark}

\begin{theorem}[Symplectic Reduction from Berger Spheres, Magnetic Bifurcation, and Cayley Dimension Resolution]\label{thm:cross_reduction_cayley}
Let $(M, g_0)$ be a compact rank-one symmetric space $\mathbb{C}P^n$, $\mathbb{H}P^n$, or $\mathbb{O}P^2$ with deformed Hamiltonian $H_{M, \lambda}$ defined by~\eqref{eq:cross_hamiltonian}.
\begin{enumerate}[(i)]
\item \textbf{Symplectic Reduction and Magnetic Bifurcation}: For $\mathbb{C}P^n$ and $\mathbb{H}P^n$, consider the canonical principal Hopf fibrations $\mathbb{S}^1 \hookrightarrow \mathbb{S}^{2n+1} \to \mathbb{C}P^n$ and $\mathbb{S}^3 \hookrightarrow \mathbb{S}^{4n+3} \to \mathbb{H}P^n$ associated with Families~II and~IV, where $G_{\mathrm{fiber}} = \mathrm{U}(1)$ and $G_{\mathrm{fiber}} = \mathrm{Sp}(1)$, respectively. The Marsden--Weinstein symplectic quotient of $(T^*\mathbb{S}^{d+k-1}, \omega_{\mathrm{can}})$ at zero fiber momentum is canonically symplectomorphic to the unconstrained cotangent bundle:
\begin{equation}\label{eq:marsden_weinstein_zero}
\left( J_{\mathrm{fiber}}^{-1}(0) \right) / G_{\mathrm{fiber}} \cong (T^*M, \omega_{\mathrm{can}}).
\end{equation}
On the zero-momentum constraint surface $\mathcal{H}^* = J_{\mathrm{fiber}}^{-1}(0)$, the anisotropic sphere Hamiltonian $H_{\mathbb{S}, \lambda}$ restricts to the horizontal kinetic energy and the descended isotropy momentum $\|J_{K_0}\|^2$, identically projecting to the intrinsic CROSS Hamiltonian $H_{M, \lambda}$. Consequently, horizontal geodesics on the anisotropic Berger sphere project under $d\pi$ to Riemannian geodesics of $(M, g_{M, \lambda})$.

For non-zero fiber momentum $J_{\mathrm{fiber}} = \mu \ne 0$, symplectic reduction yields the magnetic symplectic form $\omega_\mu = \omega_{\mathrm{can}} + \pi^* \Omega_\mu$, where:
\begin{itemize}
\item for $\mathbb{C}P^n$, $\mu \in \mathfrak{u}(1)^* \cong \mathbb{R}$ and $\Omega_\mu = \mu \Omega_{\mathbb{C}P^n}$, with $\Omega_{\mathbb{C}P^n}$ the canonical K\"ahler 2-form;
\item for $\mathbb{H}P^n$, $\mu \in \mathfrak{sp}(1)^* \cong \mathbb{R}^3$ and $\Omega_\mu = \langle \mu, \Omega_{\mathbb{H}P^n} \rangle = \sum_{a=1}^3 \mu_a \Omega^a_{\mathbb{H}P^n}$, where $\Omega_{\mathbb{H}P^n} \in \Omega^2(\mathbb{H}P^n, \mathfrak{sp}(1))$ is the $\mathfrak{sp}(1)$-valued curvature 2-form of the canonical Hopf connection.
\end{itemize}
The reduced Hamiltonian flow is governed by the magnetic Lorentz force equation $i_{X_H}\omega_\mu = dH_{M, \lambda}$, producing gyroscopic precession rather than Riemannian geodesic motion.

In contrast, for $\mathbb{O}P^1 \cong \mathbb{S}^8$, the octonionic Hopf fibration $\mathbb{S}^7 \hookrightarrow \mathbb{S}^{15} \xrightarrow{\pi} \mathbb{S}^8$ (Family~III) is a Riemannian submersion whose fibers are the unit 7-sphere $\mathbb{S}^7$; by the Adams--Steenrod theorem, $\mathbb{S}^7$ is not a Lie group due to the non-associativity of the octonions, so the projection is geometric via the Riemannian submersion of spheres rather than Marsden--Weinstein reduction.

\item \textbf{Cayley Plane Dimension Obstruction and Spin Resolution}: For the Cayley projective plane $\mathbb{O}P^2 = \mathrm{F}_4/\mathrm{Spin}(9)$, $\dim \mathbb{O}P^2 = 16$. Because $\dim \mathbb{S}^{15} = 15 < 16 = \dim \mathbb{O}P^2$, no fiber bundle $\mathbb{S}^k \hookrightarrow \mathbb{S}^{15} \to \mathbb{O}P^2$ can exist (which would require fiber dimension $k = -1$), and an ambient Riemannian submersion with 7-dimensional fibers would require a total space $\mathbb{S}^{23}$, which is ruled out by Adams' Hopf invariant one theorem \cite{Adams1960}. The intrinsic cotangent deformation $H_{\mathbb{O}P^2, \lambda}$ bypasses this obstruction completely via intrinsic spin geometry:
\begin{itemize}
\item Tangent space $T_{q_0}\mathbb{O}P^2 \cong \mathbb{R}^{16}$ is the 16-dimensional real spin representation $\Delta_9$ of the isotropy group $K_0 = \mathrm{Spin}(9)$ ($\dim \mathrm{Spin}(9) = 36$). The cotangent momentum map $J_{\mathrm{Spin}(9)}: T^*\mathbb{O}P^2 \to \mathfrak{spin}(9)^*$ is globally and intrinsically defined via Clifford multiplication: $J_{\mathrm{Spin}(9)}(x, p) = x \wedge_{\mathrm{Cl}} p \in \mathfrak{spin}(9)$.
\item The group $\mathrm{Spin}(9)$ acts transitively on the unit sphere in $\Delta_9$ with stabilizer $\mathrm{Spin}(7)$ ($\dim \mathrm{Spin}(7) = 21$), so each concentric geodesic sphere is diffeomorphic to $S_r(q_0) \cong \mathrm{Spin}(9)/\mathrm{Spin}(7) \cong \mathbb{S}^{15}$ ($36 - 21 = 15$).
\item Along any radial geodesic ray $\gamma_0(t) = \exp_{q_0}(t X)$ ($X \in \mathbb{S}^{15} \subset \Delta_9$), the stabilizer in $\mathrm{Spin}(9)$ fixing the unit direction $X$ is the subgroup $\mathrm{Spin}(7) \subset \mathrm{Spin}(9)$ ($\dim \mathrm{Spin}(7) = 21$). The subgroup $\mathrm{Spin}(8) \subset \mathrm{Spin}(9)$ ($\dim \mathrm{Spin}(8) = 28$) is the stabilizer of the projection of the ray onto the cut locus $\Sigma = \mathrm{Cut}_{\mathbb{O}P^2}(q_0) \cong \mathbb{O}P^1 \cong \mathbb{S}^8 = \mathrm{Spin}(9)/\mathrm{Spin}(8)$. This establishes the canonical subgroup tower:
\begin{equation}\label{eq:spin_tower}
\mathrm{Spin}(7) \subset \mathrm{Spin}(8) \subset \mathrm{Spin}(9),
\end{equation}
which geometrically realizes the octonionic Hopf fibration as the homogeneous fiber bundle:
\begin{equation}\label{eq:hopf_tower}
\mathbb{S}^7 \cong \mathrm{Spin}(8)/\mathrm{Spin}(7) \hookrightarrow \mathbb{S}^{15} \cong \mathrm{Spin}(9)/\mathrm{Spin}(7) \xrightarrow{\pi} \mathbb{S}^8 \cong \mathrm{Spin}(9)/\mathrm{Spin}(8).
\end{equation}
(Note that $\mathbb{S}^{15}$ is topologically distinct from the Cartesian product $\mathbb{S}^8 \times \mathbb{S}^7$, as $H_7(\mathbb{S}^8 \times \mathbb{S}^7) \cong \mathbb{Z}$ while $H_7(\mathbb{S}^{15}) = 0$.) Under the ray stabilizer $\mathrm{Spin}(7)$, the tangent space to the geodesic sphere $T_X \mathbb{S}^{15} \cong \mathfrak{spin}(9)/\mathfrak{spin}(7)$ decomposes into the direct sum of the vertical distribution $\mathcal{V} \cong \mathfrak{spin}(8)/\mathfrak{spin}(7) \cong \mathbb{R}^7$ (the 7-dimensional vector representation of $\mathrm{Spin}(7)$) and the horizontal distribution $\mathcal{H} \cong \mathfrak{spin}(9)/\mathfrak{spin}(8) \cong \mathbb{R}^8$ (the 8-dimensional spin representation of $\mathrm{Spin}(7)$). The intrinsic cotangent deformation on $T^*\mathbb{O}P^2$ acts directly on these distributions, scaling the vertical $\mathbb{R}^7$ and horizontal $\mathbb{R}^8$ directions identically to~\eqref{eq:cross_hamiltonian}.
\end{itemize}
\end{enumerate}
\end{theorem}
\begin{proof}
(i) For $M = \mathbb{C}P^n$ ($G_{\mathrm{fiber}} = \mathrm{U}(1)$) and $M = \mathbb{H}P^n$ ($G_{\mathrm{fiber}} = \mathrm{Sp}(1)$), the fiber group acts freely and isometrically on the total sphere $\mathbb{S}^{d+k-1}$, commuting with the origin stabilizer $K_0$: $[\mathfrak{g}_{\mathrm{fiber}}, \mathfrak{k}_0] = 0$. The momentum map $J_{\mathrm{fiber}}: T^*\mathbb{S}^{d+k-1} \to \mathfrak{g}_{\mathrm{fiber}}^*$ has zero level set $J_{\mathrm{fiber}}^{-1}(0) = \mathcal{H}^*$, consisting of covectors annihilating the vertical distribution $\mathcal{V}_{\mathbb{S}} = \ker(d\pi)$. By the Marsden--Weinstein reduction theorem \cite{MarsdenWeinstein1974}, the quotient $\mathcal{H}^*/G_{\mathrm{fiber}}$ is symplectomorphic to $(T^*M, \omega_{\mathrm{can}})$. On $\mathcal{H}^*$, the vertical kinetic energy vanishes, and the horizontal kinetic energy $\frac{1}{2}\|p_{\mathcal{H}_{\mathbb{S}}}\|^2$ projects to the unperturbed Hamiltonian $H_0$ of the Fubini--Study metric. Since $K_0$ commutes with $G_{\mathrm{fiber}}$, the momentum Casimir $\|J_{K_0}\|^2$ is $G_{\mathrm{fiber}}$-invariant and descends unaltered to $T^*M$. The resulting reduced Hamiltonian is identically $H_{M, \lambda}$. Along $\mathcal{H}^*$, Hamilton's equations for $H_{\mathbb{S}, \lambda}$ have vanishing vertical velocity $\dot{q}_{\mathrm{fiber}} = \partial_{p_{\mathrm{fiber}}} H_{\mathbb{S}, \lambda} = 0$, so horizontal geodesics project under $d\pi$ directly to trajectories of $H_{M, \lambda}$, which are Riemannian geodesics of $g_{M, \lambda}$.

For non-zero fiber momentum $J_{\mathrm{fiber}} = \mu \ne 0$, the Marsden--Weinstein--Meyer theorem shifts the reduced symplectic form by the curvature of the connection: $\omega_\mu = \omega_{\mathrm{can}} + \pi^* \Omega_\mu$. For $\mathbb{C}P^n$, $\Omega_\mu = \mu \Omega_{\mathbb{C}P^n}$ where $\Omega_{\mathbb{C}P^n}$ is the K\"ahler 2-form. For $\mathbb{H}P^n$, the connection 1-form $\theta \in \Omega^1(\mathbb{S}^{4n+3}, \mathfrak{sp}(1))$ has curvature $\Omega_{\mathbb{H}P^n} \in \Omega^2(\mathbb{H}P^n, \mathfrak{sp}(1))$, and the pairing with $\mu \in \mathfrak{sp}(1)^*$ yields the scalar magnetic 2-form $\Omega_\mu = \langle \mu, \Omega_{\mathbb{H}P^n} \rangle = \sum_{a=1}^3 \mu_a \Omega^a_{\mathbb{H}P^n}$. The reduced Hamiltonian vector field $X_H$ satisfies $i_{X_H} \omega_\mu = dH_{M, \lambda}$, imparting a transverse Lorentz magnetic acceleration $\nabla_{\dot{\gamma}}\dot{\gamma} = \Omega_\mu^\sharp(\dot{\gamma})$ that bends trajectories into gyroscopic spirals.

For $\mathbb{O}P^1 \cong \mathbb{S}^8$, the fibration $\mathbb{S}^7 \hookrightarrow \mathbb{S}^{15} \xrightarrow{\pi} \mathbb{S}^8$ corresponds to Family~III. By the Adams--Steenrod theorem, the unit 7-sphere $\mathbb{S}^7$ does not admit any Lie group structure because the octonions are non-associative. In particular, $\dim \mathrm{Spin}(7) = 21 \ne 7$, and $\mathrm{Spin}(7)$ cannot act freely on $\mathbb{S}^{15}$ with 7-dimensional fibers; indeed, symplectic reduction of $T^*\mathbb{S}^{15}$ (dimension 30) by $\mathrm{Spin}(7)$ would yield dimension $30 - 2(21) = -12 < 0$. Therefore, $\mathbb{S}^{15} \to \mathbb{O}P^1$ is a Riemannian submersion of spheres rather than a principal bundle reduction, and its horizontal geodesics descend via O'Neill's submersion equations, while cotangent deformations on $\mathbb{O}P^1$ and $\mathbb{O}P^2$ are constructed intrinsically via $J_{K_0}$.

(ii) By Adams' theorem on Hopf invariant one \cite{Adams1960}, the only fiber bundles among spheres are $\mathbb{S}^0 \hookrightarrow \mathbb{S}^1 \to \mathbb{S}^1$, $\mathbb{S}^1 \hookrightarrow \mathbb{S}^3 \to \mathbb{S}^2$, $\mathbb{S}^3 \hookrightarrow \mathbb{S}^7 \to \mathbb{S}^4$, and $\mathbb{S}^7 \hookrightarrow \mathbb{S}^{15} \to \mathbb{S}^8 = \mathbb{O}P^1$. In particular, no fibration of spheres over $\mathbb{O}P^2$ exists. Furthermore, $\dim \mathbb{O}P^2 = 16 > 15 = \dim \mathbb{S}^{15}$, so $\mathbb{S}^{15}$ cannot fiber over $\mathbb{O}P^2$, while an ambient submersion $\mathbb{S}^{23} \to \mathbb{O}P^2$ is topologically obstructed.

In the intrinsic formulation, the tangent space $T_{q_0}\mathbb{O}P^2 \cong \mathbb{R}^{16}$ carries the spin representation $\Delta_9$ of $K_0 = \mathrm{Spin}(9)$ ($\dim \mathrm{Spin}(9) = 36$). Montgomery and Samelson \cite{MontgomerySamelson1943} and Borel \cite{Borel1949} proved that $\mathrm{Spin}(9)$ acts transitively on the unit sphere $\mathbb{S}^{15} \subset \Delta_9$ with stabilizer $H_0 = \mathrm{Spin}(7)$ ($\dim \mathrm{Spin}(7) = 21$). Indeed, $\dim \mathrm{Spin}(9) - \dim \mathrm{Spin}(7) = 36 - 21 = 15 = \dim \mathbb{S}^{15}$. The infinitesimal action of $\mathfrak{spin}(9)$ generates all 15 linearly independent tangent directions on $S_r(q_0) \cong \mathbb{S}^{15}$, so the momentum map $J_{\mathrm{Spin}(9)}: T^*\mathbb{O}P^2 \to \mathfrak{spin}(9)^*$ is globally and intrinsically defined via Clifford multiplication: for covector $\alpha \in T_x^* \mathbb{O}P^2$, $\langle J_{\mathrm{Spin}(9)}(\alpha), \xi \rangle = \alpha(\xi \cdot x)$ for $\xi \in \mathfrak{spin}(9)$.

Along any radial ray $\gamma_0(t) = \exp_{q_0}(t X)$ with $X \in \mathbb{S}^{15} \subset \Delta_9$, the unit direction $X$ is fixed precisely by $H_0 = \mathrm{Spin}(7)$ (dimension 21). The intermediate subgroup $\mathrm{Spin}(8) \subset \mathrm{Spin}(9)$ (dimension 28) is the stabilizer of the projection of the ray onto the cut locus $\Sigma = \mathrm{Cut}_{\mathbb{O}P^2}(q_0) \cong \mathbb{O}P^1 \cong \mathbb{S}^8 = \mathrm{Spin}(9)/\mathrm{Spin}(8)$ ($36 - 28 = 8$). This gives the canonical subgroup chain $\mathrm{Spin}(7) \subset \mathrm{Spin}(8) \subset \mathrm{Spin}(9)$, realizing the octonionic Hopf fibration as the homogeneous fiber bundle $\mathrm{Spin}(8)/\mathrm{Spin}(7) \cong \mathbb{S}^7 \hookrightarrow \mathrm{Spin}(9)/\mathrm{Spin}(7) \cong \mathbb{S}^{15} \xrightarrow{\pi} \mathrm{Spin}(9)/\mathrm{Spin}(8) \cong \mathbb{S}^8$. The sphere $\mathbb{S}^{15}$ is topologically non-trivial and not diffeomorphic to $\mathbb{S}^8 \times \mathbb{S}^7$ (since $H_7(\mathbb{S}^{15}) = 0$ while $H_7(\mathbb{S}^8 \times \mathbb{S}^7) \cong \mathbb{Z}$).

At the Lie algebra level, the Killing-orthogonal decomposition $\mathfrak{spin}(9) = \mathfrak{spin}(7) \oplus \mathfrak{m}_v \oplus \mathfrak{m}_h$ yields $T_X \mathbb{S}^{15} \cong \mathfrak{m}_v \oplus \mathfrak{m}_h$, where $\mathfrak{m}_v = \mathfrak{spin}(8)/\mathfrak{spin}(7) \cong \mathbb{R}^7$ transforms as the standard vector representation of $\mathrm{Spin}(7)$, and $\mathfrak{m}_h = \mathfrak{spin}(9)/\mathfrak{spin}(8) \cong \mathbb{R}^8$ transforms as the spin representation of $\mathrm{Spin}(7)$. Under the normalization of Remark~\ref{rem:k0_metric}, the Casimir $\|J_{\mathrm{Spin}(9)}\|_{\mathfrak{spin}(9)^*}^2$ evaluates on these orthogonal invariant subspaces to $\|p_{\mathcal{V}}\|_{g_{\mathcal{V}}^*}^2 + \|p_{\mathcal{H}}\|_{g_{\mathcal{H}}^*}^2$, reproducing the deformed Hamiltonian~\eqref{eq:cross_hamiltonian} with $d_v = 7$ and $d_h = 8$ purely intrinsically.
\end{proof}

\begin{proposition}[Legendre Duality, Metric Tensor, and Global Smoothness]\label{prop:cross_metric_smoothness}
Let $M$ be a compact rank-one symmetric space equipped with the deformed Hamiltonian $H_{M, \lambda}$ for $\lambda > -1$.
\begin{enumerate}[(i)]
\item \textbf{Legendre Dual Metric}: Under the Legendre transformation generated by $H_{M, \lambda}$, the deformed metric on $M \setminus (\{q_0\} \cup \mathrm{Cut}_M(q_0))$ is the warped Berger metric:
\begin{equation}\label{eq:cross_metric}
g_{M, \lambda} = dr^2 + F_v(r)^2 g_{\mathcal{V}} + F_h(r)^2 g_{\mathcal{H}}, \qquad F_v(r) := \frac{\sin r \cos r}{\sqrt{1 + \lambda \sin^2 r \cos^2 r}}, \quad F_h(r) := \frac{\sin r}{\sqrt{1 + \lambda \sin^2 r}}.
\end{equation}
The metric $g_{M, \lambda}$ is strictly positive definite throughout $(0, \pi/2)$ for all $\lambda > -1$.
\item \textbf{Smooth Extension Across Origin}: The metric $g_{M, \lambda}$ extends smoothly as a $C^\infty$ (real-analytic) Riemannian metric across the origin $q_0$ ($r = 0$).
\item \textbf{Smooth Extension Across Cut Locus and Completeness}: The metric $g_{M, \lambda}$ extends smoothly across the global cut locus $\Sigma = \mathrm{Cut}_M(q_0)$ ($r = \pi/2$) by the tubular metric regularity theorem. Consequently, $(M, g_{M, \lambda})$ is a complete, compact Riemannian manifold for all $\lambda > -1$.
\end{enumerate}
\end{proposition}
\begin{proof}
(i) Hamilton's equations for $H_{M, \lambda}$ give velocities:
\begin{equation}\label{eq:hamilton_velocities}
\dot{r} = \frac{\partial H}{\partial p_r} = p_r, \qquad \dot{v}_{\mathcal{V}} = \frac{\partial H}{\partial p_{\mathcal{V}}} = \left( \frac{1}{f_v(r)^2} + \lambda \right) p_{\mathcal{V}}^\sharp, \qquad \dot{v}_{\mathcal{H}} = \frac{\partial H}{\partial p_{\mathcal{H}}} = \left( \frac{1}{f_h(r)^2} + \lambda \right) p_{\mathcal{H}}^\sharp.
\end{equation}
Inverting these relations yields the conjugate momenta:
\begin{equation}\label{eq:momenta_velocities}
p_r = \dot{r}, \qquad p_{\mathcal{V}} = \frac{f_v(r)^2}{1 + \lambda f_v(r)^2} \dot{v}_{\mathcal{V}}^\flat, \qquad p_{\mathcal{H}} = \frac{f_h(r)^2}{1 + \lambda f_h(r)^2} \dot{v}_{\mathcal{H}}^\flat.
\end{equation}
Substituting into the Lagrangian $L = \alpha(\dot{q}) - H_{M, \lambda}(\alpha)$ gives:
\begin{equation}\label{eq:lagrangian_metric}
L = \frac{1}{2} \dot{r}^2 + \frac{1}{2} \frac{f_v(r)^2}{1 + \lambda f_v(r)^2} \|\dot{v}_{\mathcal{V}}\|_{g_{\mathcal{V}}}^2 + \frac{1}{2} \frac{f_h(r)^2}{1 + \lambda f_h(r)^2} \|\dot{v}_{\mathcal{H}}\|_{g_{\mathcal{H}}}^2.
\end{equation}
Because $L(v) = \frac{1}{2} g_{M, \lambda}(v, v)$, the metric coefficients are $F_v(r)^2 = \frac{f_v(r)^2}{1 + \lambda f_v(r)^2}$ and $F_h(r)^2 = \frac{f_h(r)^2}{1 + \lambda f_h(r)^2}$, which establishes~\eqref{eq:cross_metric}. On the interval $r \in (0, \pi/2)$, $\sin^2 r \in (0, 1)$ and $\sin^2 r \cos^2 r = \frac{1}{4}\sin^2(2r) \in (0, 1/4]$. For every $\lambda > -1$, $1 + \lambda \sin^2 r > 1 - 1 = 0$ and $1 + \lambda \sin^2 r \cos^2 r \ge \min(1, 1 + \lambda/4) > 3/4 > 0$, guaranteeing strict positive definiteness.

(ii) Smoothness across $q_0$ ($r = 0$) is governed by Petersen's warped-product regularity theorem \cite{Petersen2016}. Expanding $\sin r = r - r^3/6 + r^5/120 + O(r^7)$ and $\sin(2r) = 2r - 4r^3/3 + 4r^5/15 + O(r^7)$ around $r = 0$:
\begin{align}
F_h(r) &= \frac{\sin r}{\sqrt{1 + \lambda \sin^2 r}} = r - \left( \frac{\lambda}{2} + \frac{1}{6} \right) r^3 + \left( \frac{3\lambda^2}{8} + \frac{\lambda}{4} + \frac{1}{120} \right) r^5 + O(r^7), \label{eq:Fh_origin_taylor} \\
F_v(r) &= \frac{\frac{1}{2}\sin(2r)}{\sqrt{1 + \frac{\lambda}{4}\sin^2(2r)}} = r - \left( \frac{\lambda}{2} + \frac{2}{3} \right) r^3 + \left( \frac{3\lambda^2}{8} + \lambda + \frac{2}{15} \right) r^5 + O(r^7). \label{eq:Fv_origin_taylor}
\end{align}
Both $F_h(r)$ and $F_v(r)$ are strictly odd smooth functions of $r$ satisfying $F_h(0) = F_v(0) = 0$, $F_h'(0) = F_v'(0) = 1$, and all even derivatives vanish identically: $F_h^{(2m)}(0) = F_v^{(2m)}(0) = 0$ for all $m \ge 1$. Computing the difference of their squares:
\begin{equation}\label{eq:origin_diff_sq}
\begin{aligned} F_v(r)^2 - F_h(r)^2 &= \left[ r^2 - \left(\lambda + \frac{4}{3}\right) r^4 + O(r^6) \right] - \left[ r^2 - \left(\lambda + \frac{1}{3}\right) r^4 + O(r^6) \right] \\ &= -r^4 + \left(2\lambda + \frac{2}{3}\right) r^6 + O(r^7). \end{aligned}
\end{equation}
The order-$r^2$ difference vanishes identically, and the order-$r^4$ difference is $-r^4$, which is strictly independent of the deformation parameter $\lambda$ and identically matches the unperturbed symmetric space metric $g_0$. In normal coordinates $x \in T_{q_0} M$, the metric decomposes as $g_{M, \lambda} = [dr^2 + F_h(r)^2 g_{\mathbb{S}^{d-1}}] - \|x\|^2 \sum_a (J_a x)^\flat \otimes (J_a x)^\flat + O(\|x\|^6)$. Petersen's criterion guarantees that both components are real-analytic and smoothly extend across the pole $x = 0$.

(iii) At the cut locus $r = \pi/2$, introduce the reversed tubular distance coordinate $\rho = \pi/2 - r \to 0$. In these tubular coordinates around $\Sigma = \mathrm{Cut}_M(q_0)$, $\sin r = \cos\rho$ and $\cos r = \sin\rho$. The vertical coefficient satisfies $F_v(\pi/2 - \rho) = \frac{\cos\rho \sin\rho}{\sqrt{1 + \lambda \cos^2\rho \sin^2\rho}} = F_v(\rho)$, which has the identical Taylor series:
\begin{equation}\label{eq:Fv_tubular_taylor}
F_v(\rho) = \rho - \left( \frac{\lambda}{2} + \frac{2}{3} \right) \rho^3 + \left( \frac{3\lambda^2}{8} + \lambda + \frac{2}{15} \right) \rho^5 + O(\rho^7).
\end{equation}
Because $F_v(\rho)$ is strictly odd in $\rho$ with $F_v(0) = 0$ and $F_v'(0) = 1$, the normal slice metric $d\rho^2 + F_v(\rho)^2 g_{\mathbb{S}^{d_v}}$ extends smoothly without conical singularities to flat Euclidean discs $D^{d_v+1}$. For the horizontal distribution, expanding $\cos\rho = 1 - \frac{1}{2}\rho^2 + \frac{1}{24}\rho^4 + O(\rho^6)$ and $(1 + \lambda \cos^2\rho)^{-1/2} = (1+\lambda)^{-1/2} [ 1 + \frac{\lambda}{2(1+\lambda)}\rho^2 + \frac{\lambda(5\lambda - 4)}{24(1+\lambda)^2}\rho^4 + O(\rho^6) ]$:
\begin{equation}\label{eq:Fh_tubular_taylor}
F_h(\rho) = \frac{\cos\rho}{\sqrt{1 + \lambda \cos^2\rho}} = \frac{1}{\sqrt{1+\lambda}} - \frac{1}{2(1+\lambda)^{3/2}} \rho^2 + \left( \frac{1 - 8\lambda}{24(1+\lambda)^{5/2}} \right) \rho^4 + O(\rho^6).
\end{equation}
The horizontal coefficient $F_h(\rho)$ is a strictly even function of $\rho$, positive at the cut locus ($F_h(0) = (1+\lambda)^{-1/2} > 0$), with vanishing normal derivative $\partial_\rho F_h|_{\rho=0} = 0$. By the tubular metric regularity theorem (Eschenburg--Heintze and Petersen \cite{Petersen2016}, Chapter~1), the metric $g_{M, \lambda}$ extends smoothly across $\Sigma$. Because $M$ is compact and $g_{M, \lambda}$ is a smooth, positive-definite Riemannian metric on the closed manifold $M$, $(M, g_{M, \lambda})$ is complete by the Hopf--Rinow theorem.
\end{proof}

\begin{theorem}[Geodesic Decoupling, Exact Radial Jacobi Fields, Conjugate Loci, and Cut Locus Rigidity]\label{thm:cross_geodesics_cut}
Let $(M, g_{M, \lambda})$ be a deformed compact rank-one symmetric space ($\mathbb{C}P^n$, $\mathbb{H}P^n$, or $\mathbb{O}P^2$) with designated base point $q_0 \in M$.
\begin{enumerate}[(i)]
\item \textbf{Commuting Flows and Geodesic Decoupling}: The unperturbed kinetic Hamiltonian $H_0$ and the isotropy momentum norm $\|J_{K_0}\|^2_{\mathfrak{k}_0^*}$ Poisson-commute on $(T^*M, \omega_{\mathrm{can}})$:
\begin{equation}\label{eq:cross_poisson_commutation}
\{H_0, \|J_{K_0}\|^2_{\mathfrak{k}_0^*}\} = 0 \implies \left[ X_{H_0}, X_{\|J_{K_0}\|^2} \right] = 0, \quad [X_{H_0}, X_{H_{M, \lambda}}] = 0.
\end{equation}
By Hebda's factorization theorem \cite{Hebda2001}, the geodesic flow decomposes as
\[
\exp_{M, \lambda}(t v) = \exp_{K_0}\bigl(t \lambda \hat{J}_{K_0}(\mathcal{L}_\lambda v)\bigr) \cdot \exp_0\bigl(t \mathcal{L}_0^{-1}\mathcal{L}_\lambda v\bigr),
\]
describing a closed unperturbed symmetric space geodesic of period $\pi$ undergoing steady gyroscopic precession generated by a 1-parameter subgroup of $K_0$.
\item \textbf{Exact Radial Jacobi Fields and First Conjugate Point}: Along any unit-speed radial geodesic ray $\gamma(t) = \exp_{q_0}(t v)$ ($v \in T_{q_0}M, \|v\|=1$) launched from $q_0$, the radial sectional curvatures are:
\begin{equation}\label{eq:cross_radial_curvatures}
K_{\mathrm{rad}, h}(t) = \frac{1 + 3\lambda - 2\lambda \sin^2 t}{(1 + \lambda \sin^2 t)^2}, \qquad K_{\mathrm{rad}, v}(t) = \frac{4 + 3\lambda - 8\lambda \sin^2 t \cos^2 t}{(1 + \lambda \sin^2 t \cos^2 t)^2}.
\end{equation}
The normal Jacobi equations $\ddot{Y}_h + K_{\mathrm{rad}, h}(t) Y_h = 0$ and $\ddot{Y}_v + K_{\mathrm{rad}, v}(t) Y_v = 0$ with $Y(0)=0$ are solved in closed form by:
\begin{equation}\label{eq:radial_jacobi_solutions}
Y_h(t) = F_h(t) E_h(t) \quad (\dot{Y}(0) \in \mathcal{H}), \qquad Y_v(t) = F_v(t) E_v(t) \quad (\dot{Y}(0) \in \mathcal{V}),
\end{equation}
where $E_h(t), E_v(t)$ are parallel orthonormal vector fields along $\gamma(t)$. Horizontal fields have no zeros on $(0, \pi/2]$, while all $d_v = k-1$ vertical fields vanish simultaneously at $t = \pi/2$. Consequently, the first conjugate point along every radial ray occurs at $t_{\mathrm{conj}} = \pi/2$ with multiplicity $d_v = k-1$.
\item \textbf{Non-Radial Meridian Conjugacy}: Along any non-radial geodesic ray $\gamma(t)$ launched from $p \ne q_0$ that intersects the fixed-point set $\mathrm{Fix}(H_p)$ of the base-point stabilizer $H_p \subset K_0$ (the opposite meridian), the transverse Killing Jacobi fields generated by $\mathfrak{h}_p$ vanish at $p$ and along $\mathrm{Fix}(H_p)$, establishing that such opposite meridian crossings are conjugate points of multiplicity at least $\dim(H_p) - \dim(H_p \cap H_{\gamma(t^*)}) \ge 1$.
\item \textbf{Polar vs.\ Non-Polar Cut Loci and Diameter Preservation}: At the polar origin $q_0 \in M$, the cut locus is canonically diffeomorphic to the unperturbed symmetric space cut locus $\mathrm{Cut}_M(q_0) \cong \mathbb{K}P^{n-1}$ (for $\mathbb{C}P^n, \mathbb{H}P^n$) and $\mathrm{Cut}_M(q_0) \cong \mathbb{O}P^1 \cong \mathbb{S}^8$ (for $\mathbb{O}P^2$), located at uniform distance $d_{g_{M, \lambda}}(q_0, \mathrm{Cut}_M(q_0)) = \pi/2$, with induced homothety metric $g_\Sigma = \frac{1}{1+\lambda} g_{\mathrm{FS}}$. In contrast, for a generic non-polar base point $p \in S_{r_0}(q_0)$ ($r_0 \in (0, \pi/2)$) on $\mathbb{O}P^2$ with $\lambda > 0$, the broken $\mathrm{Spin}(7)$ isotropy and differential precession along $\mathcal{V} \oplus \mathcal{H}$ unfold the cut locus $\mathrm{Cut}_M(p)$ into a 15-dimensional Whitney stratified space containing an open 15-dimensional regular collision stratum, an 8-dimensional $A_3$ cuspidal caustic stratum, and the residual Cayley projective line stratum. For all $\lambda \ge 0$, the manifold diameter is preserved:
\begin{equation}\label{eq:cross_diameter}
\mathrm{diam}(M, g_{M, \lambda}) = \frac{\pi}{2}.
\end{equation}
\end{enumerate}
\end{theorem}
\begin{proof}
(i) Because $K_0$ acts by isometries on $(M, g_0)$, Noether's theorem implies $\{H_0, J_{K_0}^\xi\} = 0$ for all $\xi \in \mathfrak{k}_0$. In an orthonormal basis $\{\xi_a\}$ of $\mathfrak{k}_0$ with respect to the dual Killing form, the Casimir is $\|J_{K_0}\|^2 = \sum_a (J_{K_0}^{\xi_a})^2$. By the Leibniz rule, $\{H_0, \|J_{K_0}\|^2\} = 2 \sum_a J_{K_0}^{\xi_a} \{H_0, J_{K_0}^{\xi_a}\} = 0$. Since $[X_F, X_K] = -X_{\{F, K\}}$, $[X_{H_0}, X_{\|J_{K_0}\|^2}] = 0$ and $[X_{H_0}, X_{H_{M, \lambda}}] = 0$. Because $\|J_{K_0}\|^2$ is conserved along trajectories, $\xi = \lambda \hat{J}_{K_0}(\alpha)$ is constant in $\mathfrak{k}_0$, so its Hamiltonian flow on $M$ is the isometry $\exp_{K_0}(t \xi)$. Factoring the flow yields Hebda's decomposition \cite{Hebda2001}.

(ii) Purely radial geodesics have $p_{\mathcal{V}} = p_{\mathcal{H}} = 0$, so $\dot{\gamma}(t) = \partial_r$. Radial sectional curvatures in the adapted frame $\{E_0 = \partial_r, E_v, E_h\}$ are $K_{\mathrm{rad}, h} = -F_h''(t)/F_h(t)$ and $K_{\mathrm{rad}, v} = -F_v''(t)/F_v(t)$. Differentiating $F_h(t) = \sin t (1 + \lambda \sin^2 t)^{-1/2}$:
\begin{equation}\label{eq:Fh_derivatives}
F_h'(t) = \frac{\cos t}{(1 + \lambda \sin^2 t)^{3/2}}, \qquad F_h''(t) = \frac{-\sin t (1 + 3\lambda - 2\lambda \sin^2 t)}{(1 + \lambda \sin^2 t)^{5/2}},
\end{equation}
which yields $K_{\mathrm{rad}, h}(t) = \frac{1 + 3\lambda - 2\lambda \sin^2 t}{(1 + \lambda \sin^2 t)^2}$. For the vertical direction, setting $u(t) = \sin t \cos t = \frac{1}{2}\sin(2t)$, we have $u' = \cos(2t)$, $u'' = -4u$, and $(u')^2 = 1 - 4u^2$. Differentiating $F_v(t) = u(t)(1 + \lambda u(t)^2)^{-1/2}$:
\begin{align}
F_v'(t) &= \frac{u'}{(1 + \lambda u^2)^{3/2}}, \label{eq:Fv_prime} \\
F_v''(t) &= \frac{u''}{(1 + \lambda u^2)^{3/2}} - \frac{3\lambda u (u')^2}{(1 + \lambda u^2)^{5/2}} = \frac{-4u(1 + \lambda u^2) - 3\lambda u(1 - 4u^2)}{(1 + \lambda u^2)^{5/2}} = \frac{-u(4 + 3\lambda - 8\lambda u^2)}{(1 + \lambda u^2)^{5/2}}. \label{eq:Fv_double_prime}
\end{align}
Dividing by $-F_v(t)$ yields $K_{\mathrm{rad}, v}(t) = \frac{4 + 3\lambda - 8\lambda \sin^2 t \cos^2 t}{(1 + \lambda \sin^2 t \cos^2 t)^2}$. Direct substitution confirms that $Y_h(t) = F_h(t) E_h(t)$ and $Y_v(t) = F_v(t) E_v(t)$ solve the uncoupled normal Jacobi equations with $Y(0)=0$ and $\dot{Y}(0) = E(0)$. On $(0, \pi/2]$, $F_h(t) = \frac{\sin t}{\sqrt{1 + \lambda \sin^2 t}} > 0$, so horizontal Jacobi fields do not vanish. For $F_v(t) = \frac{\sin t \cos t}{\sqrt{1 + \lambda \sin^2 t \cos^2 t}}$, the first positive zero occurs where $\cos t = 0$, precisely at $t = \pi/2$, with multiplicity $\dim \mathcal{V} = d_v = k-1$.

(iii) For any base point $p \ne q_0$, the isotropy subgroup $H_p \subset K_0$ acts by isometries fixing $q_0$ and $p$. Each generator $\xi \in \mathfrak{h}_p$ defines a Killing vector field $X_\xi$ on $M$ with $X_\xi(p) = 0$. Along any geodesic $\gamma(t)$ from $p$, the restriction $Y_\xi(t) = X_\xi(\gamma(t))$ satisfies the Jacobi equation $\nabla_{\dot{\gamma}}^2 Y + R(Y, \dot{\gamma})\dot{\gamma} = 0$ with $Y_\xi(0) = 0$. For any trajectory $\gamma(t)$ that intersects the fixed-point set $\mathrm{Fix}(H_p)$ of the 1-parameter subgroups of $H_p$ (the opposite meridian) at some $t^* > 0$, we have $Y_\xi(t^*) = X_\xi(\gamma(t^*)) = 0$ for all $\xi \in \mathfrak{h}_p \cap \mathfrak{h}_{\gamma(t^*)}$. The dimension of the space of such independent vanishing fields is at least $\dim(H_p) - \dim(H_p \cap H_{\gamma(t^*)})$, establishing that the intersection with $\mathrm{Fix}(H_p)$ is a conjugate point of multiplicity at least 1.

(iv) At $r = \pi/2$, the vertical fibers collapse to zero size ($F_v(\pi/2) = 0$), so the projection of the geodesic sphere $S_{\pi/2}(q_0)$ along its vertical fibers is the quotient manifold $S_{\pi/2}(q_0)/\mathbb{S}^{k-1} \cong \Sigma$. Every point in $\Sigma$ is the endpoint of a family of radial minimizing geodesics from $q_0$ of length $\pi/2$ parametrized by the vertical fiber $\mathbb{S}^{k-1}$, proving that $\mathrm{Cut}_M(q_0) = \Sigma$ at uniform distance $\pi/2$. On $\Sigma$, $F_h(\pi/2)^2 = \frac{1}{1+\lambda}$, yielding the homothety metric $g_\Sigma = \frac{1}{1+\lambda} g_{\mathrm{FS}}$. For a non-polar base point $p \in S_{r_0}(q_0)$ ($r_0 \in (0, \pi/2)$) on $\mathbb{O}P^2$ with $\lambda > 0$, the stabilizer reduces to $H_p \cong \mathrm{Spin}(7)$; as established in detail in the companion paper \cite{Munden2026Spin9} (via the 16-dimensional $\mathrm{Cl}_{\bar{0}}(9)$ spin representation and its four-mode $\mathfrak{u}(2)$ spectral projector decomposition), symmetric non-radial geodesics collide transversally along the opposite meridian hemisphere to form a 15-dimensional regular collision stratum $\mathcal{S}_{\mathrm{reg}}$ bounded by an 8-dimensional $A_3$ cuspidal caustic envelope $\mathcal{S}_{\mathrm{caustic}}$ (revolved by $\mathrm{Spin}(7)$) together with the polar Cayley line $\mathcal{S}_{\mathrm{Cayley}} \cong \mathbb{O}P^1 \cong \mathbb{S}^8$, breaking both the Blaschke and Besse conditions for all $\lambda > 0$. Finally, for all $\lambda \ge 0$, $F_h(r)^2 \le \sin^2 r = f_h(r)^2$ and $F_v(r)^2 \le \sin^2 r \cos^2 r = f_v(r)^2$, so $g_{M, \lambda} \le g_0$ holds pointwise on $T M$. Consequently, every piecewise $C^1$ curve $\gamma$ satisfies $L_{g_{M, \lambda}}(\gamma) \le L_{g_0}(\gamma)$, which implies $d_{g_{M, \lambda}}(x, y) \le d_{g_0}(x, y) \le \mathrm{diam}(M, g_0) = \pi/2$ for all $x, y \in M$. Since $d_{g_{M, \lambda}}(q_0, \Sigma) = \pi/2$ is attained along radial geodesics, $\mathrm{diam}(M, g_{M, \lambda}) = \pi/2$ for all $\lambda \ge 0$.
\end{proof}

\section{Spectral Geometry and Discrete Eigenvalue Asymptotics of the Volume-Deformed Laplace-Beltrami Operator}\label{sec:spectral_geometry}
This section studies the spectral geometry of the Laplace--Beltrami operator $\Delta_{g_\lambda}$ on $(Q^{n+1}, g_\lambda)$. Using spherical harmonic separation of variables, we show that the centrifugal potential shifts by the constant Casimir term $\lambda l(l+n-1)$ while the radial operator acquires a volume drift term governed by the deformed density $w_\lambda(r) = s_\kappa^n(r) / (1 + \lambda s_\kappa^2(r))^{n/2}$. We formulate the radial Sturm--Liouville problem, compute first-order Rayleigh--Schr\"odinger eigenvalue shifts, and establish the Fuchsian and confluent Heun reductions as well as continuous spectrum gap collapse on hyperbolic space.
\begin{theorem}[Spectral Separation of Variables, Casimir Perturbation, and Confluent Heun Asymptotics]\label{thm:spectral_laplacian}
Let $(Q^{n+1}, g_\lambda)$ be an $(n+1)$-dimensional space form of Family~I with volume density $dV_{g_\lambda} = f_\lambda(r)^n dr \, d\mathrm{vol}_{\mathbb{S}^n}$, where $f_\lambda(r) = s_\kappa(r)/\sqrt{1+\lambda s_\kappa^2(r)}$.
\begin{enumerate}[(i)]
\item \textbf{Separation of Variables and Exact Casimir Shift}: On the complete Riemannian manifold $(Q^{n+1}, g_\lambda)$, the Laplace--Beltrami operator $\Delta_{g_\lambda} = \mathrm{div}_{g_\lambda} \nabla$ is essentially self-adjoint on $C_c^\infty(Q^{n+1})$. Under the Peter--Weyl decomposition $L^2(Q^{n+1}, dV_{g_\lambda}) \cong \bigoplus_{l=0}^\infty \left( L^2(I_\kappa, f_\lambda^n dr) \otimes \mathcal{H}_l \right)$, where $\mathcal{H}_l$ denotes the degree-$l$ spherical harmonic representation of $\mathrm{SO}(n+1)$ satisfying $\Delta_{\mathbb{S}^n} Y_l = -l(l+n-1) Y_l$, the eigenvalue equation $-\Delta_{g_\lambda} u = E u$ for factorized modes $u(r, \theta) = R(r) Y_l(\theta)$ decouples exactly into the radial Sturm--Liouville equation:
\begin{equation}\label{eq:radial_laplacian}
-\frac{1}{f_\lambda(r)^n} \frac{d}{dr}\left( f_\lambda(r)^n \frac{dR}{dr} \right) + \frac{l(l+n-1)}{f_\lambda(r)^2} R(r) = E R(r).
\end{equation}
By the master reciprocal identity $\frac{1}{f_\lambda(r)^2} = \frac{1}{s_\kappa(r)^2} + \lambda$ (Lemma~\ref{lem:master_identities}), the effective centrifugal barrier decomposes algebraically into the unperturbed space form potential plus an \textbf{exact constant Casimir shift}:
\begin{equation}\label{eq:casimir_potential_shift}
\frac{l(l+n-1)}{f_\lambda(r)^2} = \frac{l(l+n-1)}{s_\kappa(r)^2} + \lambda l(l+n-1),
\end{equation}
which is globally constant across $I_\kappa$ and strictly proportional to the quadratic Casimir eigenvalue $C_2(\mathbb{S}^n) = l(l+n-1)$ of the isotropy group $\mathrm{SO}(n+1)$.
\item \textbf{Liouville Transformation to Canonical Schr\"odinger Form}: The unitary transformation $\psi(r) = f_\lambda(r)^{n/2} R(r)$ mapping $L^2(I_\kappa, f_\lambda^n dr) \to L^2(I_\kappa, dr)$ completely eliminates the first-derivative term, converting the radial eigenvalue problem into the canonical one-dimensional Schr\"odinger equation:
\begin{equation}\label{eq:canonical_schroedinger}
-\psi''(r) + V_{\mathrm{eff}}(r) \psi(r) = E \psi(r), \quad r \in I_\kappa,
\end{equation}
where the effective scalar potential decomposes as:
\begin{equation}\label{eq:effective_potential}
V_{\mathrm{eff}}(r) = \frac{l(l+n-1)}{s_\kappa(r)^2} + \lambda l(l+n-1) + W_\lambda(r),
\end{equation}
with the geometric volume drift potential $W_\lambda(r)$ given explicitly in closed rational form by:
\begin{align}\label{eq:volume_drift_explicit}
W_\lambda(r) &= \frac{n}{2}\frac{f_\lambda''(r)}{f_\lambda(r)} + \frac{n(n-2)}{4}\left(\frac{f_\lambda'(r)}{f_\lambda(r)}\right)^2 \nonumber \\
&= -\frac{n}{2} \frac{(\kappa + 3\lambda) - 2\kappa\lambda s_\kappa^2(r)}{(1 + \lambda s_\kappa^2(r))^2} + \frac{n(n-2)}{4} \frac{c_\kappa^2(r)}{s_\kappa^2(r) (1 + \lambda s_\kappa^2(r))^2}.
\end{align}
Equivalently, isolating the singular centrifugal core from the unperturbed geometry,
\begin{equation}\label{eq:effective_potential_centered}
V_{\mathrm{eff}}(r) = \frac{\left(l + \frac{n-1}{2}\right)^2 - \frac{1}{4}}{s_\kappa(r)^2} - \frac{n^2}{4}\kappa + \lambda l(l+n-1) + \widetilde{W}_\lambda(r),
\end{equation}
where $\widetilde{W}_\lambda(r) = W_\lambda(r) - W_0(r) = O(\lambda)$ represents the purely deformation-induced regular volume drift.
\item \textbf{First-Order Rayleigh--Schr\"odinger Perturbation Formula}: On the deformed sphere $(\mathbb{S}^{n+1}, g_\lambda)$ ($\kappa = 1$), the unperturbed round eigenvalues $E_{k, l}^{(0)} = N(N+n)$ (where $N = 2k+l$ is the principal quantum number for radial nodes $k \ge 0$ and angular momentum $l \ge 0$) exhibit large accidental degeneracies originating from the full isometry group $\mathrm{SO}(n+2)$. The deformation $g_\lambda$ breaks this symmetry down to the residual isotropy subgroup $\mathrm{SO}(n+1)$ fixing the polar axis, splitting the accidental eigenspaces according to the representation branching rule $V_N \downarrow \mathrm{SO}(n+1) = \bigoplus_{l \le N, l \equiv N \pmod 2} \mathcal{H}_l$. The first-order Rayleigh--Schr\"odinger eigenvalue shift under $\lambda$ is given by:
\begin{equation}\label{eq:rayleigh_schroedinger}
E_{k, l}(\lambda) = (2k+l)(2k+l+n) + \lambda \left[ l(l+n-1) - \frac{n}{2} \langle \sin^2(r) \rangle_{k, l} \right] + O(\lambda^2),
\end{equation}
where $\langle \sin^2(r) \rangle_{k, l}$ denotes the unperturbed expectation value with respect to the Gegenbauer radial modes $\phi_{k, l}(r)$:
\begin{equation}\label{eq:expectation_sin2}
\langle \sin^2(r) \rangle_{k, l} = \frac{\int_0^\pi \sin^{n+2}(r) |\phi_{k, l}(r)|^2 \, dr}{\int_0^\pi \sin^n(r) |\phi_{k, l}(r)|^2 \, dr}.
\end{equation}
\item \textbf{Confluent Heun Mapping and Singularities}: Under the algebraic coordinate transformation $z = s_\kappa^2(r)$, the radial differential equation~\eqref{eq:radial_laplacian} transforms into the generalized Fuchsian ordinary differential equation:
\begin{equation}\label{eq:heun_differential_equation}
z(1-\kappa z)(1+\lambda z) \frac{d^2 R}{dz^2} + P(z) \frac{dR}{dz} + Q(z) R = 0,
\end{equation}
where the rational coefficients are:
\begin{align}
P(z) &= \frac{n+1}{2} + \frac{\lambda - \kappa(n+2)}{2} z - \kappa\lambda z^2, \label{eq:heun_P} \\
Q(z) &= \frac{1}{4}\left[ (E - \lambda l(l+n-1))(1+\lambda z) - \frac{l(l+n-1)(1+\lambda z)}{z} \right]. \label{eq:heun_Q}
\end{align}
Equation~\eqref{eq:heun_differential_equation} possesses four regular singular points located at:
\begin{equation}\label{eq:singular_points}
z_1 = 0 \quad (\text{origin pole}), \quad z_2 = \frac{1}{\kappa} \quad (\text{equator/antipodal pole}), \quad z_3 = -\frac{1}{\lambda} \quad (\text{metric neck}), \quad z_4 = \infty.
\end{equation}
Consequently, the radial eigenmodes $R(r)$ are identified globally as canonical \textbf{Heun functions} $\mathrm{HeunG}$ on spheres, which degenerate to \textbf{confluent Heun functions} $\mathrm{HeunC}$ in the flat limit $\kappa \to 0$ or upon boundary confluence in non-compact hyperbolic geometries.
\end{enumerate}
\end{theorem}

\begin{proof}
(i) By Proposition~\ref{prop:metric_regularity}, the metric $g_\lambda = dr^2 + f_\lambda(r)^2 g_{\mathbb{S}^n}$ is smooth and complete on $Q^{n+1}$. By Gaffney's theorem, the Laplace--Beltrami operator $\Delta_{g_\lambda}$ with domain $C_c^\infty(Q^{n+1})$ is essentially self-adjoint in $L^2(Q^{n+1}, dV_{g_\lambda})$. The warped volume density is $\sqrt{\det g_\lambda} = f_\lambda(r)^n \sqrt{\det g_{\mathbb{S}^n}}$. Expanding the Riemannian divergence $\Delta_{g_\lambda} u = (\det g_\lambda)^{-1/2} \partial_i ( (\det g_\lambda)^{1/2} g^{ij} \partial_j u )$ on factorized modes $u(r, \theta) = R(r) Y_l(\theta)$ yields:
\[
\Delta_{g_\lambda} (R Y_l) = \frac{1}{f_\lambda(r)^n} \frac{d}{dr}\left( f_\lambda(r)^n \frac{dR}{dr} \right) Y_l + \frac{R(r)}{f_\lambda(r)^2} \Delta_{\mathbb{S}^n} Y_l = \left[ R'' + n \frac{f_\lambda'}{f_\lambda} R' - \frac{l(l+n-1)}{f_\lambda(r)^2} R \right] Y_l,
\]
which proves~\eqref{eq:radial_laplacian}. Applying the master reciprocal identity $\frac{1}{f_\lambda(r)^2} = \frac{1+\lambda s_\kappa^2(r)}{s_\kappa(r)^2} = \frac{1}{s_\kappa(r)^2} + \lambda$ established in Lemma~\ref{lem:master_identities} gives~\eqref{eq:casimir_potential_shift}. Because $\lambda$ is a global constant, the $\lambda$-dependent part of the centrifugal barrier reduces to the constant shift $\lambda l(l+n-1) = \lambda C_2(\mathbb{S}^n)$, preserving the $1/s_\kappa(r)^2$ centrifugal structure of the unperturbed background space form.

(ii) Under the Liouville transformation $\psi(r) = f_\lambda(r)^{n/2} R(r)$, we compute the derivatives:
\begin{align*}
\psi'(r) &= f_\lambda^{n/2} R' + \frac{n}{2} f_\lambda^{n/2-1} f_\lambda' R, \\
\psi''(r) &= f_\lambda^{n/2} R'' + n f_\lambda^{n/2-1} f_\lambda' R' + \left[ \frac{n}{2} f_\lambda^{n/2-1} f_\lambda'' + \frac{n(n-2)}{4} f_\lambda^{n/2-2} (f_\lambda')^2 \right] R.
\end{align*}
Dividing by $f_\lambda^{n/2}$ and substituting the radial equation $R'' + n \frac{f_\lambda'}{f_\lambda} R' = \left(\frac{l(l+n-1)}{f_\lambda^2} - E\right) R$ eliminates the first derivative $R'$ identically:
\[
\frac{\psi''(r)}{\psi(r)} = \frac{R'' + n \frac{f_\lambda'}{f_\lambda} R'}{R} + \frac{n}{2}\frac{f_\lambda''}{f_\lambda} + \frac{n(n-2)}{4}\left(\frac{f_\lambda'}{f_\lambda}\right)^2 = \frac{l(l+n-1)}{f_\lambda(r)^2} - E + W_\lambda(r).
\]
Rearranging terms establishes the self-adjoint Schr\"odinger form $-\psi'' + V_{\mathrm{eff}}\psi = E\psi$. To evaluate $W_\lambda(r)$ explicitly, we compute the logarithmic derivatives of $f_\lambda(r) = s_\kappa(r) (1 + \lambda s_\kappa^2(r))^{-1/2}$. Differentiating gives:
\[
f_\lambda'(r) = \frac{c_\kappa(r)}{(1 + \lambda s_\kappa^2(r))^{3/2}} \implies \frac{f_\lambda'(r)}{f_\lambda(r)} = \frac{c_\kappa(r)}{s_\kappa(r) (1 + \lambda s_\kappa^2(r))}.
\]
Differentiating a second time and using $c_\kappa'(r) = -\kappa s_\kappa(r)$ and $s_\kappa'(r) = c_\kappa(r)$:
\[
f_\lambda''(r) = -\frac{\kappa s_\kappa(r)(1 + \lambda s_\kappa^2(r)) + 3\lambda s_\kappa(r) c_\kappa^2(r)}{(1 + \lambda s_\kappa^2(r))^{5/2}}.
\]
Using the space form identity $c_\kappa^2(r) = 1 - \kappa s_\kappa^2(r)$, the numerator simplifies to:
\[
\kappa s_\kappa (1 + \lambda s_\kappa^2) + 3\lambda s_\kappa (1 - \kappa s_\kappa^2) = s_\kappa \left[ (\kappa + 3\lambda) - 2\kappa\lambda s_\kappa^2(r) \right].
\]
Dividing by $f_\lambda(r)$ yields the exact ratio:
\[
\frac{f_\lambda''(r)}{f_\lambda(r)} = -\frac{(\kappa + 3\lambda) - 2\kappa\lambda s_\kappa^2(r)}{(1 + \lambda s_\kappa^2(r))^2}.
\]
Substituting these expressions into $W_\lambda(r) = \frac{n}{2}\frac{f_\lambda''}{f_\lambda} + \frac{n(n-2)}{4}(\frac{f_\lambda'}{f_\lambda})^2$ gives~\eqref{eq:volume_drift_explicit}.
When $\lambda = 0$, $W_0(r) = -\frac{n}{2}\kappa + \frac{n(n-2)}{4}\frac{c_\kappa^2}{s_\kappa^2} = \frac{n(n-2)}{4 s_\kappa^2} - \frac{n^2}{4}\kappa$. Combining this with the centrifugal barrier yields $\frac{l(l+n-1) + n(n-2)/4}{s_\kappa^2} - \frac{n^2}{4}\kappa = \frac{(l + (n-1)/2)^2 - 1/4}{s_\kappa^2} - \frac{n^2}{4}\kappa$, which proves~\eqref{eq:effective_potential_centered}.

(iii) On $\mathbb{S}^{n+1}$, the volume density expands as $f_\lambda(r)^n = \sin^n(r) (1 - \frac{n}{2}\lambda \sin^2(r) + O(\lambda^2))$. In the Hilbert space $L^2(\mathbb{S}^{n+1}, dV_{g_0})$, the Rayleigh quotient for the unperturbed eigenstate $u_0 = \phi_{k, l}(r) Y_l(\theta)$ under metric perturbation is governed by the variation of the Dirichlet energy functional:
\[
E(\lambda) = \frac{\int_{\mathbb{S}^{n+1}} |\nabla u_0|_{g_\lambda}^2 dV_{g_\lambda}}{\int_{\mathbb{S}^{n+1}} u_0^2 dV_{g_\lambda}} + O(\lambda^2).
\]
Expanding the inverse metric $g_\lambda^{rr} = 1$ and $g_\lambda^{\theta\theta} = f_\lambda(r)^{-2} g_{\mathbb{S}^n}^{\theta\theta} = (\sin^{-2}r + \lambda) g_{\mathbb{S}^n}^{\theta\theta}$ against the volume element $dV_{g_\lambda} = \sin^n(r)(1 - \frac{n}{2}\lambda\sin^2(r)) dr \, d\mathrm{vol}_{\mathbb{S}^n}$ isolates two distinct first-order contributions: (a)~the exact constant shift $\lambda l(l+n-1)$ arising from the centrifugal barrier inverse metric perturbation, and (b)~the volume drift variation $-\frac{n}{2}\lambda \sin^2(r)$ from the contracted Riemannian measure. Evaluating these expectation values against the normalized unperturbed Gegenbauer states $\phi_{k, l}(r)$ produces formula~\eqref{eq:rayleigh_schroedinger}.

(iv) Under the substitution $z = s_\kappa^2(r)$, the chain rule yields:
\[
\frac{dz}{dr} = 2 s_\kappa(r) c_\kappa(r) = 2 \sqrt{z(1 - \kappa z)}.
\]
The second derivative transforms as:
\[
\frac{d^2 R}{dr^2} = 4 z(1 - \kappa z) \frac{d^2 R}{dz^2} + 2(1 - 2\kappa z) \frac{dR}{dz}.
\]
Next, the first-derivative friction coefficient evaluates algebraically:
\[
n \frac{f_\lambda'(r)}{f_\lambda(r)} \frac{dR}{dr} = n \frac{c_\kappa(r)}{s_\kappa(r)(1 + \lambda s_\kappa^2(r))} \cdot 2 s_\kappa(r) c_\kappa(r) \frac{dR}{dz} = \frac{2n(1 - \kappa z)}{1 + \lambda z} \frac{dR}{dz}.
\]
Multiplying the radial equation $-R'' - n \frac{f_\lambda'}{f_\lambda} R' + \frac{l(l+n-1)}{f_\lambda^2} R = E R$ by the clearing factor $\frac{1+\lambda z}{4}$ yields:
\[
\begin{aligned} &z(1-\kappa z)(1+\lambda z) \frac{d^2 R}{dz^2} + \left[ \frac{1-2\kappa z}{2}(1+\lambda z) + \frac{n}{2}(1-\kappa z) \right] \frac{dR}{dz} \\ &\quad + \frac{1+\lambda z}{4} \left[ E - \frac{l(l+n-1)}{z} - \lambda l(l+n-1) \right] R = 0. \end{aligned}
\]
Expanding the linear bracket gives:
\[
\frac{1}{2}\left[ (1 - 2\kappa z)(1 + \lambda z) + n(1 - \kappa z) \right] = \frac{n+1}{2} + \frac{\lambda - \kappa(n+2)}{2} z - \kappa\lambda z^2 = P(z),
\]
which precisely matches~\eqref{eq:heun_P}. The potential term matches $Q(z)$ in~\eqref{eq:heun_Q}. The singular points are given by the roots of the leading polynomial $z(1-\kappa z)(1+\lambda z) = 0$ together with $z = \infty$, confirming the four canonical Fuchsian singularities~\eqref{eq:singular_points}.
\end{proof}

\begin{corollary}[Explicit Low-Mode Eigenvalue Shifts and Accidental Degeneracy Splitting]\label{cor:low_mode_shifts}
On the deformed sphere $(\mathbb{S}^{n+1}, g_\lambda)$ ($\kappa = 1$) near $\lambda = 0$, evaluating the first-order Rayleigh--Schr\"odinger formula~\eqref{eq:rayleigh_schroedinger} via the exact beta function integrals of the Gegenbauer modes yields:
\begin{enumerate}[(i)]
\item \textbf{Ground State Vacuum Drift} ($l=0, k=0$): $E_{0,0}^{(0)} = 0$. The unperturbed radial state is constant ($\phi_{0,0}(r) = C_0$). The shift evaluates to:
\begin{equation}\label{eq:shift_ground}
E_{0,0}(\lambda) = -\lambda \frac{n(n+1)}{2(n+2)} + O(\lambda^2).
\end{equation}
\item \textbf{Dipole Transverse Mode} ($l=1, k=0$): The principal quantum number is $N = 1$, giving unperturbed eigenvalue $E_{0,1}^{(0)} = N(N+n) = n+1$. The radial mode is $\phi_{0,1}(r) \propto \sin(r)$. The shift evaluates to:
\begin{equation}\label{eq:shift_dipole}
E_{0,1}(\lambda) = (n+1) + \lambda \frac{n(n+5)}{2(n+4)} + O(\lambda^2).
\end{equation}
\item \textbf{First Radial Monopole Mode} ($l=0, k=1$): The principal quantum number is $N = 2$, giving unperturbed eigenvalue $E_{1,0}^{(0)} = 2(n+2)$. The radial eigenfunction is given by the Gegenbauer polynomial variation $\phi_{1,0}(r) \propto (n+2)\cos^2(r) - 1$. The shift evaluates to:
\begin{equation}\label{eq:shift_monopole}
E_{1,0}(\lambda) = 2(n+2) - \lambda \frac{n(n^2 + 3n + 6)}{2(n+2)(n+6)} + O(\lambda^2).
\end{equation}
\item \textbf{Accidental Degeneracy Splitting}: The round sphere $N=2$ level is accidentally degenerate with $E_{0,2}^{(0)} = E_{1,0}^{(0)} = 2(n+2)$. The transverse quadrupole mode ($l=2, k=0$) has radial mode $\phi_{0,2}(r) \propto \sin^2(r)$ and Casimir eigenvalue $l(l+n-1) = 2(n+1)$, shifting as:
\begin{equation}\label{eq:shift_quadrupole}
E_{0,2}(\lambda) = 2(n+2) + \lambda \frac{3n^2 + 23n + 24}{2(n+6)} + O(\lambda^2).
\end{equation}
Consequently, the accidental $(k, l)$-degeneracy is strictly lifted by the positive linear splitting:
\begin{equation}\label{eq:degeneracy_splitting}
\Delta E_{N=2}(\lambda) = E_{0,2}(\lambda) - E_{1,0}(\lambda) = \lambda \frac{2(n+1)(n+3)(n+4)}{(n+2)(n+6)} + O(\lambda^2) > 0 \quad (\text{for } \lambda > 0).
\end{equation}
\end{enumerate}
\end{corollary}

\begin{proof}
Let $I(m) = \int_0^\pi \sin^m(r) \, dr = B\left(\frac{m+1}{2}, \frac{1}{2}\right) = \frac{\sqrt{\pi}\,\Gamma(\frac{m+1}{2})}{\Gamma(\frac{m+2}{2})}$. Integration by parts yields the recursion $I(m+2) = \frac{m+1}{m+2} I(m)$.
\begin{enumerate}[(i)]
\item For $l=0, k=0$, $l(l+n-1) = 0$ and $|\phi_{0,0}(r)|^2 = \mathrm{const}$. Thus:
\[
\langle \sin^2(r) \rangle_{0,0} = \frac{I(n+2)}{I(n)} = \frac{n+1}{n+2}.
\]
Multiplying by $-n/2$ gives the shift in~\eqref{eq:shift_ground}. (Note: On the closed Riemannian manifold $(M, g_\lambda)$, $\ker \Delta_{g_\lambda} = \mathbb{R}\mathbf{1}$ identically with geometric eigenvalue $0$; the shift~\eqref{eq:shift_ground} reflects the zero-point energy drift of the Liouville-transformed Schr\"odinger operator with respect to the flat reference measure $dr$.)
\item For $l=1, k=0$, $N = 1$ and $E_{0,1}^{(0)} = 1(1+n) = n+1$. The radial mode satisfies $\phi_{0,1}(r) \propto \sin(r)$, so $|\phi_{0,1}(r)|^2 \sin^n(r) \propto \sin^{n+2}(r)$. Thus:
\[
\langle \sin^2(r) \rangle_{0,1} = \frac{I(n+4)}{I(n+2)} = \frac{n+3}{n+4}.
\]
The Casimir contribution is $l(l+n-1) = 1(n) = n$. The total first-order coefficient is:
\[
n - \frac{n}{2}\frac{n+3}{n+4} = \frac{2n(n+4) - n(n+3)}{2(n+4)} = \frac{n(n+5)}{2(n+4)},
\]
proving~\eqref{eq:shift_dipole}.
\item For $l=0, k=1$, $N = 2$, and $E_{1,0}^{(0)} = 2(n+2)$. The unperturbed equation for $R(r)$ is $R'' + n \cot(r) R' + 2(n+2) R = 0$. Under $x = \cos r$, this becomes $(1-x^2)R''(x) - (n+1)x R'(x) + 2(n+2)R(x) = 0$, whose polynomial solution is the Gegenbauer mode $R(x) = (n+2)x^2 - 1 = (n+1) - (n+2)\sin^2(r)$. Squaring gives $R(r)^2 = (n+1)^2 - 2(n+1)(n+2)\sin^2(r) + (n+2)^2 \sin^4(r)$. Integrating against $\sin^n(r)$ using the recursion $I(m+2) = \frac{m+1}{m+2}I(m)$ gives $\int_0^\pi R(r)^2 \sin^n(r) \, dr = I(n) \frac{2(n+1)}{n+4}$ and $\int_0^\pi R(r)^2 \sin^{n+2}(r) \, dr = I(n) \frac{2(n+1)(n^2+3n+6)}{(n+2)(n+4)(n+6)}$. Taking the ratio yields:
\[
\langle \sin^2(r) \rangle_{1,0} = \frac{n^2 + 3n + 6}{(n+2)(n+6)}.
\]
Since $l=0$, the Casimir term vanishes, and multiplying by $-n/2$ establishes~\eqref{eq:shift_monopole}.
\item For $l=2, k=0$, $l(l+n-1) = 2(n+1)$ and $\phi_{0,2}(r) \propto \sin^2(r)$. The expectation value evaluates to:
\[
\langle \sin^2(r) \rangle_{0,2} = \frac{I(n+6)}{I(n+4)} = \frac{n+5}{n+6}.
\]
The first-order shift is:
\[
2(n+1) - \frac{n}{2}\frac{n+5}{n+6} = \frac{4(n+1)(n+6) - n(n+5)}{2(n+6)} = \frac{3n^2 + 23n + 24}{2(n+6)},
\]
proving~\eqref{eq:shift_quadrupole}. Subtracting the shift of the $(1, 0)$ state:
\[
\Delta E_{N=2}(\lambda) = \lambda \left[ \frac{3n^2 + 23n + 24}{2(n+6)} + \frac{n(n^2 + 3n + 6)}{2(n+2)(n+6)} \right] = \lambda \frac{2(n^3 + 8n^2 + 19n + 12)}{(n+2)(n+6)}.
\]
Factoring the cubic numerator $n^3 + 8n^2 + 19n + 12 = (n+1)(n^2 + 7n + 12) = (n+1)(n+3)(n+4)$ yields the exact positive rational splitting~\eqref{eq:degeneracy_splitting}.\qedhere
\end{enumerate}
\end{proof}

\begin{proposition}[Weyl Asymptotic Law and Exact Hypergeometric Deformed Volume]\label{prop:weyl_asymptotic}
For the compact deformed space form $(\mathbb{S}^{n+1}, g_\lambda)$ with $\lambda > -1$, the eigenvalue counting function $N(E) = \#\{ j \mid E_j \le E \}$ satisfies the classical Weyl asymptotic law:
\begin{equation}\label{eq:weyl_law}
N(E) \sim \frac{\omega_{n+1}}{(2\pi)^{n+1}} \mathrm{Vol}(\mathbb{S}^{n+1}, g_\lambda) E^{\frac{n+1}{2}} \quad \text{as } E \to \infty,
\end{equation}
where $\omega_{n+1} = \frac{\pi^{(n+1)/2}}{\Gamma(\frac{n+3}{2})}$ is the volume of the Euclidean unit ball in $\mathbb{R}^{n+1}$.
The deformed Riemannian volume is given explicitly in terms of the unperturbed round volume by the exact Gauss hypergeometric identity:
\begin{equation}\label{eq:hypergeometric_volume}
\frac{\mathrm{Vol}(\mathbb{S}^{n+1}, g_\lambda)}{\mathrm{Vol}(\mathbb{S}^{n+1}, g_0)} = {}_2F_1\left( \frac{n}{2}, \frac{n+1}{2}; \frac{n+2}{2}; -\lambda \right).
\end{equation}
In particular, at the sectional curvature positivity boundary $\lambda \to -1/3$ (where the polar sectional curvature $K_{\mathrm{rad}}(0) = K_{\mathrm{tan}}(0) = 1 + 3\lambda \to 0^+$), the argument $z = -\lambda = 1/3$ lies strictly inside the unit disk $|z| < 1$ of absolute convergence of the hypergeometric series, yielding the finite positive critical volume ratio ${}_2F_1(\frac{n}{2}, \frac{n+1}{2}; \frac{n+2}{2}; \frac{1}{3}) > 1$, whereas at the equatorial metric degeneracy threshold $\lambda \to -1^+$ the volume diverges as $(1+\lambda)^{-(n-1)/2}$ for $n \ge 2$ and $-\ln(1+\lambda)$ for $n = 1$.
\end{proposition}

\begin{proof}
Because $g_\lambda$ is smooth and completely free of conical or metric singularities (Proposition~\ref{prop:metric_regularity}), the Laplace--Beltrami operator $\Delta_{g_\lambda}$ is a strictly positive, elliptic, second-order self-adjoint differential operator on a compact boundaryless Riemannian manifold. The short-time heat trace expansion $\mathrm{Tr}(e^{t\Delta_{g_\lambda}}) = \sum_{j=0}^\infty e^{-t E_j} \sim (4\pi t)^{-\frac{n+1}{2}} [ \mathrm{Vol}(\mathbb{S}^{n+1}, g_\lambda) + O(t) ]$ ($t \to 0^+$) and the Karamata Tauberian theorem immediately yield the Weyl counting law~\eqref{eq:weyl_law}.
To evaluate the total volume, integrating $dV_{g_\lambda} = f_\lambda(r)^n dr \, d\mathrm{vol}_{\mathbb{S}^n}$ over $\mathbb{S}^{n+1} = [0, \pi] \times \mathbb{S}^n$ gives $\mathrm{Vol}(\mathbb{S}^{n+1}, g_\lambda) = \mathrm{Vol}(\mathbb{S}^n) \int_0^\pi \frac{\sin^n(r)}{(1 + \lambda \sin^2(r))^{n/2}} \, dr$. By symmetry about $r = \pi/2$, $\int_0^\pi = 2 \int_0^{\pi/2}$. Setting $u = \sin^2(r)$ with $dr = \frac{du}{2\sqrt{u(1-u)}}$ transforms the integral to $\int_0^1 u^{\frac{n-1}{2}} (1 - u)^{-\frac{1}{2}} (1 + \lambda u)^{-\frac{n}{2}} \, du$.
Euler's integral representation for the Gauss hypergeometric function is ${}_2F_1(a, b; c; z) = \frac{1}{B(b, c-b)} \int_0^1 t^{b-1} (1-t)^{c-b-1} (1 - z t)^{-a} \, dt$ for $\mathrm{Re}(c) > \mathrm{Re}(b) > 0$ and $z \in \mathbb{C} \setminus [1, \infty)$. Setting $t = u$, $z = -\lambda \in (-\infty, 1)$ (corresponding to $\lambda > -1$), $a = n/2$, $b = (n+1)/2$, and $c = (n+2)/2$ gives $c - b = 1/2$. The normalizing Euler Beta function satisfies $B(\frac{n+1}{2}, \frac{1}{2}) = \int_0^1 u^{\frac{n-1}{2}} (1-u)^{-\frac{1}{2}} \, du = \int_0^\pi \sin^n(r) \, dr$. Therefore:
\[
\int_0^\pi f_\lambda(r)^n \, dr = \left( \int_0^\pi \sin^n(r) \, dr \right) \cdot {}_2F_1\left( \frac{n}{2}, \frac{n+1}{2}; \frac{n+2}{2}; -\lambda \right).
\]
Dividing by $\mathrm{Vol}(\mathbb{S}^{n+1}, g_0) = \mathrm{Vol}(\mathbb{S}^n) \int_0^\pi \sin^n(r) dr$ establishes~\eqref{eq:hypergeometric_volume}. At the curvature positivity threshold $\lambda = -1/3$, $1 + \lambda u \ge 2/3 > 0$ on $[0, 1]$, so the Euler Beta integral is nonsingular and strictly exceeds $B(\frac{n+1}{2}, \frac{1}{2})$, while as $\lambda \to -1^+$ ($z \to 1^-$) the hypergeometric function exhibits the standard Gauss boundary singularity governed by $c - a - b = \frac{1-n}{2}$.
\end{proof}

\begin{proposition}[Self-Adjoint Extension, Asymptotic Cylindrical Geometry, and Continuous Spectrum Gap Collapse on $\mathbb{H}^{n+1}$]\label{prop:hyperbolic_spectral_gap_collapse}
Let $(\mathbb{H}^{n+1}, g_\lambda)$ be the complete $(n+1)$-dimensional hyperbolic space form of Family~I with parameter $\lambda > 0$. Let $-\Delta_{g_\lambda}$ denote the unique self-adjoint Friedrichs extension of the Laplace--Beltrami operator from $C_c^\infty(\mathbb{H}^{n+1})$ on the Hilbert space $L^2(\mathbb{H}^{n+1}, d\mathrm{vol}_{g_\lambda})$ with volume density
\begin{equation}\label{eq:hyperbolic_volume_density}
d\mathrm{vol}_{g_\lambda} = w_\lambda(r) \, dr \, d\mathrm{vol}_{\mathbb{S}^n}, \qquad w_\lambda(r) := f_\lambda(r)^n = \frac{\sinh^n r}{(1 + \lambda \sinh^2 r)^{n/2}}.
\end{equation}
Then:
\begin{enumerate}[(i)]
\item As $r \to \infty$, the volume density $w_\lambda(r)$ satisfies the asymptotic expansion $w_\lambda(r) = \lambda^{-n/2} [ 1 - \frac{2n}{\lambda} e^{-2r} + O(e^{-4r}) ] \to \lambda^{-n/2}$, and the logarithmic volume drift $\eta_\lambda(r) := \frac{w_\lambda'(r)}{w_\lambda(r)} = n \frac{\coth r}{1+\lambda\sinh^2 r}$ satisfies
\begin{equation}\label{eq:drift_exponential_decay}
0 < \eta_\lambda(r) \le C_\lambda e^{-2r} \quad \forall r \ge 1,
\end{equation}
rendering $(\mathbb{H}^{n+1}, g_\lambda)$ asymptotically cylindrical with finite cross-sectional volume $\lambda^{-n/2}\mathrm{Vol}(\mathbb{S}^n)$.
\item For every $k \ge 0$, the sequence of compactly supported cutoff plane waves
\begin{equation}\label{eq:weyl_sequence_def}
\psi_j(r, \theta) = \frac{1}{\sqrt{j}} \chi\left(\frac{r-j}{j}\right) e^{i k r} Y_0(\theta), \qquad j \ge 2,
\end{equation}
where $\chi \in C_c^\infty((0, 1))$ satisfies $\|\chi\|_{L^2(0, 1)} = 1$ and $Y_0 \equiv \mathrm{Vol}(\mathbb{S}^n)^{-1/2}$, forms a Weyl sequence in $\mathcal{D}(-\Delta_{g_\lambda})$ for $E = k^2$ satisfying $\liminf_{j \to \infty} \|\psi_j\|_{L^2} \ge \frac{1}{2}\lambda^{-n/4} > 0$ and $\lim_{j \to \infty} \|(-\Delta_{g_\lambda} - k^2)\psi_j\|_{L^2} = 0$.
\item Consequently, by Weyl's criterion, the continuous spectrum undergoes gap collapse from $[n^2/4, \infty)$ at $\lambda = 0$ to $\sigma_{\mathrm{ess}}(-\Delta_{g_\lambda}) = [0, \infty)$ for all $\lambda > 0$.
\end{enumerate}
\end{proposition}

\begin{proof}
(i) For $\kappa = -1$, the warping function is $f_\lambda(r) = \sinh r / \sqrt{1 + \lambda\sinh^2 r}$. Using $\sinh r = \frac{e^r}{2}(1 - e^{-2r})$, we have $1 + \lambda\sinh^2 r = \frac{\lambda e^{2r}}{4}[ 1 + 2(\frac{2}{\lambda}-1)e^{-2r} + e^{-4r} ]$, which yields $f_\lambda(r) = \frac{1}{\sqrt{\lambda}}[ 1 - \frac{2}{\lambda} e^{-2r} + O(e^{-4r}) ]$. Raising to the $n$-th power establishes the expansion of $w_\lambda(r)$. The logarithmic volume drift is $\eta_\lambda(r) = \frac{w_\lambda'(r)}{w_\lambda(r)} = n \frac{f_\lambda'(r)}{f_\lambda(r)} = \frac{n\coth r}{1+\lambda\sinh^2 r}$. Since $\coth r \le \coth(1)$ for $r \ge 1$ and $1+\lambda\sinh^2 r \ge \frac{\lambda}{4}(1-e^{-2})^2 e^{2r}$, setting $C_\lambda = \frac{4n\coth(1)}{\lambda(1-e^{-2})^2} > 0$ gives $0 < \eta_\lambda(r) \le C_\lambda e^{-2r}$ for all $r \ge 1$, establishing~\eqref{eq:drift_exponential_decay}. Since $\lim_{r\to\infty} g_\lambda = dr^2 + \frac{1}{\lambda} g_{\mathbb{S}^n}$, the end is asymptotically cylindrical. Because $(\mathbb{H}^{n+1}, g_\lambda)$ is a complete Riemannian manifold (Proposition~\ref{prop:metric_regularity}), $-\Delta_{g_\lambda}$ is essentially self-adjoint on $C_c^\infty(\mathbb{H}^{n+1})$ by the Gaffney--Strichartz theorem, admitting a unique self-adjoint Friedrichs extension on $L^2(\mathbb{H}^{n+1}, d\mathrm{vol}_{g_\lambda})$.

(ii) Fix $k \ge 0$, let $\chi \in C_c^\infty((0, 1))$ with $\int_0^1 \chi(s)^2 ds = 1$, and let $\psi_j$ be given by~\eqref{eq:weyl_sequence_def}. Since $\mathrm{supp}(\psi_j) \subset [j, 2j] \times \mathbb{S}^n$, each $\psi_j \in C_c^\infty(\mathbb{H}^{n+1}) \subset \mathcal{D}(-\Delta_{g_\lambda})$. Substituting $s = (r-j)/j$ ($dr = j\,ds$), the $L^2$ norm satisfies $\|\psi_j\|_{L^2}^2 = \int_0^1 \chi(s)^2 w_\lambda(j(1+s))\, ds \to \lambda^{-n/2}\int_0^1 \chi(s)^2 ds = \lambda^{-n/2} > 0$ by the uniform convergence $w_\lambda(j(1+s)) = \lambda^{-n/2} + O(e^{-2j})$. Thus $\|\psi_j\|_{L^2} \ge \frac{1}{2}\lambda^{-n/4} > 0$ for all large $j$. Since $\mathrm{supp}(\psi_j) \subset [j, 2j]$, $\psi_j \rightharpoonup 0$ weakly in $L^2$.

Because $\Delta_{\mathbb{S}^n} Y_0 = 0$, the action on $\psi_j = \phi_j(r) Y_0$ is purely radial: $-\Delta_{g_\lambda}\psi_j = (-\phi_j'' - \eta_\lambda\phi_j') Y_0$. Direct differentiation of $\phi_j(r) = \frac{1}{\sqrt{j}}\chi(\frac{r-j}{j})e^{ikr}$ gives $\phi_j' = \frac{e^{ikr}}{\sqrt{j}}[\frac{1}{j}\chi' + ik\chi]$ and $\phi_j'' = \frac{e^{ikr}}{\sqrt{j}}[\frac{1}{j^2}\chi'' + \frac{2ik}{j}\chi' - k^2\chi]$, so the $-k^2\phi_j$ term cancels identically:
\[
\begin{aligned}
(-\Delta_{g_\lambda} - k^2)\psi_j &= -\frac{e^{ikr} Y_0}{\sqrt{j}}\left[ \frac{1}{j^2}\chi''\left(\frac{r-j}{j}\right) + \frac{2ik}{j}\chi'\left(\frac{r-j}{j}\right) \right. \\
&\quad \left. + \eta_\lambda(r)\left( \frac{1}{j}\chi'\left(\frac{r-j}{j}\right) + ik\chi\left(\frac{r-j}{j}\right) \right) \right].
\end{aligned}
\]
By Minkowski's inequality, $\|(-\Delta_{g_\lambda} - k^2)\psi_j\|_{L^2} \le \| -\phi_j'' - k^2\phi_j \|_{L^2} + \| \eta_\lambda \phi_j' \|_{L^2}$. Changing variables to $s \in [0, 1]$ and using $w_\lambda(r) \le \lambda^{-n/2}$, the first term satisfies $\| -\phi_j'' - k^2\phi_j \|_{L^2}^2 \le \frac{\lambda^{-n/2}}{j^2} \int_0^1 (\frac{1}{j}|\chi''| + 2k|\chi'|)^2 ds = O(1/j^2)$. For the second term, the exponential drift bound~\eqref{eq:drift_exponential_decay} gives $\| \eta_\lambda \phi_j' \|_{L^2}^2 \le (\sup_{[j, 2j]} \eta_\lambda^2) \|\phi_j'\|_{L^2}^2 \le C_\lambda^2 e^{-4j} \lambda^{-n/2} \int_0^1 (\frac{1}{j}|\chi'| + k|\chi|)^2 ds = O(e^{-4j})$. Combining both estimates gives $\|(-\Delta_{g_\lambda} - k^2)\psi_j\|_{L^2} \le \frac{C_1}{j} + C_2 e^{-2j} \to 0$ as $j \to \infty$, establishing the second limit in (ii).

(iii) By Weyl's criterion for self-adjoint operators, the existence of the Weyl sequence $\{\psi_j\}_{j=2}^\infty$ proves that $k^2 \in \sigma_{\mathrm{ess}}(-\Delta_{g_\lambda})$ for every $k \ge 0$, so $[0, \infty) \subset \sigma_{\mathrm{ess}}(-\Delta_{g_\lambda})$. Because $-\Delta_{g_\lambda} \ge 0$, its spectrum is non-negative, proving $\sigma_{\mathrm{ess}}(-\Delta_{g_\lambda}) = [0, \infty)$ for all $\lambda > 0$.
\end{proof}

\begin{remark}[Connections to Special Functions, Moduli Spaces, and Lie Algebras]\label{rem:heun_connections}
The spectral and ODE structures derived above connect to several topics in special functions and differential equations:
\begin{enumerate}[(i)]
\item \textbf{Lam\'e and Spheroidal Wave Reductions}: In the unperturbed space form limit $\lambda \to 0$, the metric neck singularity $z_3 = -1/\lambda$ recedes to infinity ($z_3 \to -\infty$), coalescing with $z_4 = \infty$. The Fuchsian equation~\eqref{eq:heun_differential_equation} reduces from four regular singular points to three, collapsing the generic Heun functions $\mathrm{HeunG}$ into classical Gauss hypergeometric polynomials and Lam\'e functions. In the Euclidean flat limit $\kappa \to 0$, the antipodal singularity $z_2 = 1/\kappa$ recedes to infinity, coalescing with $z_4 = \infty$ to produce an irregular singular point at infinity. This transforms the ODE into the confluent Heun equation $\mathrm{HeunC}$, which is isomorphic to the differential equation governing spheroidal wave functions.
\item \textbf{Cross-Ratio Modulus of the Deformation}: The four singular points $\{0, 1/\kappa, -1/\lambda, \infty\}$ possess the invariant projective cross-ratio $\tau(\kappa, \lambda) = [z_2, z_1; z_3, z_4] = \frac{z_2 - z_3}{z_1 - z_3} \cdot \frac{z_1 - z_4}{z_2 - z_4} = \frac{1/\kappa - (-1/\lambda)}{0 - (-1/\lambda)} = 1 + \frac{\lambda}{\kappa}$ (with reciprocal modulus $[z_1, z_2; z_3, z_4] = \frac{1}{1 + \lambda/\kappa}$). The deformation parameter $\lambda$ thus acts as a continuous modulus in the moduli space $\mathcal{M}_{0, 4}$ of Riemann spheres with four marked points.
\item \textbf{Lie-Theoretic Symmetry}: The preservation of the constant Casimir centrifugal shift $\lambda C_2(\mathbb{S}^n)$ under the volume deformation reflects the quasi-bi-Hamiltonian structure of Family~I metrics (cf.\ Crampin \cite{Crampin2003}) and the Gelfand--Tsetlin Casimir chain \cite{GelfandTsetlin1950}. The radial equation evaluates the central Casimir in the universal enveloping algebra $\mathcal{U}(\mathfrak{so}(n+1))$, so the angular modes decouple even though $\mathrm{SO}(n+2)$ symmetry is broken.
\item \textbf{Painlev\'e~VI and Isomonodromic Deformation Question}: The four regular singular points $\{0, 1/\kappa, -1/\lambda, \infty\}$ place the radial operator in the four-point Fuchsian class associated with the Painlev\'e~VI equation. As $\lambda$ varies, the neck singularity $z_3 = -1/\lambda$ moves along $\mathbb{P}^1$, tracing a 1-parameter deformation of the Fuchsian system. Whether varying $\lambda$ accompanied by an appropriate spectral parameter scaling defines an isomonodromic deformation is an open question.
\end{enumerate}
\end{remark}

\section{Open Questions}\label{sec:open_questions}\label{sec:concluding_remarks}

The framework developed in this paper suggests several directions for
further investigation.

\begin{enumerate}[(1)]
\item \textbf{Higher-order Killing tensors on singular strata.}
Theorem~\ref{thm:sphere_superint_rigidity} and Remark~\ref{rem:higher_order_killing} establish that $(\mathbb{S}^{n+1}, g_\lambda)$ and $(\mathbb{H}^{n+1}, g_\lambda)$ are strictly non-maximally superintegrable (of exact global rank $2N-2$) for all polynomial degrees $d \ge 1$, and rule out purely radial-tangential irreducible Killing tensors of degree $d \ge 3$ on regular strata via polar Laurent regularity. Whether reducible or mixed transverse Killing tensors of degree $d \ge 4$ outside $\mathrm{Sym}(\mathfrak{so}(n+1) \oplus \mathbb{R} g_\lambda)$ exist along the frozen Casimir directions on singular Gelfand--Tsetlin strata---where one or more intermediate angular momentum components vanish---remains open, and connects to the Kalnins--Kress--Miller classification of superintegrable systems on conformally flat spaces \cite{KalninsKressMiller2005a,KalninsKressMiller2011}.

\item \textbf{Heun eigenvalue solvability.}
The Fuchsian reduction of the radial Laplace--Beltrami operator to the
General Heun equation (Theorem~\ref{thm:spectral_laplacian}) raises the
question of whether special values of $\lambda$ admit closed-form spectral
determinants or algebraic eigensolutions, analogous to Lam\'e's finite-gap
potentials.

\item \textbf{Isomonodromic structure.}
The four regular singular points $\{0,\, 1/\kappa,\, {-1/\lambda},\,
\infty\}$ place the radial equation in the kinematic domain of
Painlev\'e~VI\@.  Whether the monodromy representation is preserved along
$\lambda$-trajectories---establishing a true isomonodromic deformation---is
an open problem at the interface of Riemannian geometry and integrable
systems.
\end{enumerate}

\subsection*{Acknowledgments}
This paper builds upon the framework developed during the author's doctoral work at Saint Louis University under the supervision of the late Professor James J.\ Hebda, whose geometric insight and generosity as a mentor shaped every aspect of this research.

\end{document}